\documentclass[11pt]{article}

\usepackage[english]{babel}
\usepackage[T1]{fontenc}
\usepackage[utf8]{inputenc}
\usepackage{amsbsy}
\usepackage{amsmath}
\usepackage{amssymb}
\usepackage{amsfonts}
\usepackage{amsthm}
\usepackage{bm}
\usepackage{graphicx}
\usepackage[active]{srcltx}
\usepackage{upgreek}
\usepackage{xcolor}
\usepackage[all]{xy}
\usepackage{dsfont}

\usepackage{tikz}

\newcommand{\C}{\mathbb{C}}

\newcommand{\N}{\mathbb{N}}

\newcommand{\R}{\mathbb{R}}

\newcommand{\Z}{\mathbb{Z}}

\newcommand{\gd}{\mathfrak{d}}
\renewcommand{\ge}{\mathfrak{e}}

\newcommand{\gh}{\mathfrak{h}}

\newcommand{\gu}{\mathfrak{u}}

\newcommand{\grad}{\nabla}

\newcommand{\ctilde}{\widetilde{c}}

\newcommand{\varepsilontilde}{\widetilde{\varepsilon}}
\newcommand{\etatilde}{\widetilde{\eta}}
\newcommand{\etilde}{\widetilde{e}}
\newcommand{\etildeeta}{\widetilde{e}_\eta}
\newcommand{\etildev}{\widetilde{e}_v}
\newcommand{\ftilde}{\widetilde{f}}
\newcommand{\ftildeeta}{\widetilde{f}_\eta}
\newcommand{\ftildev}{\widetilde{f}_v}
\newcommand{\htilde}{\widetilde{h}}
\newcommand{\htildeeta}{\widetilde{h}_\eta}
\newcommand{\htildev}{\widetilde{h}_v}

\newcommand{\energyset}{\mathcal{X}(\R)}

\newcommand{\energysethydro}{\mathcal{X}_{hy}(\R)}
\newcommand{\energysethydrok}{\mathcal{X}_{hy}^l(\R)}
\newcommand{\Nenergyset}{\mathcal{NX}(\R)}
\newcommand{\Nenergysethydro}{\mathcal{NX}_{hy}(\R)}

\newcommand{\Nenergysethydrokk}[1]{\mathcal{NX}_{hy}^{#1}(\R)}

\newcommand{\grandOde}[1]{\mathcal{O}\left( #1\right)}

\newcommand{\ntend}{\underset{n\rightarrow +\ii}{\longrightarrow}}

\newcommand{\ttendii}{\underset{t\rightarrow +\ii}{\longrightarrow}}

\newcommand{\ntendf}{\underset{n\rightarrow +\ii}{\rightharpoonup}}

\newcommand{\ntendfd}{\overset{d}{\underset{n\rightarrow +\ii}{\rightharpoonup}}}

\newcommand{\ntendfLdeux}{\overset{L^2}{\underset{n\rightarrow +\ii}{\rightharpoonup}}}

\newcommand{\tiitendf}{\underset{t\rightarrow +\ii}{\rightharpoonup}}

\newcommand{\tiitendfd}{\overset{d}{\underset{t\rightarrow +\ii}{\rightharpoonup}}}

\newcommand{\ntendfX}{\overset{\mathcal{X}_{hy}}{\underset{n\rightarrow +\ii}{\rightharpoonup}}}

\newcommand{\tiitendfX}{\overset{\mathcal{X}_{hy}}{\underset{t\rightarrow +\ii}{\rightharpoonup}}}

\newcommand{\ii}{\infty}
\newcommand{\dom}{\mathcal{D}om}

\newcommand{\spess}{\sigma_{ess}}

\newcommand{\loc}{\mathrm{loc}}

\newcommand{\normX}[1]{\left\Vert #1\right\Vert_{\mathcal{X}_{hy}}}

\newcommand{\normXk}[1]{\left\Vert #1\right\Vert_{\mathcal{X}_{hy}^l}}
\newcommand{\normxxk}[2]{\left\Vert #1\right\Vert_{\mathcal{X}_{hy}^{#2}}}

\newcommand{\normLii}[1]{\left\Vert #1\right\Vert_{L^\ii}}
\newcommand{\normLdeux}[1]{\left\Vert #1\right\Vert_{L^2}}
\newcommand{\normLun}[1]{\left\Vert #1\right\Vert_{L^1}}

\newcommand{\normHun}[1]{\left\Vert #1\right\Vert_{H^1}}
\newcommand{\psLdeux}[2]{\left\langle #1,#2 \right\rangle_{L^2}}
\newcommand{\psLdeuxLdeux}[2]{\left\langle #1,#2 \right\rangle_{L^2\times L^2}}
\newcommand{\vvvert}[1]{\big\vert\kern-0.25ex\big\vert\kern-0.25ex\big\vert #1\big\vert\kern-0.25ex\big\vert\kern-0.25ex\big\vert}
\newcommand{\normLiiunun}[1]{\left\Vert #1\right\Vert_{L^\ii([-1,1])}}
\newcommand{\normXboule}[2]{\left \Vert #1\right\Vert_{\mathcal{X}_{hy} (\mathcal{B}_{#2})}}

\newcommand{\spanned}[1]{\mathrm{Span}(#1)}
\newcommand{\normLdeuxLdeux}[1]{\left\Vert #1\right\Vert_{L^2\times L^2}}

\newcommand{\normR}[1]{\left| #1\right|}
\newcommand{\normXpoidsx}[3]{\left \Vert #1\right\Vert_{\mathcal{X}_{hy}^{#3} \big(|x|^{\tiny #2}dx\big)}}
\newcommand{\normXpoidsunplusx}[3]{\left \Vert #1\right\Vert_{\mathcal{X}_{hy}^{#3} \big((1+|x|^{ #2})dx\big)}}

\newcommand{\normYkk}[2]{\left \Vert #1\right\Vert_{\mathcal{Y}^{#2}}}

\newcommand{\normYpoidsx}[3]{\left \Vert #1\right\Vert_{\mathcal{Y}^{#3} \big(|x|^{\tiny #2}dx\big)}}

\newcommand{\normWkp}[3]{\left \Vert #1\right\Vert_{W^{#2,#3}}}

\newcommand{\normWkcarre}[3]{\left \Vert #1\right\Vert_{W^{#2,#3}\times W^{#2,#3}}}

\newtheorem{thm}{Theorem}[section]

\newtheorem{claim}[thm]{Claim}
\newtheorem{cor}[thm]{Corollary}
\newtheorem{lem}[thm]{Lemma}
\newtheorem{prop}[thm]{Proposition}

\newtheorem*{thm*}{Theorem}
\newtheorem*{merci}{Acknowledgments}
\newtheorem{rem}[thm]{Remark}

\theoremstyle{remark}

\title{Asymptotic stability of travelling waves for general nonlinear Schrödinger equations with non-zero condition at infinity}

\begin{document}

\date{}
\author{
\renewcommand{\thefootnote}{\arabic{footnote}}
Jordan Berthoumieu\footnotemark[1]}
\footnotetext[1]{CY Cergy Paris Universit\'e, Laboratoire Analyse, G\'eom\'etrie, Mod\'elisation, F-95302 Cergy-Pontoise, France. Website: https://berthoumieujordan.wordpress.com \ E-mail: {\tt jordan.berthoumieu@cyu.fr}}
\maketitle

\begin{abstract}
In previous works~\cite{Bert1,Bert2}, existence and uniqueness of travelling waves for the nonlinear Schrödinger equations have been shown for speeds close to the speed of sound. Furthermore, it has been proved that a chain of dark solitons of well-ordered speeds near the sound speed, taken initially apart from each other, is orbitally stable. In this article, we complete this study by proving the asymptotic stability of these travelling waves, namely that a solution initially close to a travelling wave eventually converges towards a travelling wave of close speed. This article generalizes the stability results for the Gross-Pitaevskii equation addressed by F. Béthuel, P. Gravejat and D. Smets in~\cite{BetGrSm2} and relies on the methods used in the latter article and first introduced by Y. Martel and F. Merle in~\cite{MartMer7}.  
\end{abstract}

\section{Introduction}\label{section: intro}

The equation of main interest reads
\begin{equation}\label{NLS}\tag{$NLS$}
    i\partial_t \Psi +\partial_x^2 \Psi +\Psi f(|\Psi|^2)=0 \text{ for }(t,x)\in\R\times\R.
\end{equation}

This equation is a generalization of the well-known Gross-Pitaevskii equation in the sense that we impose that $f(1)=0$ and the natural condition at infinity is then
\begin{equation*}\label{nonvanishing condition à l'infini}
    |\Psi(t,x)|\underset{|x|\rightarrow +\ii}{\longrightarrow}1.
\end{equation*}

The Gross-Pitaevskii case corresponds to the nonlinearity $f(s)=1-s$. It is a relevant model in condensed matter physics and nonlinear optics. In particular, it arises in the context of the Bose-Einstein condensation or of superfluidity (see~\cite{AbiHuMeNoPhTu,Coste1}). In the context of nonlinear optics, the non-vanishing condition at infinity expresses the presence of a nonzero background. Besides in the latter context, experiments highlight a notch in the density distribution. This density wave propagates in a bright background, which explains the denomination of "dark" soliton. We refer to~\cite{Bestesipa} for a significant example and more details on these features. Our study relates to defocusing nonlinearities, namely functions satisfying 
\begin{equation*}\label{f'(1) <0}
    f'(1)<0.
\end{equation*}

The Hamiltonian structure of~\eqref{NLS} is given by the generalized Ginzburg-Landau energy,
\begin{equation*}\label{définition de l'energie de Ginzburg Landau}
    E(\Psi):=E_k(\Psi)+E_p(\Psi):=\dfrac{1}{2}\int_\R |\partial_x \Psi|^2+ \dfrac{1}{2}\int_\R F(|\Psi|^2),
\end{equation*}

with
\begin{equation*}\label{definition de F}
    F(\rho):=\int_\rho^1 f(r)dr.
\end{equation*}

In~\cite{Chiron7}, many different behaviors have been highlighted according to the characteristics of $f$. Then in~\cite{Bert1}, special examples are taken to investigate more precisely the existence and orbital stability of the special solutions that are analyzed in this article.

In the sequel, we shall rely on the hydrodynamical formulation which gives a more convenient way to study the Hamiltonian structure of the equation. This framework is related to new variables given by $Q=(\eta,v):=(1-|\Psi|^2,-\partial_x \varphi)$ for non-vanishing functions that we lift as $\Psi=|\Psi| e^{i\varphi}$. Plugging the new variables into~\eqref{NLS}, we obtain the hydrodynamical form of~\eqref{NLS}, that is
\begin{equation}\label{NLShydro}\tag{$NLS_{hy}$}
    \left\{
\begin{array}{l}
    \partial_t\eta =-2\partial_x \big(v(1-\eta)\big), \\
    \partial_t v =-\partial_x \bigg( f(1-\eta)-v^2-\dfrac{\partial_x^2 \eta}{2(1-\eta)}-\dfrac{(\partial_x\eta)^2}{4(1-\eta)^2} \bigg). \\
\end{array}
\right.
\end{equation}

The main difference with Schrödinger equations with solutions vanishing at infinity is that the trivial solution becomes $\Psi\equiv 1$ instead of $0$. When we linearize around the hydrodynamic trivial solution $Q=(\eta,v)\equiv (0,0)$, we obtain, in the long wave approximation, the free wave equation with the sound speed 
\begin{equation}\label{relation entre vitesse c_s et f'(1)}
c_s = \sqrt{-2 f'(1)}.
\end{equation}

The linearized equation indeed differs from the case of null condition at infinity, in the sense that the resulting dispersion relation is not the same: the plane wave solutions of this linearized system satisfy the dispersion relation
\begin{equation*}
    \omega(\xi)=\pm\sqrt{\xi^4+c_s^2 \xi^2}.
\end{equation*}

For the sake of clarity, we now define the momentum\footnote{We refer to Subsection~\ref{section: Liouville property around a travelling wave} for rigorous definitions and regularity statements regarding the momentum $p$, and the energy $E$.}, in the hydrodynamical variables, as
\begin{equation}\label{expression du moment}
    p(Q)=p(\eta,v):=\int_\R\pi(\eta,v):=\dfrac{1}{2}\int_\R \eta v .
\end{equation}

The main result of stability is stated in the classical variables, and since we deal with infinite dimensional dynamical systems, we need to find a suitable functional setting. We just recall that the natural setting for taking care of classical solutions is the energy set 
\begin{equation*}
    \energyset:=\big\lbrace \psi\in H^1_{\loc}(\R) \big| \psi' \in L^2 (\R), F(|\psi|^2) \in L^1(\R) \big\rbrace 
\end{equation*}

endowed with the distance
\begin{equation}\label{métrique d}
    d(\psi_1,\psi_2)=\left\Vert \psi_1-\psi_2\right\Vert_{L^\ii([-1,1])}+\normLdeux{\eta_1-\eta_2}+\normLdeux{\psi_1'-\psi_2'}.
\end{equation}

In a dedicated appendix, we will see that the Cauchy problem is well-posed in $\mathcal{X}(\R)$ and that the energy is conserved along the flow. We shall also need to analyze the continuity of the flow, and the Cauchy problem for~\eqref{NLShydro}. We refer to Appendix~\ref{appendix: flow maps and weak continuity} for more details. For now, we just mention that the Cauchy problem for~\eqref{NLS} has been handled under several conditions on $f$ that are described below and which we shall assume throughout the article. Observing moreover that $\mathcal{X}(\R)+H^1(\R)\subset \mathcal{X}(\R)$, we can state the following result.
\begin{thm}[\cite{Gallo1}]\label{théorème: local global well-posedness of cauchy problem}
Let $u_0\in \mathcal{X}(\R)$. Take $f$ in $ C^2(\R)$ satisfying \eqref{hypothèse de croissance sur F minorant intermediaire} below. In addition, assume that there exist $\alpha_1\geq 1$ and $C_0 >0$ such that for all $\rho\geq 1$, 
\begin{equation*}\label{theoreme de gallo condition de croissance sur f''}\tag{H0}
|f''(\rho)|\leq \dfrac{C_0}{\rho^{3-\alpha_1}}.
\end{equation*}

If $\alpha_1 > \frac{3}{2}$, assume moreover that there exists $\alpha_2\in [ \alpha_1-\frac{1}{2},\alpha_1]$ such that for $\rho\geq 2$, $C_0 \rho^{\alpha_2} \leq F(\rho)$.\\
There exists a unique function $w \in  C^0\big(\R,H^1(\R)\big)$ such that $u:=u_0+w$ solves \eqref{NLS}. Moreover, the solution continuously depends  on the initial condition, and the energy $E$ and the momentum $p$ are conserved by the flow.
\end{thm}

In this article, we aim at understanding the dynamics of special solutions of the form $\Psi(t,x)=\gu_c(x-ct)$. These are called travelling waves and their profiles $\gu_c$ satisfy the ordinary differential equation 
\begin{equation}\label{TWC}\tag{$TW_c$}
    -ic\gu_c'+\gu_c''+\gu_c f(|\gu_c|^2)=0.
\end{equation}

In the Gross-Pitaevskii case, where $c_s=\sqrt{2}$, travelling waves exist with a one-to-one correspondence with the speed $c$ in the entire range $c\in(-c_s,c_s)$. In the case of a generic nonlinearity $f$, it is a much more subtle problem. Nonetheless, there can only exist nontrivial finite energy travelling waves with speed in $(-c_s,c_s)$ (see Theorem~5.1 and the remark just after in~\cite{Maris6}). Under suitable conditions on the nonlinearity $f\in C^3(\R_+)$, we have existence and uniqueness of transonic speeds. We assume that these conditions hold true in the sequel. Namely
\begin{itemize}
    \item For all $\rho\in\R$,
\begin{equation}\label{hypothèse de croissance sur F minorant intermediaire}\tag{H1}
    \dfrac{c_s^2}{4} (1-\rho)^2 \leq F(\rho).
\end{equation}
    \item There exist $M\geq 0$ and $q_*\in [2, +\ii)$ such that for all $\rho \geq 2$,
\begin{equation}\label{hypothèse de croissance sur F majorant}\tag{H2}
    F(\rho)\leq M|1-\rho|^{q_*}.
\end{equation}
    \item \begin{equation}\label{condition suffisante pour la stabilité orbitale sur f''(1)+6f'(1)>0}\tag{H3}
    k:=2f''(1)+6f'(1)< 0.
\end{equation}
\end{itemize}

Hypothesis~\eqref{condition suffisante pour la stabilité orbitale sur f''(1)+6f'(1)>0} can be relaxed to the condition $k \neq 0$. Indeed, the sign of some quantity $\xi_c$, introduced in Appendix~\ref{appendix: transonic régime}, is opposite\footnote{We refer to Subsection 2.2 in~\cite{Bert2} for more details.} to that of $k$. Moreover, under Hypothesis~\eqref{hypothèse de croissance sur F minorant intermediaire}, one must have\footnote{We refer to Remark 1.11 in~\cite{Bert1} for more details.}\ $\xi_c>0$ for all admissible speeds which implies $k<0$. Consequently, within this framework,~\eqref{condition suffisante pour la stabilité orbitale sur f''(1)+6f'(1)>0} is equivalent to $k\neq 0$.

Regarding the existence and properties of travelling waves for~\eqref{NLS}, we state a result proved in~\cite{Bert2}.
\begin{thm}[\cite{Bert2}]\label{théorème: il existe une branche C^1 de solitons proche de c_s}
    There exists a critical speed $c_0 >0$ such that for $c\in (c_0,c_s)$, there exists a non-constant and non-vanishing smooth solution $\gu_c$ of \eqref{TWC}, that is unique up to translations and constant phase shifts. 
    The hydrodynamical variables $Q_c:=(\eta_c,v_c)$ and classical variables are related by the formula
    \begin{equation}\label{formule de gu_c classique en fonction eta_c et v_c hydrodynamique}
        \gu_c(x)=\sqrt{1-\eta_c(x)}e^{-i\int_0^x v_c(r)dr}.
    \end{equation}
Moreover the hydrodynamical variables are smooth and $c\mapsto Q_c$ belongs to $ C^2\big((c_0,c_s),\mathcal{NX}_{hy}^2(\R)\big)$, with
\begin{equation}\label{théorème: hydrodynamique condition de grillakis sans le signe}
    \dfrac{d}{dc}\big(p(Q_c)\big) < 0.
\end{equation}

Finally there exist $a_{d},K_{d}>0$ independent of $c\in (c_0,c_s)$ and $x\in\R$ such that,
    \begin{equation}\label{estimée décroissance exponentielle à tout ordre pour eta_c et v_c}
\sum_{\substack{0\leq m_1\leq 3\\ 0\leq m_2\leq 2 \\ 0\leq j\leq 2}}(c_s^2-c^2)^{j-1}\Big(|\partial_c^j\partial_x^{m_1} \eta_c(x)|+c^{1+2j+2m_2}|\partial_c^j\partial_x^{m_2} v_c(x)|\Big)\leq K_{d}e^{-a_{d}\sqrt{c_s^2-c^2} |x|}.
\end{equation}
\end{thm}

One can now wonder whether these solutions are stable or not, and precise the sense of "stable". Taking an arbitrary initial condition close to a travelling wave $\gu_c$, then the solution associated with this initial condition will certainly not stay close to $\gu_c(.-ct)$ uniformly in time. Indeed, consider another travelling wave as an initial condition, with different speed $\widetilde{c}$. As close as these speeds are taken, the solutions associated with these initial conditions will separate into two different bubbles, which will be more and more apart the further they travel. To illustrate this claim, one can show that there exists a positive constant $\varepsilon_{c,\ctilde}$ only depending on $c$ and $\ctilde$ such that
\begin{equation}\label{dbig( gu_c(.-ct),gu_widetildec(.-widetildect)big)ttendii +ii}
    \sup_{ t\geq 0} d\big( \gu_c(.-ct),\gu_{\widetilde{c}}(.-\widetilde{c}t)\big)\geq \varepsilon_{c,\ctilde},
\end{equation}

thus the Lyapunov stability does not stand a fair chance of success. However, intermediate notions of stability exist for travelling waves, such as \textit{orbital} and \textit{asymptotic} stability. In this article, we will prove that these notions are well-suited and why a stronger notion could not be true. Recall that the asymptotic stability of travelling waves in the Gross-Pitaevskii case has been thoroughly investigated in~\cite{BetGrSm2}, even for the speed $c=0$, in~\cite{GravSme1}. Our study refrains from handling this speed for a significant reason, the black soliton (of speed $c=0$) vanishes and thus prevents us from considering any smooth lifting of $\gu_0$. However, we have to restrict ourselves to speeds even closer to $c_s$ and we assume in addition that the hypothesis of Theorem~\ref{théorème: local global well-posedness of cauchy problem}, as well as~\eqref{hypothèse de croissance sur F minorant intermediaire},~\eqref{hypothèse de croissance sur F majorant} and~\eqref{condition suffisante pour la stabilité orbitale sur f''(1)+6f'(1)>0}, are satisfied to state our main result.

\begin{thm}\label{théorème: stab asymp classical}
    Assume that~\eqref{theoreme de gallo condition de croissance sur f''}-\eqref{condition suffisante pour la stabilité orbitale sur f''(1)+6f'(1)>0} hold. There exists an integer $l_*$ and $c_1\in (0,c_s)$ such that if $f\in C^{l_*}(\R_+)$, then the travelling wave with a given speed in $(c_1,c_s)$ is asymptotically stable in the following sense. Let $c^*\in (c_1,c_s)$. There exists $\delta_{*} >0$ such that for any $\psi_0\in\energyset$ satisfying $d(\psi_0,\gu_{c^*})\leq \delta_{*}$, then $\psi(t)$ the solution associated with the initial condition $\psi_0$ is globally defined and there exist $c^*_0\in (c_1,c_s) $ and two functions $(b,\theta)\in C^1(\R_+,\R^2)$ such that 

    \begin{equation}\label{e^-itheta(t) psibig( t, .+b(t)big)tiitendfd gu_c^*_0}
        e^{-i\theta(t)} \psi\big( t, .+b(t)\big)\tiitendfd \gu_{c^*_0},
    \end{equation}

    and
    \begin{equation*}
        b'(t)\ttendii c^*_0\quad\quad\text{ and }\quad\quad \theta'(t)\ttendii 0.
    \end{equation*}

Here the notation $\Psi(t)\tiitendfd\Psi_\ii$ stands for the fact that all three following convergences hold \begin{equation}\label{définition convergence faible d}
\left\{\begin{array}{l}
    1-|\Psi(t)|^2\tiitendf 1-|\Psi_\ii|^2\quad\text{in }L^2(\R),\\
    \partial_x \Psi(t)\tiitendf\partial_x \Psi_\ii\quad\text{in }L^2(\R), \\
    \Psi(t)\ttendii\Psi_\ii\quad\text{in }L^\ii_\loc(\R) . \\
\end{array}
\right.
\end{equation}
\end{thm}

The proof of this theorem makes crucial use of hydrodynamical variables and relies on a method developed in~\cite{BetGrSm2} that itself is based on ideas introduced in~\cite{MartMer7} for the Korteweg-de Vries equation. This article follows the footsteps of the preceding articles~\cite{Bert1,Bert2} regarding the same kind of general nonlinearities. Indeed, it is of interest to identify the conditions that the nonlinearity $f$ needs to satisfy for the travelling wave to be asymptotically stable. More generally, this series of articles aim at understanding the behaviour of solutions for large times, while not using scattering methods, as it is possible in the integrable case of the Gross-Piteasvkii equation~\cite{GeraZha1,CuccJen1}. Unlike the latter articles, the proof relies here on PDEs arguments, already highlighted in integrable cases as in~\cite{BetGrSm2} and then for non-integrable systems in~\cite{MartMer7,MarMeTs2,MartMer2,MartMer4}.

\begin{rem}
    Unlike the Gross-Pitaevskii case, we need to restrict ourselves to small travelling waves. This corresponds to travelling waves with speeds taken close to $c_s$, the study of which is developed in Appendix~\ref{appendix: transonic régime}. The special choice of these travelling waves is consistent with the branch of solitons already exhibited in Theorem~\ref{théorème: il existe une branche C^1 de solitons proche de c_s}. A further restriction will provide positivity on crucial quantities, which are necessary for proving Proposition~\ref{lemme: rigidity property proche d'un soliton, version algébrique} below.
\end{rem}

\begin{rem}
The modulation parameters in~\eqref{e^-itheta(t) psibig( t, .+b(t)big)tiitendfd gu_c^*_0} are closely related to the invariances of~\eqref{NLS}. More precisely, the asymptotic stability takes into account the symmetries that prevent the solution from being stable in the Lyapunov sense, that are translations and constant phase shifts.
\end{rem}

\begin{rem}
In contrast with the Korteweg-de Vries equation, negative speeds are allowed and the asymptotic stability of the corresponding travelling waves can be reached by the same method by simply taking the complex conjugate.
\end{rem}

\begin{rem}\label{remarque: value of l_*}
The value of $l_*$ is fixed during the proof (see Subsection~\ref{section: Convergence along the evolution trightarrow +ii}), so that~\eqref{NLS} is sufficiently differentiable with respect to the spatial variable. This differentiability is used to obtain a sufficient decay for the solution $\psi$ that is crucial for characterizing a travelling wave as the limit object in~\eqref{e^-itheta(t) psibig( t, .+b(t)big)tiitendfd gu_c^*_0}.
\end{rem}

The main ingredients in the proof of Theorem~\ref{théorème: stab asymp classical} stand in three main ideas. Orbital stability allows us to construct a compact profile with speed $c_0^*$ towards which the solution $\psi$ weakly tends as time goes to $+\ii$. Furthermore, as previously analyzed in works such as~\cite{EsKePoV5}, the structure of~\eqref{NLShydro} and the exponential decay of the hydrodynamical travelling waves imply a smoothing effect for this profile, as well as some uniform exponential decay. This property is crucial to prove a rigidity theorem for solutions taken close to a stable travelling wave, which eventually provides Theorem~\ref{théorème: stab asymp classical}.

\subsection{The hydrodynamic framework}
We first recall that a travelling wave of speed $c$ close to $c_s$ neither vanishes nor does any perturbation of such a wave, as ensured by the Sobolev embedding theorem. Non-vanishing solutions can be lifted as $\psi=|\psi|e^{i\varphi}$, hence we can construct the hydrodynamical variables given by $Q=(\eta,v)=\big(1-|\Psi|^2,-\partial_x \varphi\big)$. In this framework, the natural energy set appears to be $\energysethydro:=H^1(\R)\times L^2(\R)$. This set corresponds to the particular case $l=0$ in the class of spaces $\energysethydrok=H^{l+1}(\R)\times H^l(\R)$ for $l\geq 0$ endowed with the euclidean norms
\begin{equation*}
    \normXk{Q}^2=\Vert \eta\Vert_{H^{l+1}}^2 + \Vert v\Vert_{H^l}^2.
\end{equation*}

One major difference compared to the original setting is that~\eqref{NLShydro} is locally well-posed on a vector space (see Appendix~\ref{appendix: flow maps and weak continuity} for more details on the Cauchy problem) instead of a complete metric subspace. We also define the $\mathcal{X}$-weak convergence as follows. For $Q_n,Q\in\energysethydro$, we naturally write that \begin{equation*}
    Q_n\ntendfX Q,
\end{equation*}
when
\begin{equation*}\label{définition convergence faible X}
\left\{\begin{array}{l}
    \eta_n\ntendf \eta\quad\text{in }H^1(\R),\\
    v_n\ntendf v\quad\text{in }L^2(\R).\\
\end{array}
\right.
\end{equation*}

Note that the first convergence in~\eqref{définition convergence faible d} is a straightforward consequence of the weal convergence in $\energysethydro$. We can now state the asymptotic stability in the hydrodynamic frameworl.

\begin{thm}\label{théorème: stab asymp hydro}
    Let $c^*\in (c_1,c_s)$. There exists $\beta_{*} >0$ such that for any $Q_0\in\energysethydro$ satisfying $d(Q_0,Q_{c^*})\leq \beta_{*} $, then $Q(t)$ the solution associated with the initial condition $Q_0$ is globally defined and there exist $c_0^*\in (c_1,c_s) $ and a function $a\in C^1(\R_+)$ such that 

    \begin{equation*}
        Q\big(t,.+a(t)\big)\tiitendfX Q_{c^*_ 0}.
    \end{equation*}
\end{thm}

We will eventually derive asymptotic stability in the original frameworl, that is Theorem~\ref{théorème: stab asymp classical}, from Theorem~\ref{théorème: stab asymp hydro}.

\subsection{Orbital stability of a travelling wave}\label{section: minimizing property of the travelling waves}
We have seen that the notion of stability needs to be properly defined. Typically, to understand the instability phenomenon described by~\eqref{dbig( gu_c(.-ct),gu_widetildec(.-widetildect)big)ttendii +ii}, one has to take the symmetries of the equation into account. In other words, one can view the invariances of~\eqref{NLS} as degrees of freedom regarding the possible different types of evolution for a solution. Orbital stability is designed in that purpose. We refer to~\cite{Bert1,Bert2} where this question is solved, for more details regarding the following statements. Under the previous conditions~\eqref{hypothèse de croissance sur F minorant intermediaire},~\eqref{hypothèse de croissance sur F majorant} and~\eqref{condition suffisante pour la stabilité orbitale sur f''(1)+6f'(1)>0} on the nonlinearity $f\in C^3(\R_+)$ and for given speed $c^*\in (c_0,c_s)$, the hydrodynamic travelling wave $Q_{c^*}$ is orbitally stable in the following sense.

\begin{thm}[\cite{Bert2}]\label{théorème: stabilité orbitale hydro}
    Let $c^*\in (c_0,c_s)$. There exists $\alpha_*,A_*>0$, such that the following holds. If $Q_0=(\eta_0,v_0)\in\energysethydro$ is such that for some $a^0\in\R$ ,

\begin{equation}\label{alpha_0:=normXQ_0-R_c^*,a^0leq alpha_*}
    \alpha_0:=\normX{Q_0-Q_{c^*}(.-a^0)}\leq \alpha_*,
\end{equation}

then, the unique solution $Q(t)=\big(\eta(t),v(t)\big)$ to~\eqref{NLShydro} associated with the initial datum $(\eta_0,v_0)$ is globally defined. Moreover there exists $(a,c)\in C^1(\R_+,\R^{2})$ such that for any $t\in\R_+$, $Q(t)$ can be decomposed as 
\begin{equation}\label{modulation décomposition}    Q\big(t,.+a(t)\big)=Q_{c(t)}+\varepsilon(t),
\end{equation}

where 
\begin{equation}\label{normXvarepsilon(t+normRc(t)-c^*+normRa(t)-a^leq A_*alpha_0}
    \normX{\varepsilon(t)}+\normR{c(t)-c^*}\leq A_*\alpha_0,
\end{equation}

and
\begin{equation}\label{normXQ(t)-Q_c^*,a(t)^2+normR a'(t)-c(t)^2+normR c'(t)leq A_* alpha_0^2}
    \normR{ a'(t)-c(t)}^2+\normR{ c'(t)}\leq A_* \normX{\varepsilon(t)}^2.
\end{equation}
\end{thm}

Furthermore, recall that we construct the functions $a$ and $c$ in Theorem~\ref{théorème: stabilité orbitale hydro} by an implicit function argument such that we have both the orthogonality conditions
    \begin{equation}\label{condition d'orthogonalité sur varepsilon dans décomposition orthogonale}
    \psLdeuxLdeux{\varepsilon(t)}{\partial_x Q_{c(t)}}=\nabla p(Q_{c(t)}).\varepsilon(t) =0,
\end{equation}

in the sense of the definition of the momentum given in~\eqref{expression du moment}. From a spectral analysis of the linearized operator $\mathcal{H}_c$, we also recall that under the previous orthogonality conditions, there exists a positive constant $l_{c(t)}$ such that 
    \begin{equation}\label{coercivité de la forme quadratique H_c sous condition d'orthogonalité}
    \psLdeuxLdeux{\mathcal{H}_{c(t)}\big(\varepsilon(t)\big)}{\varepsilon(t)}\geq l_{c(t)}\normX{\varepsilon}^2.
\end{equation}

\begin{rem}
This result expresses the fact that we can modulate a perturbation of a travelling wave $Q_c$ so that it stays close to this travelling wave as time evolves. The modulation consists of modifying the speed and translating the travelling wave along the evolution, and the so-called parameters are given by the functions $a$ and $c$ above. Moreover,~\eqref{normXvarepsilon(t+normRc(t)-c^*+normRa(t)-a^leq A_*alpha_0} and~\eqref{normXQ(t)-Q_c^*,a(t)^2+normR a'(t)-c(t)^2+normR c'(t)leq A_* alpha_0^2} provide some control over these parameters. As a straightforward consequence of this result, a compact profile emerges and then appears to be a suitable candidate to feature the asymptotic behavior of the solution $Q(t)$.
\end{rem}

\begin{rem}
    Taking the limit $\alpha_0\rightarrow 0$ in both~\eqref{normXvarepsilon(t+normRc(t)-c^*+normRa(t)-a^leq A_*alpha_0} and~\eqref{normXQ(t)-Q_c^*,a(t)^2+normR a'(t)-c(t)^2+normR c'(t)leq A_* alpha_0^2} formally implies that the translation parameters behave like a linear function $a(t)\sim mt$ and then $Q(t)\sim Q_{c^*}(.-mt)$. It is of interest to know whether the translation parameter $a(t)$ behaves as an actual linear function of $t$. However, this does not stand a fair chance of success since it was shown in~\cite{MartMer4} for the generalized Korteweg-de Vries equation, that the analogue translation parameters $a(t)$ can behave, up to a renormalization, as $mt +\sqrt{\ln(t)}$.
\end{rem}

\begin{rem}
    The only difference with the original result in~\cite{Bert2} lies in the quadratic control in~\eqref{normXQ(t)-Q_c^*,a(t)^2+normR a'(t)-c(t)^2+normR c'(t)leq A_* alpha_0^2}. This is justified in Appendix~\ref{appendix: propriété du linéarizé}.
\end{rem}

In the sequel, we introduce the main ideas of the proof of Theorem~\ref{théorème: stab asymp classical} and how they are articulated.

\subsection{Sketch of the proof}

\subsubsection{Construction of the profile in hydrodynamic variables}

Using orbital stability, we can construct a limit profile, up to a subsequence, by using~\eqref{normXvarepsilon(t+normRc(t)-c^*+normRa(t)-a^leq A_*alpha_0}. Indeed, the functions $\varepsilon$ and $c$ are bounded respectively in $\energysethydro$ and $\R$. Given a speed $c^*\in (c_0,c_s)$, an initial condition $Q_0$ satisfying~\eqref{alpha_0:=normXQ_0-R_c^*,a^0leq alpha_*} and a sequence of time $(t_n)$ tending to $+\ii$, we can suppose, up to a subsequence, that there exist $\widetilde{\varepsilon}_0\in\energysethydro$ and $c_0^*\in\R$ such that 
\begin{equation}\label{convergence faible vers le profil limite}
    \left\{\begin{array}{l}
    \varepsilon(t_n)\ntendfX \widetilde{\varepsilon}_0,\\
    c(t_n)\ntend c_0^*. \\
\end{array}
\right.
\end{equation}

Now, by local well-posedness, we can also construct $\widetilde{Q}(t)$ the unique solution of~\eqref{NLShydro} with initial condition $\widetilde{Q}_0:=Q_{c_0^*}+\widetilde{\varepsilon}_0$. Using both convergences~\eqref{convergence faible vers le profil limite}, then the control~\eqref{normXQ(t)-Q_c^*,a(t)^2+normR a'(t)-c(t)^2+normR c'(t)leq A_* alpha_0^2} with a potentially smaller $\alpha_*$, and the Lipschitz continuity of $c\mapsto Q_c$ as stated in Lemma A.1 in~\cite{Bert2}, we obtain
\begin{align}
    \normX{\widetilde{Q}_0-Q_{c^*}}&\leq \normX{\widetilde{\varepsilon}_0}+K_{lip}|c^*-c_0^*|\notag\\ 
    &\leq \liminf_{n\rightarrow +\ii}\normX{\varepsilon(t_n)} + K_{lip}\limsup_{n\rightarrow +\ii}|c^*-c(t_n)|\leq \alpha_* \label{widetildeQ_0-Q_c^* leq alpha_*}.
\end{align}

Therefore, by orbital stability, $\widetilde{Q}(t)$ is globally defined and there exist two functions $\widetilde{c},\widetilde{a}$ such that for any $t\in\R$, we have the decomposition 

\begin{equation}\label{décomposition orthogonale de widetile Q}
    \widetilde{Q}(t,.+\widetilde{a}(t))=Q_{\widetilde{c}(t)}+\widetilde{\varepsilon}(t),
\end{equation}

with for all $t\in\R$,
\begin{equation}\label{controle de widetilde varepsilon et ctilde -c_0^* par alpha_*}
    \normX{\widetilde{\varepsilon}(t)}+\normR{\widetilde{c}(t)-c_0^*}\leq A_*\normX{\widetilde{Q}_0-Q_{c^*}}\leq A_*\alpha_*,
\end{equation}

and
\begin{equation}\label{controle des parametres de modulation}
    \normR{ \widetilde{a}'(t)-\ctilde(t)}^2+\normR{ \ctilde'(t)}\leq A_* \normX{\varepsilontilde(t)}^2.
\end{equation}

For $c\in [0,c_s]$, we define
\begin{equation*}
    \nu_c:=\sqrt{c_s^2-c^2}.
\end{equation*}
From both previous estimates~\eqref{controle de widetilde varepsilon et ctilde -c_0^* par alpha_*} and~\eqref{controle des parametres de modulation}, and up to taking a smaller $\alpha_*$ in~\eqref{widetildeQ_0-Q_c^* leq alpha_*}, we can deduce some uniform control in time for several crucial quantities, such as $\nu_c$. In this sense, we state the following result, and we refer to Corollary~1.12 in~\cite{Bert2} for the proof.

\begin{prop}\label{proposition: borne uniforme sur les quantités nu_c etc}
    There exists $\iota_*>0$ and $I_*\in (0,\frac{3}{2})$ depending only on $c^*$ such that for any $t\in\R$
\begin{equation*}
    \min\left\{\widetilde{c}(t),\nu_{\widetilde{c}(t)},\inf_\R(1-\eta_{\widetilde{c}(t)}),\inf_\R (1-\etatilde),l_{\ctilde(t)},-\dfrac{d}{dc}\big(p(Q_{\ctilde(t)}\big),-\dfrac{d}{dc}\big(p(Q_{c(t)}\big)\right\}\geq \iota_*,
\end{equation*}
and
\begin{equation*}
    \max\left\{\normLii{\widetilde{\eta}(t)},|\widetilde{a}'(t)-c^*_0|,\normLii{\widetilde{\psi}(t)}\right\}\leq I_* .
\end{equation*}
\end{prop}

Moreover, we shall use the weak continuity of the flow in order to obtain some convergences similar to~\eqref{convergence faible vers le profil limite}, but along the time evolution. That is the aim of the next proposition, the proof of which can be found at the end of Appendix~\ref{appendix: flow maps and weak continuity}.

\begin{prop}\label{proposition: convergence faible vers le profil limite le long de l'évolution}
Let $t\in\R$. Then we have 
\begin{equation*}\label{Qbig(t_n+t,.+a(t_n)big)ntendfX widetilde{Q}(t)}
    Q\big(t_n+t,.+a(t_n)\big)\ntendfX \widetilde{Q}(t),
\end{equation*}

as well as
\begin{equation*}\label{a(t_n+t)-a(t_n)ntend a^*(t)}
    a(t_n+t)-a(t_n)\ntend \widetilde{a}(t)\quad\text{and}\quad c(t_n+t)\ntend \widetilde{c}(t).
\end{equation*}
\end{prop}

\begin{rem}\label{remarque: widetildea(0)=0 and widetildec(0)=c_0^*}
Proposition~\ref{proposition: convergence faible vers le profil limite le long de l'évolution} implies that $\widetilde{a}(0)=0$ and $\widetilde{c}(0)=c_0^*$.
\end{rem}

\subsubsection{Liouville property around a travelling wave}\label{section: Liouville property around a travelling wave}
The following step of the proof is to establish a rigidity property for smooth decaying solutions close to travelling waves. In fact, we shall prove that such solutions have no other choice but to be exactly a travelling wave. In the sequel, the previous profile $\widetilde{Q}$ shall be shown to decay sufficiently fast, fulfilling the conditions~\eqref{lemme: decaying property in rigidity} in the next proposition. Eventually, since $\widetilde{Q}$ meets the assumptions of this rigidity property, this will impose that it is equal to a travelling wave.

\begin{prop}\label{lemme: rigidity property proche d'un soliton, version algébrique}
There exist an integer $l_0$ and a speed $c_1\in (0,c_s)$ such that if $c_0^*\in (c_1,c_s)$, then the following holds. If there exist $r > 5/2$ and $C>0$ such that for any $(t,x)\in\R^2$ and any integer $l\in\{0,...,l_0\}$, we have
\begin{equation}\label{lemme: decaying property in rigidity}
    \left|\partial_x^{l+1}\widetilde{\eta}\big(t,x+\widetilde{a}(t)\big)\right|+\left|\partial_x^{l}\widetilde{\eta}\big(t,x+\widetilde{a}(t)\big)\right|+\left|\partial_x^{l}\widetilde{v}\big(t,x+\widetilde{a}(t)\big)\right|\leq \dfrac{C}{1+|x|^r},
\end{equation}
then we obtain $\widetilde{a}(t)=c_0^*t,\ \widetilde{c}(t)=c^*_0$ and for any $t\in\R$,\begin{equation*}
    \widetilde{Q}(t)=Q_{c^*_0}\big(.-c^*_0t\big).
\end{equation*}
\end{prop}

The rigidity result is stated in the frame of the decomposition with modulation parameters. In its proof, we make crucial use of the smallness of the considered travelling waves. This is related to the bound from below in Lemma~\ref{lemme: monotonicity formula sur le viriel} below, that implies Proposition~\ref{lemme: rigidity property proche d'un soliton, version algébrique} and that we are able to prove only for $c$ close enough to $c_s$. We refer to Claim~\ref{claim: spectre essentiel de T_c} and the remark just above to understand how the restriction to transonic speeds is crucial in the proof of Lemma~\ref{lemme: monotonicity formula sur le viriel}.

\begin{rem}
In the sequel, we will show that it is sufficient to take $l_0=12$. Then in a second step and consistently with Remark~\ref{remarque: value of l_*}, the choice of $l_*$ for which the main result holds will depend on the value of $l_0$.
\end{rem}

\begin{rem}\label{remarque: décroissance sur widetilde Q implique celle de varepsilontilde}
    Combining both the exponential decay in~\eqref{estimée décroissance exponentielle à tout ordre pour eta_c et v_c} of the travelling waves, and the decay~\eqref{lemme: decaying property in rigidity} implies that $\widetilde{\varepsilon}$ and its derivatives decay similarly to~\eqref{lemme: decaying property in rigidity}. Indeed, we infer that there exists $C_*>0$ such that for any $(t,x)\in\R^2$and $l\in\{0,...,l_0\}$, 
\begin{align*}\label{decaying property in rigidity for varepsilon^*}
    |\partial_x^{l+1}\widetilde{\varepsilon_\eta}(t,x)|+|\partial_x^l\widetilde{\varepsilon_\eta}(t,x)|+|\partial_x^l\widetilde{\varepsilon_v}(t,x)|&\leq \dfrac{\widetilde{C}}{1+|x|^r}+K_d e^{-a_d\iota_{*}|x|}\leq\dfrac{C_*}{1+|x|^r}.
\end{align*}
\end{rem}

\begin{rem}
    Once $c_1$ is exhibited (see the proof of Lemma~\ref{lemme: monotonicity formula sur le viriel}), the condition $c_0^*\in (c_1,c_s)$ can be imposed by taking $c^*\in (c_1,c_s)$ and potentially reducing the value of $\alpha_*$. 
\end{rem}

\begin{rem}
    Proposition~\ref{lemme: rigidity property proche d'un soliton, version algébrique} provides a classification for solutions of~\eqref{NLShydro} initially taken sufficiently close to a stable travelling wave and enjoying the decaying property~\eqref{lemme: decaying property in rigidity}. Furthermore, we think that this classification can be proved if the next less restrictive property is satisfied. Given $\epsilon >0$, there exists $R>0$ such that for any $t\in\R$,
\begin{equation*}
    \int_{|x-\widetilde{a}(t)|\geq R}\pi\big(\widetilde{Q}(t,x)\big)dx \leq \epsilon,
\end{equation*}

where $\pi$ is the momentum density defined in~\eqref{expression du moment}. This uniform compactness condition is set consistantly with the assumptions of the analogous Liouville theorem stated in~\cite{Martel2,MartMer2}. Under this condition, Proposition~\ref{lemme: rigidity property proche d'un soliton, version algébrique} prescribes the structure of a solution taken close to a travelling wave and that does not disperse. There is no other choice for these solutions to be an exact travelling wave.
\end{rem}

\begin{rem}
One contribution of this article lies in the fact that the rigidity property can be shown for solutions that have algebraic decay instead of exponential decay. Precisely, for higher dimensions, the travelling waves in the Gross-Pitaevskii case are known to decay with an algebraic rate~\cite{Graveja3}.
\end{rem}

The proof of Proposition~\ref{lemme: rigidity property proche d'un soliton, version algébrique} can be divided into several steps. We first focus on the equation satisfied by $\widetilde{\varepsilon}(t)$. The Hamiltonian structure of~\eqref{NLShydro} can be formally expressed for a generic non-vanishing solution $Q$ as

\begin{equation}\label{equation hamiltonienne vérifiée par Q avec J et grad E}
    \partial_t Q= J\grad E(Q),
\end{equation}
where
\begin{equation}\label{expression de J}
    \begin{array}{l}
    JQ=-2S\partial_x Q \\
\end{array}
\quad\text{with}\quad S:=\begin{pmatrix}
        0 & 1\\
        1 & 0
    \end{pmatrix}
\end{equation}

and where the expression of the hydrodynamic energy is  
\begin{equation}\label{expression de l'énergie}
    E(Q)=E(\eta,v):=\int_\R e(\eta,v):=\dfrac{1}{8}\int_\R\dfrac{(\partial_x\eta)^2}{1-\eta}+\dfrac{1}{2}\int_\R (1-\eta)v^2+\dfrac{1}{2}\int_\R F(1-\eta).
\end{equation}

We recall that the expression of the momentum is given by 
\begin{equation*}
    p(Q)=p(\eta,v):=\int_\R\pi(\eta,v):=\dfrac{1}{2}\int_\R \eta v .
\end{equation*}

Both previous quantities are conserved along the flow, according to the local well-posedness results in Appendix~\ref{appendix: flow maps and weak continuity}. Moreover, provided that $f\in C^k(\R_+)$, then the energy is $C^{k+1}$ on $\Nenergysethydro$ and the momentum $p$ is $C^\ii$ on $L^2(\R)^2$ with $\grad p(\eta,v)=\frac{1}{2}(v,\eta)$, regardless from the regularity of $f$. In the hydrodynamic framework, we recall that~\eqref{TWC} reads
\begin{equation}\label{nabla (E-cp)=0}
    \nabla (E-cp)(Q_c)=0,
\end{equation} 

and can be viewed as the Euler-Lagrange equation associated with a minimization of the energy with fixed momentum (see~\cite{Bert1}). Using~\eqref{nabla (E-cp)=0}, we can linearize~\eqref{equation hamiltonienne vérifiée par Q avec J et grad E} around the travelling wave of speed $c$ and this leads to
\begin{equation*}\label{linearisation de nlshydro autour d'un travelling wave}
    \partial_t \varepsilon = J\mathcal{H}_c (\varepsilon)
\end{equation*}
where $\mathcal{H}_c=\nabla^2 (E-cp)(Q_c)$ is the linearized operator around the travelling wave $Q_c$. Among the properties of the linearized operator that are analyzed in~\cite{Bert2}, we recall that $\mathcal{H}_c$ is a self-adjoint operator with domain $H^2(\R)\times L^2(\R)$. It owns a unique negative direction, a unique vanishing direction, and $\ker(\mathcal{H}_c)=\spanned{Q_c'}$. Now, using~\eqref{equation hamiltonienne vérifiée par Q avec J et grad E},\eqref{nabla (E-cp)=0} and the decomposition~\eqref{modulation décomposition} with $c=\widetilde{c}(t)$, we obtain the equation satisfied by $\widetilde{\varepsilon}(t)=\widetilde{Q}\big(t,.+\widetilde{a}(t)\big)-Q_{\widetilde{c}(t)}$, that is

\begin{align*}
    \partial_t\widetilde{\varepsilon}  = &  J\Big(\grad E\big(\widetilde{Q})-\grad E(Q_{\widetilde{c}(t)}) -\grad^2 E(Q_{\widetilde{c}(t)}).\widetilde{\varepsilon} \Big) + J\mathcal{H}_{\widetilde{c}(t)}(\widetilde{\varepsilon}) \\
    &+ \widetilde{c}(t)J\Big(\grad p \big(Q_{\widetilde{c}(t)}\big)+\grad^2 p\big(Q_{\widetilde{c}(t)}\big).\varepsilontilde\Big)
    + \widetilde{a}'(t)\partial_x \widetilde{Q}\big(t,.+\widetilde{a}(t)\big)- \widetilde{c}'(t)\partial_c Q_{\widetilde{c}(t)}\\
\end{align*}

Due to the very generic identities for any $Q,\varepsilon\in\energysethydro$, \begin{equation*}
    J\grad p(Q)=-\partial_x Q\text{ and }J\grad^2 p(Q).\varepsilon=-\partial_x\varepsilon ,
\end{equation*} 

we are led to
\begin{equation}\label{equation sur variable varepsilon}
     \partial_t\widetilde{\varepsilon} = J\mathcal{H}_{\widetilde{c}(t)}(\widetilde{\varepsilon})+J\mathcal{R}_{\widetilde{c}(t)}(\widetilde{\varepsilon})  + \big(\widetilde{a}'(t)-\widetilde{c}(t)\big)\big(\partial_x Q_{\widetilde{c}(t)}+\partial_x\widetilde{\varepsilon}(t)\big)-\widetilde{c}'(t)\partial_c Q_{\widetilde{c}(t)},
\end{equation}

where \begin{equation}\label{définition mathcal R_c}
    \mathcal{R}_{\widetilde{c}(t)}(\widetilde{\varepsilon})=\grad E\big(Q_{\widetilde{c}(t)}+\widetilde{\varepsilon}\big)-\grad E(Q_{\widetilde{c}(t)}) -\grad^2 E(Q_{\widetilde{c}(t)}).\widetilde{\varepsilon} .
\end{equation}

Taking the first-order approximation $\mathcal{O}\big(\normX{\widetilde{\varepsilon}}^2\big)=0$ and imposing $|\widetilde{a}'(t)-\widetilde{c}(t)|^2+|\widetilde{c}'(t)|=\mathcal{O}\big(\normX{\widetilde{\varepsilon}(t)}^2\big)$ consistantly to the orbital stability, we obtain the approximate linearized equation
\begin{equation}\label{linéarization de nlshydro dans l'approximation au premier ordre}
    \partial_t\widetilde{\varepsilon} \approx J\mathcal{H}_{\widetilde{c}(t)}(\widetilde{\varepsilon}) + \big(\widetilde{a}'(t)-\widetilde{c}(t)\big)\partial_x Q_{\widetilde{c}(t)}.
\end{equation}

In order to obtain suitable estimates for $\widetilde{\varepsilon}$, we introduce a dual problem. Formally applying $\mathcal{H}_{\widetilde{c}(t)}$ to~\eqref{linéarization de nlshydro dans l'approximation au premier ordre} and using that the null-space of $\mathcal{H}_{\widetilde{c}(t)}$ is spanned by $\partial_x Q_{\widetilde{c}(t)}$, we get $\mathcal{H}_{\widetilde{c}(t)}(\partial_t\widetilde{\varepsilon}) \approx \mathcal{H}_{\widetilde{c}(t)}J\mathcal{H}_{\widetilde{c}(t)}(\widetilde{\varepsilon})$. This leads to the definition of the new variable $\widetilde{e}(t):=S\mathcal{H}_{\widetilde{c}(t)}\big(\widetilde{\varepsilon}(t)\big)$, which we also denote $\widetilde{e}=(\widetilde{e}_\eta,\widetilde{e}_v)$. The main properties of this formulation is stated below.

\begin{lem}\label{lemme: equation sur la variable duale epsilon^*}
We assume that the hypothesis of Proposition~\ref{lemme: rigidity property proche d'un soliton, version algébrique} hold true.\\
Then $\widetilde{\varepsilon}\in C^0\big(\R,\Nenergysethydrokk{4}\big)\cap C^1\big(\R,\Nenergysethydrokk{2}\big)$ and $\widetilde{e}\in C^0\big(\R,\Nenergysethydrokk{2}\big)\cap C^1\big(\R,\Nenergysethydro\big)$. Furthermore, we have the equation 
\begin{align}\label{equation sur variable duale epsilon}
    \partial_t \widetilde{e} = &-2S\mathcal{H}_{\widetilde{c}(t)}(\partial_x \widetilde{e})+S\mathcal{H}_{\widetilde{c}(t)}(J\mathcal{R}_{\widetilde{c}(t)}\widetilde{\varepsilon}) - \widetilde{c}'(t)S\mathcal{H}_{\widetilde{c}(t)}(\partial_c Q_{\widetilde{c}(t)})\notag\\ 
    &+\widetilde{c}'(t)S\partial_c\mathcal{H}_{\widetilde{c}(t)}(\widetilde{\varepsilon})+\big((\widetilde{a}'(t)-\widetilde{c}(t)\big)S\mathcal{H}_{\widetilde{c}(t)}(\partial_x\widetilde{\varepsilon})\notag.
\end{align}

In addition, there exists $D_*>0$ only depending on $c^*$ such that for any $t\in\R$,
\begin{equation}\label{controle de varepsilon^* par widetildee}
    \normX{\widetilde{\varepsilon}(t)}\leq D_*\normX{\widetilde{e}(t)}.
\end{equation}
\end{lem}

Proposition~\ref{lemme: rigidity property proche d'un soliton, version algébrique} next relies on the introduction of a time-dependent virial argument related to a matrix-valued function taking the form
\begin{equation}\label{définition N(t)}
    N(t)=\begin{pmatrix}
        0 & x \\
        x & 0
    \end{pmatrix} + \gamma^* M_{\widetilde{c}(t)},
\end{equation}
where \begin{equation}
    M_{\ctilde}(x):=\begin{pmatrix}
        m_{1,\ctilde} & m_{2,\ctilde}\\
        m_{2,\ctilde} & 0
    \end{pmatrix} \quad \text{ with }m_{1,\ctilde}:=-\dfrac{\eta_{\ctilde}'}{\eta_{\ctilde}}\text{ and }m_{2,\ctilde}:=-\dfrac{{\ctilde}\eta_{\ctilde}'}{2(1-\eta_{\ctilde})^2}.
\end{equation}

\begin{rem}
The matrix $M_c$ is set so that we have $\psLdeuxLdeux{M_{\ctilde} Q_{\ctilde}}{\etilde(t)}=0$. Indeed, $M_{\ctilde}Q_{\ctilde}=-S\partial_x Q_{\ctilde}$ which is orthogonal to $S\mathcal{H}_{\ctilde}\big(\varepsilon(t)\big)$.
\end{rem}

We claim that the matrix $M_{\ctilde}$ is uniformly bounded with respect to $c$ and $x$. Indeed, the function $m_{2,\ctilde}$ is bounded using the exponential decay properties of the travelling waves~\eqref{estimée décroissance exponentielle à tout ordre pour eta_c et v_c} and the uniform lower bound in Proposition~\ref{proposition: borne uniforme sur les quantités nu_c etc}. As for $m_{1,\ctilde}$, possibly taking $c_1$ even closer to $c_s$ so that Proposition~\ref{proposition: borne inférieure eta_c(x)} applies and using once again Proposition~\ref{proposition: borne uniforme sur les quantités nu_c etc}, we infer that there exists a bound $M_*$, only depending on $c^*$ such that $\normLii{M_{\ctilde(t)}}\leq M_*$. The virial quantity is then defined as
\begin{equation*}
    n(t):=\psLdeuxLdeux{N(t)\widetilde{e}(t)}{\widetilde{e}(t)}.
\end{equation*} 

The first matrix in~\eqref{définition N(t)} gives some coercivity at infinity for the time derivative $n'(t)$. The matrix $M_{\ctilde(t)}$ is designed to get some local coercivity for $n'(t)$. The constant $\gamma_*$ is eventually chosen to obtain the full coercivity in the next lemma. We point out that this argument relies significantly on the smallness of the travelling waves involved, which is why we must restrict to speeds close to $c_s$.

\begin{lem}\label{lemme: monotonicity formula sur le viriel}
    We assume that the hypothesis of Proposition~\ref{lemme: rigidity property proche d'un soliton, version algébrique} hold true. Then $n$ is well-defined and differentiable on $\R$. Furthermore, there exists $c_1\in (0,c_s)$ such that the following holds. If $c_0^*\in (c_1,c_s)$ then there exists a constant $\Xi_*>0$ only depending on $c^*_0$ such that for any $t\in\R$,
    \begin{equation*}\label{monotonie sur viriel}
        n'(t)\geq \Xi_*\normX{\widetilde{e}(t)}^2.
    \end{equation*}
\end{lem}

We deduce the proof of Proposition~\ref{lemme: rigidity property proche d'un soliton, version algébrique} from both Lemmas~\ref{lemme: equation sur la variable duale epsilon^*} and~\ref{lemme: monotonicity formula sur le viriel}.

\begin{proof}[Proof of Proposition~\ref{lemme: rigidity property proche d'un soliton, version algébrique}]

By the Cauchy-Schwarz inequality, Proposition~\ref{proposition: propriété du linéarisé mathcal H_c} and then Remark~\ref{remarque: décroissance sur widetilde Q implique celle de varepsilontilde}, we obtain
\begin{align}
    \big|\psLdeuxLdeux{N(t)\widetilde{e}(t)}{\widetilde{e}(t)}\big|&\leq 2\normLdeuxLdeux{\sqrt{|x|}\widetilde{e}(t)}^2+2\gamma_* M_*\normLdeuxLdeux{\widetilde{e}(t)}^2\notag\\
    &\leq 2K_*^2\int_\R \left((\partial^2_x\varepsilontilde_\eta)^2+(\partial_x\varepsilontilde_\eta)^2+\varepsilontilde_\eta^2 + (\partial_x\varepsilontilde_v)^2+\varepsilontilde_v^2\right)(1+|x|)dx\notag\\
    &\leq 2C_*^2 K_*^2\normLun{\dfrac{1+|x|}{(1+|x|^r)^2}}.\label{leq widetildC_d^2left(2normLdeux}
\end{align}

Since $r>1$, the upper bound in the previous estimate is finite. Combining both the bound~\eqref{leq widetildC_d^2left(2normLdeux} and Lemma~\ref{lemme: monotonicity formula sur le viriel} yields to the fact that the integral 
\begin{equation*}
    \int_\R \normX{\etilde(t)}^2 dt
\end{equation*}
is also finite. As a consequence of this, there exists a sequence of times $(t_n)$ tending to $+\ii$ such that \begin{equation*}
    \normX{\widetilde{e}(t_n)}^2\ntend 0,
\end{equation*}

and as a consequence of~\eqref{controle de varepsilon^* par widetildee}, 
\begin{equation*}\label{varepsilon^*(t_n)^2ntend 0}
    \normX{\widetilde{\varepsilon}(t_n)}^2\ntend 0.
\end{equation*}

To conclude, we use the orbital stability reversely in time by considering $Q^n_0:=\widetilde{Q}\big(t_n,.+\widetilde{a}(t_n)\big)$ as initial data, writing $Q^n(t)=\widetilde{Q}\big(t_n+t,.+\widetilde{a}(t_n)\big)$ the corresponding solution to~\eqref{NLShydro} and $\varepsilon^n(t)=\widetilde{Q}(t+t_n,.+\widetilde{a}(t+t_n)\big)-Q_{\widetilde{c}(t_n+t)}$. We have $\varepsilon^n(-t_n)=\widetilde{\varepsilon}_0$ by Remark~\ref{remarque: widetildea(0)=0 and widetildec(0)=c_0^*}. Moreover, by~\eqref{normXvarepsilon(t+normRc(t)-c^*+normRa(t)-a^leq A_*alpha_0}, we have $\normX{\varepsilon^n(t)}\leq\alpha_*$ for any $t\in\R$. Applying~\eqref{normXvarepsilon(t+normRc(t)-c^*+normRa(t)-a^leq A_*alpha_0} once more, we can make $\normX{\varepsilontilde_0}$ as small as $\normX{\varepsilontilde(t_n)}$ can be. Thus, $\widetilde{\varepsilon}_0=0$ and in particular, $\widetilde{Q}_0=Q_{c^*_0}$. This necessarily leads to $\widetilde{Q}(t)=Q_{c^*_0}(.-c_0^* t)$ and then for any $t\in\R,\varepsilontilde(t)=0$. Finally, we derive from~\eqref{controle des parametres de modulation} and Remark~\ref{remarque: widetildea(0)=0 and widetildec(0)=c_0^*} that $\widetilde{c}(t)=c_0^*$ and $\widetilde{a}(t)=c_0^* t$.
\end{proof}

\subsubsection{Algebraic decay of the solution $\widetilde{Q}$}\label{section: Convergence along the evolution trightarrow +ii}

In this section, we verify that $\widetilde{Q}$ satisfies the assumptions of Proposition~\ref{lemme: rigidity property proche d'un soliton, version algébrique}. Furthermore, we shall exhibit the suitable value of $l_*$ for Theorem~\ref{théorème: stab asymp classical} and eventually conclude the proof of Theorem~\ref{théorème: stab asymp hydro}. To do this, we rely on two propositions stated below. The first one is an almost monotonicity formula for a localized version of the momentum. The second proposition states a remarkable smoothing property satisfied by the solutions of~\eqref{NLS}. Both Propositions~\ref{proposition: monotonie formula} and~\ref{proposition: controle L^2_loc L^2_x partial_x^l psi} are proven in Section~\ref{section: monotonie et smoothing effect} and depend crucially on the smoothness of $f$.

Recall first that the momentum is conserved along the flow, and let us introduce the following localized version for $\widetilde{Q}$
\begin{equation}\label{définition localized momentum}
    \widetilde{p}_R(t)=\int_\R \pi\big(\widetilde{Q}(t)\big)\chi_{R+\widetilde{a}(t)}=\dfrac{1}{2}\int_\R \etatilde(t)\widetilde{v}(t)\chi_{R+\widetilde{a}(t)},
\end{equation}
where $\chi_R(x)=\chi\big(x-R)$ and
\begin{equation*}
    \chi(x)=\dfrac{1}{2}\Big( 1 + \tanh\left(\frac{\tau_*}{2} x\right)\Big).
\end{equation*}

Here, the value $\tau_*$ is to be determined in the proof of Proposition~\ref{proposition: monotonie formula} and depends only on $c^*$. The quantity~\eqref{définition localized momentum} measures the amount of momentum of $\widetilde{Q}\big(t,.+\widetilde{a}(t)\big)$ on the right of a given real number $R\in\R$. The evolution of this quantity is prescribed according to the next proposition. This evolution also holds for the localized momentum associated with the original solution $Q$ (see Section~\ref{section: monotonie et smoothing effect})\footnote{
Throughout the article, we shall write $g(.,c)=\grandOde{f(.,c)}$ or $g(.,c)\lesssim f(.,c)$ whenever there exists a constant $K>0$ only depending on $c^*$ such that for any $(x,c)\in\R\times (0,c_s),\left| f(x,c)\right|\leq K \left|g(x,c)\right|$.
}.
\begin{prop}\label{proposition: monotonie formula}
There exist positive constants $\kappa_*,\tau_*,\sigma_*$ only depending on $c_*$ such that for any $t\in\R$ and $\sigma\in [-\sigma_*,\sigma_*]$,
\begin{equation*}
    \dfrac{d}{dt}\Big(\widetilde{p}_{R+\sigma t}(t)\Big)\geq \kappa_*\int_\R \Big((\partial_x\widetilde{\eta})^2 + \etatilde^2 + \widetilde{v}^2\Big)\big(t,.+\widetilde{a}(t)\big)\chi'(.-R-\sigma t) + \grandOde{e^{-\tau_* |R+\sigma t|}}.
\end{equation*}
\end{prop}

Arguing exactly as in~\cite{BetGrSm2}, we deduce the next corollary.

\begin{cor}\label{corollaire: controle du moment localisé en |x|^delta}
   For any $\rho \geq 0$, there exists $\kappa_{\rho}$ such that for any $t\in\R$, we have
\begin{equation*}
    \int_{t}^{t+1}\int_\R \Big((\partial_x\widetilde{\eta})^2 + \etatilde^2 + \widetilde{v}^2\Big)\big(s,x+\widetilde{a}(s)\big)|x|^\rho ds \leq \kappa_{\rho}.
\end{equation*} 
\end{cor}

Corollary~\ref{corollaire: controle du moment localisé en |x|^delta} holds for any $\rho\geq 0$, so can successively be applied with any $\rho\in [0,r]$ for some fixed $r>\frac{5}{2}$. It provides a uniform $\kappa:= \max_{\rho\in [0,r]}\kappa_\rho$ that shall be useful to control the $L^2$-norm of $1-|\widetilde{\psi}|^2$ (see~\eqref{justification de choix de kappa}). Next, we deduce that the solution $\widetilde{\psi}$ of~\eqref{NLS} associated with the variables $(\widetilde{\eta},\widetilde{v})$ can be differentiated a certain number of times (depending on $l_*$) and enjoys polynomial decay at infinity. This decay follows from some smoothing effect that was already involved in~\cite{BetGrSm2} (see also~\cite{EsKePoV5} for related results). Henceforth, we impose the condition
\begin{equation*}
    l_*\geq l_0+1,
\end{equation*}
with $l_0$ given by Proposition~\ref{lemme: rigidity property proche d'un soliton, version algébrique} and verify that the following proposition is sufficient for the Liouville property to apply to $\widetilde{Q}$.

\begin{prop}\label{proposition: controle L^2_loc L^2_x partial_x^l psi}
Let $f\in C^{l_*}(\R_+)$. For any $l\in\{1,...,l_*\}$, there exists $C_{l}>0$ only depending on $c^*$ such that for any $\Lambda\in\{0,1\},t\in\R$ and any $x\in\R$, 
\begin{equation*}
    \int_{t}^{t+1}\int_\R \Big|\partial_s^\Lambda\big[\partial_x^l \widetilde{\psi}\big(s,x+\widetilde{a}(s)\big)\big]\Big|^2(1+|x|^{r}) dxds\leq C_{l}.
\end{equation*}
\end{prop}

Having admitted that Proposition~\ref{proposition: controle L^2_loc L^2_x partial_x^l psi} holds, we complete the proof of Theorem~\ref{théorème: stab asymp hydro}. Set $\widetilde{\Psi}(t,x)=\widetilde{\psi}\big(t,x+\widetilde{a}(t)\big)$. Assuming the preceding proposition, we can show that $\widetilde{\psi}$ possesses the suitable decay. We first prove that for any $(t,x)\in\R^2$ and any $l\in\{1,...,l_*\}$, and up to a larger constant $C_{l}$, we have
\begin{equation}\label{big|partial_x^l widetildepsi(t,xwidetildea(t))big|leqdfracC_*1+|x|^rho}
    \big|\partial_x^l \widetilde{\Psi}(t,x)\big|\leq\dfrac{C_{l}}{1+|x|^r}.
\end{equation}

Indeed, we set $\widetilde{\Psi}_l:=\sqrt{1+|x|^r}\ \widetilde{\Psi}$ and by Proposition~\ref{proposition: controle L^2_loc L^2_x partial_x^l psi}, we can check that for any $l\in\{1,...,l_*-1\},\widetilde{\Psi}_l\in H^1\big(I,H^1(\R)\big)$ with $I=[t-1,t+1]$. More precisely, using the one-dimensional Sobolev embedding $H^1\hookrightarrow L^\ii$ simultaneously in time and in space, we obtain for any $l\in\{1,...,l_*-1\},t\in\R$ and $x\in\R$,
\begin{equation*}
    (1+|x|^{\widetilde{r}}) \left|\partial_x^l\widetilde{\Psi}\right|^2\leq \left\Vert \widetilde{\Psi}_l\right\Vert_{L^\ii\big(I,L^\ii(\R)\big)}^2\lesssim\left\Vert\widetilde{\Psi}_l\right\Vert_{H^1\big(I,H^1(\R)\big)}^2\leq 8C_{l},
\end{equation*}

This implies~\eqref{big|partial_x^l widetildepsi(t,xwidetildea(t))big|leqdfracC_*1+|x|^rho}. Next, we recall that we can recover $v$ by the formula
\begin{equation*}
    \widetilde{v}=\dfrac{\partial_x\widetilde{\psi}.i\widetilde{\psi}}{|\widetilde{\psi}|^2},
\end{equation*}

where $.$ denotes the usual real scalar product defined by $z_1.z_2 = \mathrm{Re}(z_1\overline{z_2}), \forall z_1,z_2\in\C$. Hence by using~\eqref{big|partial_x^l widetildepsi(t,xwidetildea(t))big|leqdfracC_*1+|x|^rho} and Proposition~\ref{proposition: borne uniforme sur les quantités nu_c etc} we show recursively that for any $l\in\{0,...,l_*\}$, $(t,x)\in\R^2$ and up to taking a larger constant $C_{l}$,
\begin{equation}\label{big|partial_x^l widetildev(t,xwidetildea(t))big|leqdfracC_*1+|x|^rho}
    \big|\partial_x^l\widetilde{v}\big(t,x+\widetilde{a}(t)\big)\big|\leq \dfrac{ C_{l}}{1+|x|^r}.
\end{equation}

On the other hand, from~\eqref{big|partial_x^l widetildepsi(t,xwidetildea(t))big|leqdfracC_*1+|x|^rho}, we also derive
\begin{equation}\label{big|partial_x^l widetildeeta(t,xwidetildea(t))big|leqdfracC_*1+|x|^rho}
    \big|\partial_x^l\etatilde\big(t,x+\widetilde{a}(t)\big)\Big|\leq \dfrac{C_{l}}{1+|x|^r},
\end{equation}
 for any $l\in\{0,...,l_*\}$. Let us handle the cases $l\in\{1,...,l_*\}$ first, writing $\etatilde=1-\widetilde{\psi}.\widetilde{\psi}$, we show using~\eqref{big|partial_x^l widetildepsi(t,xwidetildea(t))big|leqdfracC_*1+|x|^rho} and Proposition~\ref{proposition: borne uniforme sur les quantités nu_c etc} that 
\begin{equation*}
    \big|\partial_x^l\etatilde\big(t,x+\widetilde{a}(t)\big)\Big|\leq \sum_{m=0}^l \begin{pmatrix}
        l \\
        m
    \end{pmatrix}\big|\partial_x^m\widetilde{\Psi}(t,x)\big|\big|\partial_x^{l-m}\widetilde{\Psi}(t,x)\big|\leq \dfrac{C_l}{(1+|x|^r)^2} .
\end{equation*}

Then it remains to deal with the case $l=0$. By the one-dimensional Sobolev embedding, it is sufficient to prove that \begin{equation*}
    \left\Vert 1-|\widetilde{\Psi}|^2\right\Vert_{H^1\big(I,L^\ii(|x|^r dx)\big)}\leq C_{l} .
\end{equation*}

By Corollary~\ref{corollaire: controle du moment localisé en |x|^delta}, and the special choice of $\kappa$, we derive \begin{equation}\label{justification de choix de kappa}
    \left\Vert 1-|\widetilde{\Psi}|^2\right\Vert_{L^2\big(I,L^\ii(|x|^r dx)\big)}\leq \kappa,
\end{equation}

and by~\eqref{NLShydro} we have
\begin{equation*}
    \partial_t\big(1-|\widetilde{\Psi}|^2)(t,x)=\widetilde{a}'(t)\partial_x\etatilde\big(t,x+\widetilde{a}(t)\big)-2\partial_x\big(\widetilde{v}(1-\etatilde)\big)(t,x).
\end{equation*} 
Therefore, by~\eqref{big|partial_x^l widetildev(t,xwidetildea(t))big|leqdfracC_*1+|x|^rho},~\eqref{big|partial_x^l widetildeeta(t,xwidetildea(t))big|leqdfracC_*1+|x|^rho} and once more Proposition~\ref{proposition: borne uniforme sur les quantités nu_c etc}, we obtain
\begin{equation*}
    \left\Vert \partial_t (1-|\widetilde{\Psi}|^2)\right\Vert_{L^2\big(I,L^\ii(|x|^r dx)\big)}\leq \kappa .
\end{equation*}

\underline{Conclusion.} Provided that both Propositions~\ref{proposition: monotonie formula} and~\ref{proposition: controle L^2_loc L^2_x partial_x^l psi} hold and since $l_*\geq l_0+1$, we have proven that~\eqref{lemme: decaying property in rigidity} is satisfied by the limit profile $\widetilde{Q}$. Then the Liouville property in Proposition~\ref{lemme: rigidity property proche d'un soliton, version algébrique} applies and for any $t\in\R$ 
\begin{equation*}
    \widetilde{Q}(t)=Q_{c^*_0}(.-c_0^* t).
\end{equation*}
This yields, by~\eqref{convergence faible vers le profil limite}, to
\begin{equation*}
        Q\big(t_n,.+a(t_n)\big)\ntendfX Q_{c^*_0}.
    \end{equation*}

It remains to see that the convergence does not depend on the choice of the subsequence. This also relies on the monotonicity property in Proposition~\ref{proposition: monotonie formula} applied to $Q$. More precisely, it is a consequence of~\eqref{limite widetilde p_R vers 0} in Section~\ref{section: monotonie et smoothing effect}. In particular, it can be proved the same way as it was already done in Subsection~1.3.4 in~\cite{BetGrSm2}, given that this part does not depend on the nonlinearity $f$. This concludes the proof of Theorem~\ref{théorème: stab asymp hydro}.

\subsubsection{Proof of the asymptotic stability in the classical variables}
We have just proved the asymptotic stability of a hydrodynamical soliton $Q_{c}=(\eta_{c},v_{c})$ for ${c}\in (c_0,c_s)$. Let us finish the proof and show that the travelling wave $\gu_{c}$ is asymptotically stable in the original setting. Invoking Lemma~\ref{lemme: équivalence des distances classique et hydro proche d'un soliton}, we fix $\delta_c$ small enough such that $d(\psi_0,\gu_c)\leq \delta_c$ and $\normX{Q_0-Q_{c}}\leq \beta_{c}$. Therefore Theorem~\ref{théorème: stab asymp hydro} applies and there exist $c^*\in (c_0,c_s)$ and a function $a$ such that whenever $t_n\rightarrow +\ii$, $Q_{t_n}=(\eta_{t_n},v_{t_n}):=Q\big(t_n,.+a(t_n)\big)$ satisfies
\begin{equation*}\label{eta_t_nntendf eta_gcquadtextin L^2(R)}
    Q_{t_n}\ntendfX Q_{c^*}.
\end{equation*}

In view of the previous convergence, we claim that

\begin{equation}\label{psibig(t,.+b(t)big)tiitendfd gu_c^*}
    \psi\big(t,.+a(t)\big)\tiitendfd \gu_{c^*}.
\end{equation}

We now verify that all three convergences in~\eqref{psibig(t,.+b(t)big)tiitendfd gu_c^*} are true. Setting the notation $\psi_{t}:=\psi(.+a(t),t)$, we have that $(\psi_{t_n})_n$ is bounded in $H^1_{\loc}(\R)$. Then up to a subsequence, we get that

\begin{equation}\label{Psi_t_n tendf vers Psi ii dans H^1_loc}
\psi_{t_n}\ntendf \psi_\ii\quad\text{in }H^1_\loc(\R),
\end{equation}

for some function $\psi_\ii\in H^1_{\loc}(\R)$. By the Rellich theorem and up to taking a further subsequence, convergence~\eqref{Psi_t_n tendf vers Psi ii dans H^1_loc} holds strongly in $L^\ii_\loc(\R)$ as well, thus we have the pointwise equality of the functions $|\psi_\ii|$ and $|\gu_{c^*}|$. In particular, we have the decomposition $\psi_\ii=|\psi_\ii|e^{i\varphi_\ii}$. However, the hydrodynamic framework omits the phase of the classical variable $\psi$, so that we need to modulate while taking care of this degree of freedom. To take this into account, we shall construct a time-dependent phase shift $\theta$ such that the underlying solution converges as desired. This modulates the evolution so that the phases of both $\psi_\ii$ and $\gu_{c^*}$ eventually match for any time $t$.

\begin{claim}\label{claim: régularity of theta}
    There exist a function $\theta\in C^1(\R,\R/2\pi\Z)$ such that $\theta'$ is bounded and a bump function $\chi\in C^\ii_0(\R)$ such that
    \begin{equation*}
        \gd_{c^*}:=\int_\R\chi\gu_{c^*}\neq 0,
    \end{equation*}
    and for any $t\in\R$, 
    \begin{equation*}
        e^{-i\theta(t)}\int_\R \chi(x)\psi(x+a(t),t)dx \in \gd_{c^*}\R_+^*.
        \end{equation*}
\end{claim}

Replacing~\eqref{Psi_t_n tendf vers Psi ii dans H^1_loc} by $e^{-i\theta(t_n)}\psi_{t_n}$ does not change the pointwise equality $|\psi_\ii|=|\gu_c|$ so that we can assume that
\begin{equation*}\label{e^it_n Psi_t_n tendf vers Psi ii dans H^1_loc}
e^{-i\theta(t_n)}\psi_{t_n}\ntendf \psi_\ii\quad\text{in }H^1_\loc(\R).
\end{equation*}

On the other hand we have  
\begin{equation*}
    v_{t_n}\ntendf v_{c^*}\quad\text{in }L^2(\R),
\end{equation*}

while
\begin{equation*}
    v_{t_n}=\dfrac{i\psi_{t_n}.\partial_x\psi_{t_n}}{|\psi_{t_n}|^2}\ntendf \dfrac{i\psi_{\ii}.\partial_x\psi_{\ii}}{|\psi_{\ii}|^2}=\varphi'_\ii\quad\text{in }L^2_{\loc}(\R).
\end{equation*}

Then by integrating between $0$ and $x$, we obtain for any $x\in\R$ that $\varphi_{c^*}(x)-\varphi_{c^*}(0)=\varphi_\ii(x)-\varphi_\ii(0)$ and therefore the expression $\psi_\ii=\gu_{c^*} e^{i\big(\varphi_\ii(0)-\varphi_{c^*}(0)\big)}$. Hence, the convergence 

\begin{equation*}
    e^{-i\theta(t_n)}\int_\R \psi_{t_n}\chi \ntend e^{i(\varphi_\ii(0)-\varphi_{c^*}(0))}\int_\R \gu_{c^*}\chi. 
\end{equation*}

Now Claim~\ref{claim: régularity of theta} forces the difference of phases to behave as follows,
\begin{equation*}
    \varphi_\ii(0)=\varphi_{c^*}(0)\quad\text{in }\R/2\pi\Z .
\end{equation*}

In summary, we have proved that
\begin{equation*}\label{e-itheta(t_n)Psi(.+b(t_n),t_n)ntendfd gu_gc}
    e^{-i\theta(t_n)}\psi(.+a(t_n),t_n)\ntendfd \gu_{c^*} .
\end{equation*}

In view of what has been done in the hydrodynamical variable, we can show that the convergence does not depend on the choice of the sequence $(t_n)$ and thus the convergence holds as $t\rightarrow +\ii$ and one can write 
\begin{equation*}\label{e-itheta(t)Psi(.+b(t),t) ntendfd gu_gc}
    e^{-i\theta(t)}\psi(.+a(t),t)\tiitendfd \gu_{c^*} .
\end{equation*}

Furthermore, this convergence still holds along the variable $t+T$ as $t\rightarrow +\ii$. This gives
\begin{equation}\label{e-itheta(t+T)Psi(.+b(t+T),t+T) ntendfd gu_gc}
    e^{-i\theta(t+T)}\psi(.+a(t+T),t+T)\tiitendfd \gu_{c^*} .
\end{equation}

From the weak continuity of the flow in standard variables in Proposition~\ref{proposition: weak continuity of the flow in original setting}, we deduce that for any $T\in\R$, 
\begin{equation}\label{e-itheta(t)Psi(.+b(t),t+T) ntendfd gu_gc}
    e^{-i\theta(t)}\psi(.+a(t),t+T)\tiitendfd \gu_{c^*}(.-{c^*} T) .
\end{equation}

On the other hand, by Proposition~\ref{proposition: convergence faible vers le profil limite le long de l'évolution}, we infer that \begin{equation*}\label{b(t+T)-b(t)ttendii gc T}
    a(t+T)-a(t)\ttendii {c^*} T.
\end{equation*}
Convolving $t\mapsto a(t)$ with any fixed mollifier provides the suitable function $b$ satisfying for any $T\in\R$
\begin{equation}\label{b'(t) tend vers c^* et}
    b(t)-a(t)\ttendii 0,\quad b'(t)\ttendii c^*\quad\text{and}\quad b(t+T)-b(t)\ttendii c^* T.
\end{equation}

Now, by the triangle inequality, we have 
\begin{align}
    \left|\left(e^{i\big(\theta(t)-\theta(t+T)\big)}-1\right)\gu_{c^*}\right|&\leq \left|\gu_{c^*}-e^{-i\theta(t)}\psi(.+b(t)+{c^*} T,t+T)\right|\label{1ère inégalité dans e ibig(theta(t)-theta(t+T)big)}\\
    &+\Big| \psi(.+b(t)+{c^*} T,t+T)-\psi(.+b(t+T),t+T)\Big|\label{2eme inégalité dans e ibig(theta(t)-theta(t+T)big)}\\
    &+\left|e^{-i\theta(t+T)}\psi(.+b(t+T),t+T)-\gu_{c^*}\right|\label{3eme inégalité dans e ibig(theta(t)-theta(t+T)big)}.
\end{align}

By~\eqref{e-itheta(t)Psi(.+b(t),t+T) ntendfd gu_gc} resp.~\eqref{e-itheta(t+T)Psi(.+b(t+T),t+T) ntendfd gu_gc}, we deduce that the right-hand terms in~\eqref{1ère inégalité dans e ibig(theta(t)-theta(t+T)big)} resp.~\eqref{3eme inégalité dans e ibig(theta(t)-theta(t+T)big)} converge towards $0$ uniformly on any compact set of $\R$ as $t$ tends to $+\ii$. As for the term~\eqref{2eme inégalité dans e ibig(theta(t)-theta(t+T)big)}, it converges weakly to $0$ in $H^1(\R)$ due to the right-hand side of~\eqref{b'(t) tend vers c^* et}. This yields 
\begin{equation*}
    e^{i\big(\theta(t)-\theta(t+T)\big)}\gu_{c^*}\ttendii \gu_{c^*}\quad\text{in }L^\ii_\loc(\R).
\end{equation*}

Now, take a representative of the equivalence class in the torus, still denoted by $\theta$. By the previous convergence, we necessarily have
\begin{equation*}
    \theta(t+T)-\theta(t)\ttendii 2\pi l_T,
\end{equation*}
for some $l_T\in\Z$. Let us verify that $l_T=0$. By Claim~\ref{claim: régularity of theta}, we can first write $\big|\theta(t+T)-\theta(t)\big|\leq \normLii{\theta'}T$. Provided that $T$ is small enough so that $\normLii{\theta'}T < 2\pi$, we deduce by taking the limit $t\rightarrow +\ii$ that $l_T=0$. We can now adapt this argument for any $T\in\R$, by taking $\widetilde{T}=\frac{T}{n}$ with an integer $n$ large enough, so that $l_{\widetilde{T}}=0$. We conclude by writing 
\begin{equation*}
    \big|\theta(t+T)-\theta(t)\big|\leq \sum_{l=1}^{n}\big|\theta\big(t+(l-1)\widetilde{T}+\widetilde{T}\big)-\theta\big(t+(l-1)\widetilde{T}\big)\big|.
\end{equation*}
and by passing to the limit $t+(l-1)\widetilde{T}\rightarrow +\ii$ as $t\rightarrow +\ii$. 
We can finally convolve such a representative $\theta$ with a fixed mollifier as for constructing the function $b$ above, and keeping the notation $\theta$ for such a function, we therefore obtain a modified $C^1$ function $\theta$ satisfying
\begin{equation*}
    \theta'(t)\ttendii 0,
\end{equation*}

which concludes the proof of Theorem~\ref{théorème: stab asymp classical}. It remains to prove Claim~\ref{claim: régularity of theta}.

\begin{proof}[Proof of Claim~\ref{claim: régularity of theta}]
The function $\gu_{c^*}$ does not vanish on $\R$, so that in particular, $\gu_{c^*}(0)\neq 0$. Given any smooth function $\chi$, compactly supported in $[-1,1]$ such that $\chi(0)=1$, and up to shrinking the length of its support by taking $\chi_n=\chi(\frac{.}{n})$ and $n$ large enough, we can assume that $\int_\R \chi\gu_{c^*}\neq 0$. To construct $\theta$ for any $t\in\R$, we first write that for any $\widetilde{\varphi}\in\R$ that does not depend on $x$,
\begin{align*}
    \left|\int_\R\chi\psi_t\right|&\geq \left|\int_\R e^{i\widetilde{\varphi}}\chi\gu_{c^*}\right|-\left|\int_\R \chi\left(\psi_t-e^{i\widetilde{\varphi}}\gu_{c^*}\right)\right|\\
    &\geq |\gd_{c^*}|-\normLii{\psi_t-e^{i\widetilde{\varphi}}\gu_{c^*}}.
\end{align*}

Lifting the non-vanishing solution $\psi$ as $\psi=|\psi|e^{i\varphi}$, setting $\varphi_t:= \varphi(.+b(t),t)$ and $v_t:=-\partial_x\varphi_t$, and by the special choice of adding a phase shift $\widetilde{\varphi}=\varphi(b(t),t)-\varphi_{c^*}(0)$ (that only depends on $t$), we then derive \begin{align*}
    \normLii{\psi_t-e^{i\widetilde{\varphi}}\gu_{c^*}}& = \normLii{|\psi_t |e^{i\varphi_t}-|\gu_{c^*}|e^{i\varphi_t}}+\normLii{|\gu_{c^*}|(e^{i\varphi_t}-e^{i(\varphi_{c^*}+\widetilde{\varphi})})}\\
    &\leq \normLii{\big|\psi_t\big|-\big|\gu_{c^*}|}+ \normLii{\left(e^{i\varphi_t}-e^{i\big(\varphi_{c^*}-\varphi_{c^*}(0)+\varphi(b(t),t)\big)}\right)|\gu_{c^*}|}\\
    & \leq \normLii{\dfrac{\big|\psi_t\big|^2-\big|\gu_{c^*}|^2}{\big|\psi_t\big|+\big|\gu_{c^*}|}}+\normLii{\left(\varphi_t-\big(\varphi_{c^*}-\varphi_{c^*}(0)+\varphi(b(t),t)\big)\right)|\gu_{c^*}|}\\
    &\leq \dfrac{1}{\inf_\R|\gu_{c^*}|}\normLii{\eta_t-\eta_{c^*}}+\normLdeux{v_t-v_{c^*}}\normLii{\sqrt{|x|}\gu_{c^*}}.\\
\end{align*}

By orbital stability, provided that $\beta_{c}$ is taken small enough, the norm $\normX{Q_t-Q_{c}}$ can be taken as small as we wish, uniformly in time. Therefore, by both previous estimates, we can assume that for any $t\in\R$, \begin{equation}\label{int_R chi psi_t geq}
    \left|\int_\R\chi\psi_t\right|\geq \dfrac{|\gd_{c^*}|}{2} >0.
\end{equation}

Thus we can lift the previous quantity and there exists a unique function $\theta:\R\rightarrow\R/2\pi\Z$ such that for any $t\in\R$, 
\begin{equation}\label{e^-itheta int chi psi_t}
    e^{-i\theta(t)}\int_\R\chi\psi_t\in\gd_{c^*}\R^*_+.
\end{equation}

Regarding the smoothness of $\theta$, we implement an implicit function type argument. Consider the map $X$ defined for any $(t,\theta)\in\R^2$ by
\begin{equation*}
    X(t,\theta)=\mathfrak{Im}\left(e^{-i\arg(\gd_{c^*})}e^{-i\theta}\int_\R\chi\psi_t\right),
\end{equation*}
where $\arg(\gd_{c^*})$ designates the principal argument of $\gd_{c^*}$. By construction, $X\big(t,\theta(t)\big)=0$ for any $t\in\R$. Moreover, $\partial_2 X\big(t,\theta(t)\big)\neq 0$, then by uniqueness, we derive that $\theta$ is smooth and that $\theta$ satisfies the ordinary differential equation
\begin{equation*}
    \theta'(t)=\dfrac{-\mathfrak{Im}\left(ie^{-i\arg(\gd_{c^*})}e^{-i\theta(t)}\displaystyle\int_\R\big(\chi''+ib'(t)\chi'+f(|\psi_t|^2)\chi\big)\psi_t\right)}{X\big(t,\theta(t)+\frac{\pi}{2})\big)}.
\end{equation*}

Applying the Sobolev embedding theorem to $\psi_t$, and using~\eqref{int_R chi psi_t geq} and the definition of $\theta(t)$, we observe that this derivative is uniformly bounded in time.
\end{proof}

\numberwithin{equation}{section}

\section{Monotonicity and smoothing properties of~\eqref{NLS}}\label{section: monotonie et smoothing effect}
First, we prove the monotonicity formula.
\begin{proof}[Proof of Proposition~\ref{proposition: monotonie formula}]
We recall from~\cite{Bert2} the following lemma.
\begin{lem}\label{lemme: d/dt ( int psi eta v)}
Let $(\eta,v)\in C^0([0,T],\Nenergysethydro)$ be a solution of~\eqref{NLShydro}. For any function $\widetilde{\chi}\in C^0\big([0,T], C^3_b(\R)\big)\cap C^1\big([0,T], C^0_b(\R)\big)$, then $t\mapsto\psLdeux{\widetilde{\chi}\eta}{ v}$ is differentiable and its derivative is
\begin{align}
    \dfrac{d}{dt}\left(\int_\R \widetilde{\chi}\eta v\right)=&\int_\R \partial_t\widetilde{\chi} \eta v+ \int_\R \partial_x\widetilde{\chi}\Big((1-2\eta)v^2 + \widetilde{F}(\eta)+\dfrac{(3-2\eta)(\partial_x\eta)^2}{4(1-\eta)^2}\Big) \notag \\ 
    & + \dfrac{1}{2}\int_\R \partial_x^3\widetilde{\chi} \big(\eta + \ln(1-\eta)\big)\label{d/dt ( int psi eta v)},
\end{align}
where \begin{equation}\label{définition de widetilde F}
    \widetilde{F}(\rho)=\rho f(1-\rho)-F(1-\rho).
\end{equation}
\end{lem}

Applying this to $\widetilde{\chi}(t,x)=\chi_{R+\sigma t}(x)$ provides
\begin{align*}
    \dfrac{d}{dt}\left(\widetilde{p}_{R+\sigma t}(t)\right) = &\int_\R \pi_{R+\sigma t}(\etatilde,\widetilde{v})\big(t,.+\widetilde{a}(t)\big),
\end{align*}

where
\begin{align*}
    \pi_{R+\sigma t}(\etatilde,\widetilde{v})=\dfrac{1}{2}\chi'_{R+\sigma t}\bigg((1-2\etatilde)\widetilde{v}^2 + \widetilde{F}(\etatilde)+\dfrac{(3-2\etatilde)(\partial_x\etatilde)^2}{4(1-\etatilde)^2} -&\big(\widetilde{a}'(t)+\sigma\big) \etatilde\widetilde{v}\bigg)  \\
    &+\dfrac{1}{4}\chi'''_{R+\sigma t}\big(\etatilde+\ln(1-\etatilde)\big).
\end{align*}

On the one hand, we have by~\eqref{controle des parametres de modulation}
\begin{equation*}
    \big|\widetilde{a}'(t)+\sigma\big|\leq |\widetilde{a}'(t)-c_0^*|+c_0^* + \sigma_*\leq \mu_*,
\end{equation*} 
where $\mu_*:=A_*\alpha_*+c_0^*+\sigma_*$, and where the values of $\alpha_*$ and $\sigma_*$ can be decreased so that $\mu_*<c_s$. In addition, we have $|\chi'''_R|\leq \tau_*^2\chi_R'$ and by~\eqref{définition de widetilde F} and some standard real analysis, there exists a positive constant $C_{\ln}$ such that
\begin{equation*}
    \dfrac{1}{2}\chi_R' \widetilde{F}(\etatilde) + \dfrac{1}{4}\chi'''_{R}\big(\etatilde+\ln(1-\etatilde)\big) \geq \dfrac{1}{2}\chi_R'\left(-\int_0^1 rf'(1-r\etatilde)dr -\tau_*^2 C_{\ln}\right)\etatilde^2.
\end{equation*}

As a consequence of the previous work and Proposition~\ref{proposition: borne uniforme sur les quantités nu_c etc}, we obtain
\begin{equation*}
    \pi_R(\etatilde,\widetilde{v}) \geq \dfrac{1}{2}\chi_R'\Big(\dfrac{(3-2I_*)(\partial_x\etatilde)^2}{8(1+I_*)^2}+  q(\etatilde,\widetilde{v})\Big),
\end{equation*}

where \begin{align*}
    q(\etatilde,\widetilde{v})&=(1-2\etatilde)\widetilde{v}^2-\mu_*|\etatilde\widetilde{v}| +\left(-\int_0^1 rf'(1-r\etatilde)dr -\tau_*^2 C_{\ln}\right)\etatilde^2 .
\end{align*}

We now separate the real line into two disjoint parts that are an interval $\mathcal{I}=[-R_0,R_0]$, where the travelling wave is small, and its complement $\mathcal{I}^c$ where $\chi'$ is small. Regarding the region where $x\in \mathcal{I}^c$, we state the following claim.

\begin{claim}\label{coercivité dans la forme quad dans la monotonie du moment}
Up to taking $R_0$ large enough, and possibly shrinking the value of $\tau_*$, there exists $\kappa_2>0$ such that on $\mathcal{I}^c$,
    \begin{equation*}
        q(\etatilde,\widetilde{v})\geq \kappa_2 \big(\etatilde^2+\widetilde{v}^2\big).
    \end{equation*}
\end{claim}

Assuming that the preceding claim holds true, we conclude that, up to shrinking the value of $\kappa_2$, that on $\mathcal{I}^c$ and for $\sigma\in [-\sigma_*,\sigma_*]$,
\begin{equation*}
    \pi_{R+\sigma t}(\etatilde,\widetilde{v}) \geq \kappa_2\chi_{R+\sigma t}'\Big( (\partial_x\etatilde)^2 + \etatilde^2+ \widetilde{v}^2\Big).
\end{equation*}

Therefore, 
\begin{align*}
    \dfrac{d}{dt}\left(\widetilde{p}_{R+\sigma t}(t)\right) \geq \int_\mathcal{I} &\pi_{R+\sigma t}(\etatilde,\widetilde{v}) + \kappa_2\int_{\mathcal{I}^c} \Big((\partial_x\etatilde)^2+\etatilde^2 + \widetilde{v}^2\Big)\big(t,.+\widetilde{a}(t)\big)\chi'_{R+\sigma t}\\
    =\int_\mathcal{I} &\Big(\pi_{R+\sigma t}(\etatilde,\widetilde{v}) - \kappa_2\big((\partial_x\etatilde)^2+\etatilde^2 + \widetilde{v}^2\big)\chi'_{R+\sigma t}(x)\Big) \big(t,x+\widetilde{a}(t)\big)dx\\
    &+\kappa_2\int_{\R} \Big((\partial_x\etatilde)^2+\etatilde^2 + \widetilde{v}^2\Big)\big(t,.+\widetilde{a}(t)\big)\chi'_{R+\sigma t}.
\end{align*}

On the other hand, and from the previous analysis and the Cauchy-Schwarz inequality, we have for any $x\in\R$, 
\begin{align*}
    \big|\pi_{R + \sigma t}(\etatilde,\widetilde{v}) \big|\leq \dfrac{1}{2}\chi'_{R+\sigma t}\Big( \big(1+2I_*+\frac{\mu_*}{2}\big)\widetilde{v}^2 + \big(\Vert f''\Vert_{L^\ii([0,2])}+C_{ln}+\frac{\mu_*}{2}\big)\etatilde^2+\dfrac{3+2I_*}{4 \iota_*^2}\tau_*^2(\partial_x\etatilde)^2\Big).
\end{align*}

Furthermore, using~\eqref{hypothèse de croissance sur F minorant intermediaire} and Proposition~\ref{proposition: borne uniforme sur les quantités nu_c etc}, we have the useful bounds
\begin{equation}\label{controle des termes du moment par l'énergie}
    (\partial_x\etatilde)^2\leq 8(1+I_*)\dfrac{(\partial_x\etatilde)^2}{8(1-\etatilde)},\quad \etatilde^2\leq \dfrac{4}{c_s^2}\dfrac{F(1-\etatilde)}{2}\quad\text{and}\quad \widetilde{v}^2\leq \dfrac{(1-\etatilde)\widetilde{v}^2}{\iota_*},
\end{equation}

so that there exists a positive constant $\kappa_1$ such that on $\R$ now,
\begin{equation*}
    \big|\pi_{R + \sigma t}(\etatilde,\widetilde{v}) - \kappa_2\big((\partial_x\etatilde)^2+\etatilde^2 + \widetilde{v}^2\big)\chi'_{R+\sigma t}\big|\leq \kappa_1 e(\etatilde,\widetilde{v})\big|\chi'_{R+\sigma t}\big|,
\end{equation*}

where $e(\etatilde,\widetilde{v})$ is the energy density defined in~\eqref{expression de l'énergie}. Now, since we have for any $x\in \mathcal{I}$,
\begin{equation*}
    \chi'_{R+\sigma t}(x)\leq \tau_* e^{-\tau_*|x-(R+\sigma t)|}\leq \tau_* e^{\tau_* R_0}e^{-\tau_*|R+\sigma t|}=\grandOde{e^{-\tau_*|R+\sigma t|}},
\end{equation*}

we obtain the control
\begin{align*}
    \dfrac{d}{dt}\big(\widetilde{p}_{R+\sigma t}(t)\big)\geq -\kappa_1 E(\etatilde,\widetilde{v})e^{-\tau_*|R+\sigma t|} + \kappa_2\int_\R \Big((\partial_x\etatilde)^2+\etatilde^2 + \widetilde{v}^2\Big)\big(t,.+\widetilde{a}(t)\big)\chi'_{R+\sigma t},
\end{align*}
which concludes the proof by conservation of the energy.
\end{proof}

It remains to prove the claim.
\begin{proof}[Proof of Claim~\ref{coercivité dans la forme quad dans la monotonie du moment}]

We write \begin{align*}
    q(\etatilde,\widetilde{v})&=A_1 \etatilde^2 - A_2|\etatilde\widetilde{v}| + A_3 \widetilde{v}^2,
\end{align*}
with 
\begin{align*}
    A_1 &= -\int_0^1 rf'(1-r\etatilde)dr -\tau_*^2 C_{\ln},\\
    A_2&=\mu_* ,\\
    A_3&=1-2\etatilde,
\end{align*}

and obtain the expression
\begin{equation*}
    q(\etatilde,\widetilde{v})=\widetilde{A}_{1,1}\left(|\etatilde|-\dfrac{A_2}{2\widetilde{A}_{1,1}}|\widetilde{v}|\right)^2 + \widetilde{A}_{1,2}\etatilde^2 + \left(A_3-\dfrac{A_2^2}{4\widetilde{A}_{1,1}}\right)\widetilde{v}^2,
\end{equation*}
where we have set $\widetilde{A}_{1,1}=\frac{\mu_*^2+c_s^2}{4}$, and $\widetilde{A}_{1,2}:=A_1-\widetilde{A}_{1,1}$. We recall that $\mu_*<c_s$, so that $\widetilde{A}_{1,1}<\frac{c_s^2}{2}$. In addition, we set $\gd:=\frac{c_s^2}{2(\mu_*^2+c_s^2)}$. As a consequence of the preceding choices, by the decomposition~\eqref{décomposition orthogonale de widetile Q}, then the exponential decay~\eqref{estimée décroissance exponentielle à tout ordre pour eta_c et v_c} and Proposition~\ref{proposition: borne uniforme sur les quantités nu_c etc}, we have for any $(t,x)\in\R\times \mathcal{I}^c$,
\begin{equation*}
    |\widetilde{\eta}\big(t,x+\widetilde{a}(t)\big)|\leq |\eta_{\ctilde(t)}(x)|+|\varepsilontilde_\eta(t,x)|\leq K_d e^{-a_d \nu_{\ctilde(t)} |x|}+\normHun{\varepsilontilde_\eta}\leq K_d e^{-\iota_* R_0}+A_*\alpha_* .
\end{equation*}

Now, since $-2f'(1-\xi)$ tends to $c_s^2$ as $\xi$ tends to $0$, there exists $R_0$ large enough and $\tau_*$ small enough such that, up to shrinking once again the value of $\alpha_*$,
\begin{equation*}
    |\etatilde| \leq \dfrac{\gd}{2}\quad\text{and}\quad\widetilde{A}_{1,1}< A_1 <\dfrac{c_s^2}{2}.
\end{equation*}

This implies that $A_3-\frac{A_2^2}{4\widetilde{A}_{1,1}}\geq \gd$ and $\widetilde{A}_{1,2} \geq \frac{c_s^2}{8}$ which provides the suitable constant $\kappa_*$.
\end{proof}

\begin{rem}
    The articulation of the previous proof does not depend on the special solution $\widetilde{Q}$, so that it still holds for a generic solution $Q$ of~\eqref{NLShydro} that is orbitally close to $Q_{c^*}$.
\end{rem}

Before passing to the proof of Corollary~\ref{corollaire: controle du moment localisé en |x|^delta}, we deduce from the monotonicity formula that
\begin{equation}\label{limite widetilde p_R vers 0}
    \sup_{t\in\R}\big|\widetilde{p}_R(t)\big|\underset{R\rightarrow +\ii}{\longrightarrow}0\quad\text{and}\quad \sup_{t\in\R}\big|\widetilde{p}_R(t)-p\big(\widetilde{Q}(t)\big)\big|\underset{R\rightarrow +\ii}{\longrightarrow}0.
\end{equation}

This can be proved exactly the same way than it is proved in~\cite{BetGrSm2} (see Proposition~$3$ in the latter article). The final contradiction in the proof of this proposition still holds for the general nonlinearity $f$ since we can control $\widetilde{p}_R$ in terms of the energy similarly as in~\eqref{controle des termes du moment par l'énergie} so that we get
\begin{equation*}
    \big|\widetilde{p}_R(t)\big|\leq\dfrac{1}{2}\int_\R|\etatilde \widetilde{v}|\leq \dfrac{4}{c_s^2}\int_\R F(1-\etatilde)+\dfrac{1}{\iota_*}\int_\R (1-\etatilde)\widetilde{v}^2 \lesssim E(\widetilde{Q}).
\end{equation*}

\begin{proof}[Proof of Corollary~\ref{corollaire: controle du moment localisé en |x|^delta}]
Arguing exactly the same way that it is done in the proof of Proposition~4 in~\cite{BetGrSm2}, we deduce the existence of a constant such that
\begin{equation*}
    \int_{t}^{t+1}\int_\R \Big((\partial_x\widetilde{\eta})^2 + \etatilde^2 + \widetilde{v}^2\Big)\big(s,x+\widetilde{a}(s)\big)\chi'(x-R) dx ds = \grandOde{e^{-\tau_* R}}.
\end{equation*}

Now by the special choice of $\chi$ and Proposition~\ref{proposition: borne uniforme sur les quantités nu_c etc}, we derive 
\begin{equation}\label{monotonie du moment en e^N_*|x| + grandO}
     \int_{t}^{t+1}\int_\R \Big((\partial_x\widetilde{\eta})^2 + \etatilde^2 + \widetilde{v}^2\Big)\big(s,x+\widetilde{a}(s)\big)e^{\tau_*|x|} dx ds =\grandOde{1}.
\end{equation}

Since, for any $\rho\geq 0$, there exists a constant $\kappa_{\rho}$ such that for any $x\in\R$,
\begin{equation}\label{|x|^rho leq k_rho e^tau x}
    1+|x|^\rho\leq \kappa_{\rho} e^{\tau_* |x|},
\end{equation}
this is enough to conclude the proof.
\end{proof}

\begin{proof}[Proof of Proposition~\ref{proposition: controle L^2_loc L^2_x partial_x^l psi}]
Recall that $\widetilde{\Psi}(t,x)=\widetilde{\psi}\big(t,x+\widetilde{a}(t)\big)$. Using the equation~\eqref{NLS}, we obtain at least formally
\begin{align*}
    \Big\Vert \partial_s\partial_x^l \widetilde{\Psi} &\Big\Vert_{L^2\big(I,L^2(|x|^rdx)\big)}^2 \leq 2\int_I\int_\R\bigg( \left|\partial_x^l\partial_t\widetilde{\psi}\big(s,x+\widetilde{a}(s)\big)\right|^2 + \big|\widetilde{a}'(s)\big|^2\left|\partial_x^{l+1}\widetilde{\psi}\big(s,x+\widetilde{a}(s)\big)\right|^2 \bigg)|x|^r dxds\\
    &\leq 4\Big(\left\Vert \partial_x^{l+2} \widetilde{\Psi} \right\Vert_{L^2\big(I,L^2(|x|^rdx)\big)}^2 +\left\Vert \partial_x^l \left(\Psi f(|\Psi|^2)\right)\right\Vert_{L^2\big(I,L^2(|x|^rdx)\big)}^2 +\left\Vert \widetilde{a}' \partial_x^l \widetilde{\Psi}\right\Vert_{L^2\big(I,L^2(|x|^rdx)\big)}^2\Big).
\end{align*}

We are going to prove by induction that there exist positive constants $C_{1,l}$ and $C_{2,l}$ such that for any positive integer $l\in\{1,...,l_*\}$ and any compact interval $I$ of length $1$,
\begin{equation}\label{controle par C_1,l}
    \left\Vert \partial_x^l \left(\Psi f(|\Psi|^2)\right)\right\Vert_{L^2\big(I,L^2(e^{\tau_*|x|}dx)\big)} \leq C_{1,l}
\end{equation}
and
\begin{equation}\label{controle par C_2,l}
    \left\Vert \partial_x^{l} \widetilde{\Psi} \right\Vert_{L^2\big(I,L^2(e^{\tau_*|x|}dx)\big)}\leq C_{2,l},
\end{equation}
which, by~\eqref{|x|^rho leq k_rho e^tau x}, will complete the proof of both cases $\Lambda\in\{0,1\}$ simultaneously. The induction principle is to be proved rather for weights of the form $e^{\tau_*|x|}$ since it relies drastically on a smoothing effect that is stated for exponential weights (see Lemma~\ref{lemme: smoothing effect de schrodinger} below) and that can be iterated until the maximum regularity of $f$ is reached.

For $l=1$, we can control the first derivative of $\psi$ by the suitable derivatives of the hydrodynamic variables. Indeed,
\begin{equation*}
    \big|\partial_x\widetilde{\psi}\big|^2=\dfrac{(\partial_x\etatilde)^2}{4(1-\etatilde)} + (1-\etatilde)\widetilde{v}^2.
\end{equation*}
Thus, by Proposition~\ref{proposition: borne uniforme sur les quantités nu_c etc} 
\begin{equation}\label{controle par C_2,1}
    \big|\partial_x\widetilde{\Psi}(t,x)\big|^2\leq \dfrac{\big(\partial_x\etatilde(t,x+\widetilde{a}(t))\big)^2}{4\iota_*} +(1+I_*)  \big(\partial_x\widetilde{v}(t,x+\widetilde{a}(t))\big)^2.
\end{equation}
Therefore,~\eqref{monotonie du moment en e^N_*|x| + grandO} provides the desired constant $C_{2,1}$. On the other hand, we compute 
\begin{equation*}
    \Big|\partial_x\left(\widetilde{\psi} f(|\widetilde{\psi}|^2)\right)\Big|\leq |\partial_x\widetilde{\psi}| \big|f(|\widetilde{\psi}|^2)\big| + 2|\widetilde{\psi}|^2|\partial_x\widetilde{\psi}| \big|f'(|\widetilde{\psi}|^2)\big|.
\end{equation*}
Since $f$ and $f'$ are continuous, and $\Vert\widetilde{\psi}\Vert_{L^\ii}\leq I_*$, we derive the existence of a constant $\widetilde{C}_{1,1}$ such that 
\begin{equation*}
    \Big|\partial_x\left(\widetilde{\Psi} f(|\widetilde{\Psi}|^2)\right)\Big|\leq \widetilde{C}_{1,1}|\partial_x\widetilde{\Psi}|.
\end{equation*}
Combining the previous estimate with~\eqref{controle par C_2,l}, we derive the desired constant $C_{1,1}$.

Now assume that there exist constants $C_{1,l}$ and $C_{2,l}$ such that estimates~\eqref{controle par C_1,l} and~\eqref{controle par C_2,l} hold true for $1\leq l\leq m$ with $m\leq l_*-1$. To prove the statement for the next integer, we shall rely on a special smoothing property of~\eqref{NLS} whose understanding dates back to pioneering articles such as~\cite{EsKePoV5}. We refer to Proposition~5 in~\cite{BetGrSm2} for its proof.
\begin{lem}\label{lemme: smoothing effect de schrodinger}
    Let $\lambda\in\R$. Consider a solution $w\in C^0\big(\R,L^2(\R)\big)$ to 
    \begin{equation*}
        i\partial_t w + \partial_x^2w =F,
    \end{equation*}
    with $F\in L^2\big(\R,L^2(\R)\big)$. Then, there exists a positive constant $K_\lambda$, depending only on $\lambda$, such that for any $T>0$,
    \begin{equation*}
        \lambda^2\int_{-T}^T\int_\R |\partial_x u(t,x)|^2e^{\lambda x}dxdt\leq K_\lambda \int_{-T-1}^{T+1}\int_\R\Big(|u(t,x)|^2+|F(t,x)|^2\Big)e^{\lambda x}dxdt.
    \end{equation*}
\end{lem}

We now apply Lemma~\ref{lemme: smoothing effect de schrodinger} with $u(s,x)=\partial_x^m \widetilde{\Psi}\big(s+t+\frac{1}{2},x\big),T=\frac{1}{2}$ and successively $\lambda=\pm \tau_*$. There exists $\widetilde{K}_*$ depending only on both $K_{\pm\tau_*}$ such that
\begin{equation*}
    \int_{t}^{t+1}\int_\R \big|\partial_x^{m+1}\widetilde{\Psi}(s,x)\big|^2 e^{\tau_*|x|}dxds\leq \widetilde{K}_*\int_{t-1}^{t+2}\int_\R \Big(\big|\partial_x^m\widetilde{\Psi} \big|^2+\Big|\partial_x^m\left(\widetilde{\Psi} f(|\widetilde{\Psi}|^2)\right)\Big|^2\Big)(s,x)e^{\tau_*|x|}dxds.
\end{equation*}

Then 
\begin{equation*}
    \int_{t}^{t+1}\int_\R \big|\partial_x^{m+1}\widetilde{\Psi}(s,x)\big|^2 e^{\tau_*|x|}dxds\leq C_{2,m+1},
\end{equation*}
with $C_{2,m+1}:=3\widetilde{K}_*(C_{1,m}+C_{2,m})$.

Regarding the other term, we deduce from~\eqref{NLS} and the fact that$C_{2,m+1}$ is independent of $t$, that for any $l\in\{0,...,m-1\}$, $\partial_x^l\widetilde{\psi}\in L^2\big(\R,H^2(\R)\big)\cap H^1\big(\R,L^2(\R)\big)$. Thus, up to taking a larger constant $C_{2,m+1}$, we have \begin{equation}\label{Vertpartial_x^lwidetildepsirightVert_L^ii(R^2)leq C_2,m+1}
    \left\Vert\partial_x^l\widetilde{\psi}\right\Vert_{L^{\ii}(\R^2)}\leq C_{2,m+1}.
\end{equation}

We shall use the following claim, the proof of which can be derived from the Faà di Bruno's formula.
\begin{claim}\label{claim: faa di bruno}
For $g\in C^{L}(\R_+)$ with $L\geq 0$ and any generic function $\Psi\in H^L(\R)$,
\begin{equation*}
    \Big|\partial_x^L\left(\Psi g(|\Psi|^2)\right)\Big|\leq K\sum_{n=0}^L \left|g^{(n)}(|\Psi|^2)\right|\sum_{\substack{i_0+i_1+...+i_L=2n+1\\ i_1+2i_2+...+Li_L=L}}|\Psi|^{i_0}|\partial_x\Psi|^{i_1}...|\partial_x^L\Psi|^{i_L}.
\end{equation*}
\end{claim}

We shall use the previous claim with the nonlinearity $f\in C^L(\R)$ with $L=m+1\leq l_*$. We separate the argument into several cases. First, the case $m=1$ being solved already, we can assume $m\geq 2$. Let $n\in\{0,...,m+1\}$, we have $i_{m+1}\in\{0,1\}$. Assume that $i_{m+1}=0$ then either $i_m=1$ or $i_m=0$. If $i_m=1$, we use~\eqref{Vertpartial_x^lwidetildepsirightVert_L^ii(R^2)leq C_2,m+1} so that
\begin{align*}
    \left|f^{(n)}(|\widetilde{\psi}|^2)\right||\widetilde{\psi}|^{i_0}|\partial_x\widetilde{\psi}|^{i_1}...|\partial_x^{m+1}\widetilde{\psi}|^{i_{m+1}} &\leq \left\Vert f^{(n)}(|\widetilde{\psi}|^2)\right\Vert_{L^{\ii}(\R^2)}\Vert{\widetilde{\psi}}\Vert_{L^{\ii}(\R^2)}^{i_0}...\Vert{\partial_x^{m-1}\widetilde{\psi}}\Vert_{L^{\ii}(\R^2)}^{i_{m-1}}\left|\partial_x^m\widetilde{\psi}\right|,
\end{align*}

hence, up to taking a larger $K$,
\begin{equation*}
    \left\Vert \partial_x^{m+1} \left(\widetilde{\psi} f(|\widetilde{\psi}|^2)\right)\right\Vert_{L^2\big(I,L^2(e^{\tau_*|x|}dx)\big)} \leq K C_{2,m+1}^{2m+2} \max_{0\leq n\leq m+1} \left\Vert f^{(n)}(|\widetilde{\psi}|^2)\right\Vert_{L^{\ii}(\R^2)} \left\Vert \partial_x^{m} \widetilde{\psi} \right\Vert_{L^2\big(I,L^2(e^{\tau_*|x|}dx)\big)},
\end{equation*}

which provides the desired constant $C_{1,m+1}$, depending on $C_{2,m+1},C_{2,m}$ and the derivatives of $f$. The case $i_m=0$ can be dealt with similarly to~\eqref{controle par C_2,1} by arguing on derivatives of lower order than $i_m$. Finally, if $i_{m+1}=1$, there is no other choice than $i_1=...=i_{m}=0$ and then $i_0=2n$. Therefore for any $n$, the suitable control in this case reads 
\begin{align*}
    \left|f^{(n)}(|\widetilde{\psi}|^2)\right||\widetilde{\psi}|^{i_0}|\partial_x\widetilde{\psi}|^{i_1}...|\partial_x^{m+1}\widetilde{\psi}|^{i_{m+1}} &\leq \left\Vert f^{(n)}(|\widetilde{\psi}|^2)\right\Vert_{L^{\ii}(\R^2)}\Vert{\widetilde{\psi}}\Vert_{L^{\ii}(\R^2)}^{2n}\left|\partial_x^{m+1}\widetilde{\psi}\right|,
\end{align*}

and then, up to taking a further $C_{2,m+1}$,
\begin{equation*}
    \left\Vert \partial_x^{m+1} \left(\widetilde{\psi} f(|\widetilde{\psi}|^2)\right)\right\Vert_{L^2\big(I,L^2(e^{\tau_*|x|}dx)\big)} \leq  C_{2,m+1} \max_{0\leq n\leq m+1} \left\Vert f^{(n)}(|\widetilde{\psi}|^2)\right\Vert_{L^{\ii}(\R^2)} \left\Vert \partial_x^{m+1} \widetilde{\psi} \right\Vert_{L^2\big(I,L^2(e^{\tau_*|x|}dx)\big)}.
\end{equation*}

This implies~\eqref{controle par C_1,l} for $l=m+1$. In conclusion, we have proved that~\eqref{controle par C_1,l} and~\eqref{controle par C_2,l} hold true for $1\leq l\leq l_*$, which ends the proof.
\end{proof}

\section{Proof of the Liouville property}

Before handling the proof of both Lemmas~\ref{lemme: equation sur la variable duale epsilon^*} and~\ref{lemme: monotonicity formula sur le viriel}, we assert a lemma taking care of potential $L^2$-norms of high derivatives of $\varepsilontilde$. Since both these lemma rely on the assumptions of Proposition~\ref{lemme: rigidity property proche d'un soliton, version algébrique} we assume that these latter conditions are fulfilled, namely we suppose that~\eqref{lemme: decaying property in rigidity} holds. According to what is done in the proof of Lemma~\ref{lemme: monotonicity formula sur le viriel}, see~\eqref{value of l_0}, we set the value $l_0:=12$. By Remark~\ref{remarque: décroissance sur widetilde Q implique celle de varepsilontilde}, we recall that there exists $r>5/2$ such that for any integer $l\in\{0,...,l_0\}$ and any $(t,x)\in\R^2$, we have
\begin{equation}\label{decay sur varepsilon dans section liouville}
    \left|\partial_x^{l+1}\widetilde{\varepsilon}_\eta\big(t,x+\widetilde{a}(t)\big)\right|+\left|\partial_x^{l}\widetilde{\varepsilon}_\eta\big(t,x+\widetilde{a}(t)\big)\right|+\left|\partial_x^{l}\widetilde{\varepsilon}_v\big(t,x+\widetilde{a}(t)\big)\right|\leq \dfrac{C_*}{1+|x|^r}.
\end{equation}

To simplify the notations, we get rid of the time-dependence and write for instance $\varepsilontilde$ instead of $\varepsilontilde(t)$. We also define some weighted norms that will be useful throughout the article. For $\rho\geq 0$ and any integer $l\geq -1$, we write 
\begin{equation*}
    \mathcal{X}_\rho^{l}:=\left\lbrace Q=(\eta,v)\in\mathcal{X}^l(\R) \Big|  \normXpoidsx{Q}{\rho}{l} < +\ii\right\rbrace ,
\end{equation*}

with
\begin{equation*}
    \normXpoidsx{Q}{\rho}{l}^2 = \sum_{m=0}^{l+1}\int_\R \big(\partial_x^m\eta(x)\big)^2|x|^\rho dx+\sum_{m=0}^{l}\int_\R \big(\partial_x^m v(x)\big)^2|x|^\rho dx,
\end{equation*}

with the convention $\mathcal{X}^{-1}(\R)=L^2(\R)\times L^2(\R)$ and $\normXpoidsx{Q}{\rho}{-1}$ the corresponding weighted $\big(L^2(|x|^\rho dx\big)^2$-norm.

\begin{lem}\label{lemme: controle des dérivées de varepsilontilde par normX varepsilon}
    There exists a positive constant $B_*$ such that for any given $(l,d)\in \{0,...,\frac{l_0}{2}\}\times\R_+$ such that either $d=0$, or $d\neq 0$ and $l\leq d< r-\frac{1}{2}$, then
    \begin{equation*}
        \normXpoidsx{\varepsilontilde}{^d}{l}\leq B_* \normXpoidsx{\varepsilontilde}{d}{\!}^{\frac{1}{2}} .
    \end{equation*}
\end{lem}

\begin{proof}
Consider first the case $d=0$. The case $l=0$ is straightforward so we consider a positive integer $l$. By doing enough integrations by parts, and using the Cauchy-Schwarz inequality, we get
\begin{equation*}
    \normLdeux{\partial_x^l\varepsilontilde_\eta}^2\leq\normLdeux{\partial_x\varepsilontilde_\eta}\normLdeux{\partial_x^{2l-1}\varepsilontilde_\eta}.
\end{equation*}

From~\eqref{decay sur varepsilon dans section liouville} and the fact that $2l-1\leq l_0-1\leq l_0+1$ by hypothesis, there exists a constant $B_*$ such that
\begin{equation*}
    \normxxk{\varepsilontilde}{l}\leq B_*\normX{\varepsilontilde}^\frac{1}{2}.
\end{equation*}

Now if $d>0$, assume moreover that $l\leq d<r-\frac{1}{2}$. From standard arguments, $x\mapsto|x|^d$ is $l$ times differentiable in a distributional sense and \begin{equation*}
    \Big|\partial_x^l\big(|x|^d\big) \Big|\leq  d(d-1)...(d-l+1)|x|^{d-l}.
\end{equation*}

We wish to evaluate the quantity $\normXpoidsx{\varepsilontilde}{d}{l}$. The highest order of differentiation in the $\mathcal{X}^l$-norm holds on the first component of $\varepsilontilde$ and the corresponding weighted $L^2$-norm reads 
\begin{align*}
    \int_\R (\partial^{l+1}_x\varepsilontilde_\eta)^2|x|^{d}dx \leq \normLdeux{\partial_x\varepsilontilde_\eta}\normLdeux{\partial_x^{l}\big(|x|^{d}\partial_x^{l+1}\varepsilontilde)}&\leq \normLdeux{\partial_x\varepsilontilde_\eta}\sum_{m=0}^{l}\begin{pmatrix}
        l\\
        m
    \end{pmatrix}\normLdeux{|x|^{d-m}\partial_x^{2l+1-m}\varepsilontilde_\eta}.
\end{align*}

Using~\eqref{decay sur varepsilon dans section liouville} and that $2l+1-m\leq l_0+1$ for any $m\in\{0,...,l\}$, and up to enlarging the value of $C_*$, we infer
\begin{equation*}
    \int_\R (\partial^{l+1}_x\varepsilontilde_\eta)^2|x|^{d}dx \leq C_*\normX{\varepsilontilde}\sum_{m=0}^{l} \normLdeux{\dfrac{|x|^{d-m}}{1+|x|^r}}.
\end{equation*}
In view of the special choice of $d$, each of the $L^2$-norms in the sum above is finite. Thus, up to a larger $B_*$, we derive
\begin{equation*}
        \normXpoidsx{\varepsilontilde}{^d}{l}^2\leq B_*^2\normXpoidsx{\varepsilontilde}{d}{\!} .
    \end{equation*}
\end{proof}

In the first place, let us prove Lemma~\ref{lemme: equation sur la variable duale epsilon^*}.
\begin{proof}[Proof of Lemma~\ref{lemme: equation sur la variable duale epsilon^*}]
In view of~\eqref{lemme: decaying property in rigidity}, $\widetilde{Q}_0\in\Nenergysethydrokk{4}$. By the local well-posedness result in Theorem~\ref{théorème: local well posedness in hydro} and by orbital stability implying the global well-posedness of $\widetilde{Q}$, we deduce $\widetilde{Q}\in C^0\big(\R,\Nenergysethydrokk{4}\big)$. By the equation~\eqref{NLShydro}, this implies that $\widetilde{Q}\in C^1\big(\R,\Nenergysethydrokk{2}\big)$ and by the decomposition~\eqref{décomposition orthogonale de widetile Q} and the smoothness of the modulation parameters, we obtain that $\varepsilontilde \in C^0\big(\R,\Nenergysethydrokk{4}\big)\cap C^1\big(\R,\Nenergysethydrokk{2}\big)$. Now, by definition of $\etilde$ and Proposition~\ref{proposition: propriété du linéarisé mathcal H_c}, we deduce that $\etilde\in C^0\big(\R,\Nenergysethydrokk{2}\big)\cap C^1\big(\R,\Nenergysethydro\big)$ and we can use the equation on $\varepsilontilde$ that is~\eqref{equation sur variable varepsilon} and the fact that $\partial_x Q_{\widetilde{c}(t)}\in\ker\big(\mathcal{H}_{\widetilde{c}(t)}\big)$ so that we get the desired equation on $\etilde$.
Finally, we deduce~\eqref{controle de varepsilon^* par widetildee} from using~\eqref{coercivité de la forme quadratique H_c sous condition d'orthogonalité} and Proposition~\ref{proposition: borne uniforme sur les quantités nu_c etc}, as well as the estimate
\begin{align*}
    \iota_*\normX{\varepsilontilde}^2\leq \psLdeuxLdeux{S\mathcal{H}_{\ctilde}(\varepsilontilde)}{S\varepsilontilde}=\psLdeuxLdeux{\etilde}{S\varepsilontilde}\leq \normX{\etilde}\normX{\varepsilontilde}.
\end{align*}
\end{proof}

Now that all the quantities are smooth and well-defined, we can prove Lemma~\ref{lemme: monotonicity formula sur le viriel}.

\begin{proof}[Proof of Lemma~\ref{lemme: monotonicity formula sur le viriel}]
Firstly, let us justify that $n$ is properly defined and differentiable. We define the scalar function $\mu(x)=x$ and 
we write $\widetilde{c}$ instead of $\widetilde{c}(t)$ throughout the section. Recall that
\begin{equation*}
    M_{\ctilde}(x)=\begin{pmatrix}
        m_{1,\ctilde} & m_{2,\ctilde}\\
        m_{2,\ctilde} & 0
    \end{pmatrix} \quad \text{ with }m_{1,\ctilde}:=-\dfrac{\eta_{\ctilde}'}{\eta_{\ctilde}}\text{ and }m_{2,\ctilde}:=-\dfrac{{\ctilde}\eta_{\ctilde}'}{2(1-\eta_{\ctilde})^2}.
\end{equation*}

Since $M_{\ctilde(t)}$ is bounded, $\partial_t\etilde\in C^0(\R,\Nenergysethydro)$, by Lemma~\ref{lemme: equation sur la variable duale epsilon^*} and the fact that $\widetilde{c}$ is differentiable, we deduce that $t\mapsto\psLdeuxLdeux{M_{\widetilde{c}(t)}\widetilde{e}(t)}{\widetilde{e}(t)}$ is well-defined and even differentiable.

Now, let us deal with the part involving the multiplication by $\mu S$. The fact that the scalar product $\psLdeuxLdeux{\mu S\etilde(t)}{\etilde(t)}$ is well-defined relies on previous estimates, independent of the articulation of the proof, that are given in~\eqref{leq widetildC_d^2left(2normLdeux}. Now to differentiate, we need to verify that every term in the right-hand side of 
\begin{align}\label{equation sur variable duale epsilon preuve monotonie de n}
    \partial_t \widetilde{e} = &-2S\mathcal{H}_{\ctilde}(\partial_x \widetilde{e})+S\mathcal{H}_{\ctilde}(J\mathcal{R}_{\ctilde}\widetilde{\varepsilon}) - \widetilde{c}'(t)S\mathcal{H}_{\ctilde}(\partial_c Q_{\ctilde})\\ \notag
    &+\widetilde{c}'(t)S\partial_c\mathcal{H}_{\ctilde}(\widetilde{\varepsilon})+\big((\widetilde{a}'(t)-\ctilde\big)S\mathcal{H}_{\ctilde}(\partial_x\widetilde{\varepsilon})
\end{align}

in Lemma~\ref{lemme: equation sur la variable duale epsilon^*} can be integrated against the function $\mu S \etilde(t)$. By Corollary~\ref{corollaire: controle de la norme de mathcal H_c J mathcal R_c}, we obtain 
\begin{align}
    \left|\psLdeuxLdeux{-2S\mathcal{H}_{\ctilde}(\partial_x \widetilde{e})}{\mu S\etilde}\right|&\leq \normLdeuxLdeux{\sqrt{|x|}\mathcal{H}_c(\partial_x\etilde)}\normLdeuxLdeux{\sqrt{|x|}\etilde}\\
    &\leq K_*^2 \normXpoidsx{\varepsilontilde}{\!}{3}\normXpoidsx{\varepsilontilde}{\!}{1}.
\end{align}

Using now the decay~\eqref{decay sur varepsilon dans section liouville}, we obtain that the right-hand term in the latter estimate is bounded uniformly in time. The other terms can be dealt similarly by using Proposition~\ref{proposition: propriété du linéarisé mathcal H_c} with different $m$ and $l$. Therefore one can differentiate the function $n$ and we get
\begin{align*}
    n'(t)=\dfrac{d}{dt}\psLdeuxLdeux{N(t)\widetilde{e}(t)}{\widetilde{e}(t)}
    &=n_1(t)+\gamma_* \big(n_2(t)+n_3(t)\big),
\end{align*}

with 
\begin{align}
    n_1(t):= & 2\psLdeuxLdeux{\mu S\partial_t\etilde}{\etilde}\notag,\\
    n_2(t):=&   2\psLdeuxLdeux{ M_{\ctilde}\partial_t\etilde}{\etilde}\notag,\\
    n_3(t):=&\widetilde{c}'(t)\psLdeuxLdeux{\partial_c M_{\ctilde}\etilde}{\etilde}\label{terme n_3(t)}.
\end{align}

We compute
\begin{align}
    n_1(t)=& -4\psLdeuxLdeux{\mu\mathcal{H}_{\ctilde}\big(\partial_x\etilde\big)}{\etilde}+2\psLdeuxLdeux{\mu\mathcal{H}_{\ctilde}J\mathcal{R}_{\ctilde}\big(\varepsilontilde\big)}{\etilde}\label{deux premiers termes dans n_1(t)}\\
    &-2\widetilde{c}'(t)\psLdeuxLdeux{\mu\mathcal{H}_{\ctilde}(\partial_c Q_{\ctilde})}{\etilde}+2\widetilde{c}'(t)\psLdeuxLdeux{\mu\partial_c \mathcal{H}_{\ctilde}\big(\varepsilontilde\big)}{\etilde}\label{troisième et quatrième terme dans n_1(t)}\\
    &+2\big(\widetilde{a}'(t)-\ctilde\big)\psLdeuxLdeux{\mu \mathcal{H}_{\ctilde}\big(\partial_x\varepsilontilde\big)}{\etilde}\label{dernier terme dans n_1(t)}, 
\end{align}

and
\begin{align}
    n_2(t)=& \quad 2\big(\widetilde{a}'(t)-\ctilde\big)\psLdeuxLdeux{M_{\ctilde}S\mathcal{H}_{\ctilde}\partial_x\varepsilontilde}{\etilde}+2\psLdeuxLdeux{M_{\ctilde}S\mathcal{H}_{\ctilde}J\mathcal{R}_{\ctilde}\big(\varepsilontilde\big)}{\etilde}\label{deux premiers termes dans n_2(t)}\\
    &-2\widetilde{c}'(t)\psLdeuxLdeux{M_{\ctilde}S\mathcal{H}_{\ctilde}\partial_c Q_{\ctilde}}{\etilde}+2\widetilde{c}'(t)\psLdeuxLdeux{\partial_c M_{\ctilde}S\mathcal{H}_{\ctilde}\varepsilontilde}{\etilde}\label{troisième et quatrième terme dans n_2(t)}\\
    &-4\psLdeuxLdeux{M_{\ctilde}S\mathcal{H}_{\ctilde}\partial_x\etilde}{\etilde}\label{dernier terme dans n_2(t)}.
\end{align}

We are going to deal with each part separately.\\
\underline{Step 1.} Regarding the term~\eqref{dernier terme dans n_2(t)}. We prove that there exists $\Xi_1 >0$ depending continuously on $c$ such that \begin{equation*}
    -4\psLdeuxLdeux{M_{\ctilde}S\mathcal{H}_{\ctilde}\partial_x\etilde}{\etilde}\geq \Xi_1\int_\R \eta_{\ctilde}\left((\partial_x\widetilde{e}_\eta)^2 + \etildeeta^2 + \etildev^2\right).
\end{equation*}

To prove Step 1, we rely on the following claim. This claim relies on tedious computations. We refer to the end of this section for its proof. From now on and until the end of Step 1, we write $c$ instead of $\ctilde$.

\begin{claim}\label{claim: forme quadratique du terme fastidieux}
We have
\begin{align*}
    -4\psLdeuxLdeux{M_{c}S\mathcal{H}_{c}\partial_x\etilde}{\etilde}=& \int_\R\bigg(q_{1,c}\Big(\etildev+\dfrac{q_{2,c}}{2q_{1,c}} \etildeeta+\dfrac{q_{3,c}}{2q_{1,c}}\partial_x\etildeeta\Big)^2+\widetilde{q_{1,c}}\Big(\partial_x\etildeeta +m_{1,c}\etildeeta\Big)^2\bigg)
\end{align*}

where \begin{align*}
    B_c&=1-\eta_c,\\
    q_{1,c}&=2(m_{1,c} B_c)',\\
    q_{2,c}&=-2cm_{1,c}',\\
    q_{3,c}&=4 m_{2,c} B_c,\\
    \widetilde{q_{1,c}}&=-\dfrac{q_{3,c}^2}{4q_{1,c}}+\dfrac{m_{1,c}'}{B_c}+\Big(\dfrac{m_{1,c}}{2B_c}\Big)'.\\
\end{align*}

\end{claim}

Once we have shown Claim~\ref{claim: forme quadratique du terme fastidieux}, we use Proposition~\ref{proposition: signe de q_1,c et q_2,c} so that all the quantities in the following claim are well-defined and moreover $q_{1,c}>0$ and $\widetilde{q_{1,c}}>0$. In particular, the quadratic form in the integral is non-negative and we check that its null-space is given by $\spanned{Q_c}$. By~\eqref{développement asymptotique de q_1,c en grand o de eta_c^2}, we notice that 
\begin{equation}\label{limite à l'infini de q_1,c/eta_c}
    \dfrac{q_{1,c}(x)}{\eta_c(x)}\underset{|x|\rightarrow +\ii}{\longrightarrow}k_0:=2\nu_c^2-\dfrac{k}{3}>0, 
\end{equation}
so that it is natural to substitute the variable $\etilde$ by $\ftilde:=\sqrt{\eta_c}\ \etilde$ in order to provide
\begin{align*}
    -4\psLdeuxLdeux{M_{c}S\mathcal{H}_{c}\partial_x\etilde}{\etilde}= \int_\R\bigg(\dfrac{q_{1,c}}{\eta_c}\Big(\ftildev+\Big(\dfrac{q_{2,c}}{2q_{1,c}} +\dfrac{q_{3,c}m_{1,c}}{4q_{1,c}}\Big)&\ftildeeta+\dfrac{q_{3,c}}{2q_{1,c}}\partial_x\ftildeeta\Big)^2\\
    &+\dfrac{\widetilde{q_{1,c}}}{\eta_c}\Big(\partial_x\ftildeeta +\dfrac{3m_{1,c}}{2}\ftildeeta\Big)^2\bigg).
\end{align*}

Substituting the old variables once again by the new pair $\htilde=(\htildeeta,\htildev):=(\ftildeeta,\ftildev+\frac{q_{3,c}}{2q_{1,c}}\partial_x\ftildeeta)$, so that 
\begin{equation*}
    -4\psLdeuxLdeux{M_{c}S\mathcal{H}_{c}\partial_x\etilde}{\etilde}=\psLdeuxLdeux{T_c(\htilde)}{\htilde},
\end{equation*}

with \begin{equation*}
    T_c(\htilde):=\begin{pmatrix}
        -\partial_x\left(\dfrac{\widetilde{q_{1,c}}\partial_x \htildeeta}{\eta_c}\right) + \left(\dfrac{\eta_c\widetilde{q_{3,c}}^2}{ q_{1,c}}-\dfrac{3}{2}\partial_x \left(\dfrac{\widetilde{q_{1,c}}m_{1,c}}{\eta_c}\right)+\dfrac{9\widetilde{q_{1,c}}m_{1,c}^2}{4\eta_c}\right)\htildeeta + \widetilde{q_{3,c}}\htildev\\
        \widetilde{q_{3,c}}\htildeeta + \dfrac{q_{1,c}}{\eta_c}\htildev
        
    \end{pmatrix},
\end{equation*}

where $\widetilde{q_{3,c}}:=\frac{2q_{2,c}+q_{3,c}m_{1,c}}{4\eta_c}$, $T_c$ is non-negative and we can introduce 
\begin{equation*}
    \gh_c=\left(\eta_c^\frac{3}{2},-\dfrac{\widetilde{q_{3,c}}\eta_c^\frac{5}{2}}{q_{1,c}}\right)\in\energysethydro,
\end{equation*}
such that $\ker(T_c)=\spanned{\gh_c}$. One of the crucial use of the transonic limit is made here. It is the fact that the essential spectrum of $T_c$ in contained in $\R_+^*$, then assuring the coercivity of $T_c$, up to taking elements orthogonal to the kernel. It is summarized in the following claim, the proof of which is at the end of this section.
\begin{claim}\label{claim: spectre essentiel de T_c}
Up to taking $c_1$ closer to $c_s$, we have $\spess(T_c)=[\tau_c,+\ii)$ where $\tau_c>0$ depends continuously on $c\in (c_1,c_s)$, and there exists a positive constant $\Lambda_0$, that does not depend on $c$, such that $\tau_c\underset{c\rightarrow c_s}{\sim}\Lambda_0\nu_c^2$.
\end{claim}

Thus by the previous claim, there exists a constant $\Lambda_1>0$ depending continuously on $c$ such that for any $h$ orthogonal to $\gh_c$, we have 
\begin{equation*}
    \psLdeuxLdeux{T_c(h)}{h}\geq \Lambda_1 \normLdeuxLdeux{h}^2.
\end{equation*}

Moreover we have,
\begin{equation*}
    \left|\psLdeuxLdeux{T_c(h)}{h}-\int_\R\dfrac{\widetilde{q_{1,c}}}{\eta_c}\big|\partial_xh_\eta\big|^2\right|\leq A_1\normLdeuxLdeux{h}^2,
\end{equation*}

so that taking $\tau\in (0,1)$ and using Proposition~\ref{proposition: signe de q_1,c et q_2,c}, we have 
\begin{align}
    \psLdeuxLdeux{T_c(h)}{h}&\geq (1-\tau)\psLdeuxLdeux{T_c(h)}{h}+\tau\normLdeux{\left(\dfrac{\widetilde{q_{1,c}}}{\eta_c}\right)^\frac{1}{2}\partial_xh_\eta}^2 -A_1\tau \normLdeuxLdeux{h}^2\label{apparition de tau dans estimation de T_c h , h}\\
    &\geq \Lambda_1(1-\tau-A_1\tau)\normLdeuxLdeux{h}^2 + \inf_{(c,x)\in (c_1,c_s)\times\R}\left|\dfrac{\widetilde{q_{1,c}}(x)}{\eta_c(x)}\right| \normLdeux{\partial_x h_\eta}^2.
\end{align}

Reducing the value of $\tau$ enough, we obtain a constant $\Lambda_2>0$ depending continuously on $c$ such that for any $h$ orthogonal to $\gh_c$, we have
\begin{equation*}
    \psLdeuxLdeux{T_c(h)}{h}\geq \Lambda_2\normX{h}^2.
\end{equation*}

We express the previous estimate in terms of the original variable. First notice that $\htilde$ is orthogonal to $\gh_c$ if and only if $\etilde$ is orthogonal to $\ge_c$, where
\begin{equation*}
    \ge_c=\left(\eta_c-\dfrac{q_{3,c}\widetilde{q_{3,c}}\eta_c\eta_c'}{4q_{1,c}^2}+\partial_x\Big(\dfrac{q_{3,c}\widetilde{q_{3,c}}\eta_c^2}{2q_{1,c}^2}\Big),-\dfrac{\widetilde{q_{3,c}}\eta_c^2}{q_{1,c}}\right).
\end{equation*}

Furthermore, we have
\begin{equation*}
    -4\psLdeuxLdeux{M_{c}S\mathcal{H}_{c}\partial_x\etilde}{\etilde} \geq \Lambda_2 \left(\normHun{\htildeeta}^2 +\normLdeux{\htildev}^2\right).
\end{equation*}

First observe that
\begin{align*}
    \normHun{\htildeeta}^2 =\int_\R \eta_c\left((\partial_x\etildeeta)^2+\Big(1+\dfrac{(\eta_c')^2}{4\eta_c^2}-\dfrac{\eta_c''}{2\eta_c}\Big)\etildeeta^2 \right).
\end{align*}
Now by Proposition~\ref{proposition: borne inférieure eta_c(x)} and restricting to speeds $c\in (c_1,c_s)$, with $c_1$ even closer to $c_s$, which will be satisfied by $\ctilde$ up to shrinking the value of $\alpha_*$, we can assume that, uniformly in $x\in\R$, we have
\begin{equation*}
    1+\dfrac{(\eta_c')^2}{4\eta_c^2}-\dfrac{\eta_c''}{2\eta_c}\geq \dfrac{1}{2}.
\end{equation*}

Proceeding similarly than for~\eqref{apparition de tau dans estimation de T_c h , h} but with the standard inequality for $a,b\in L^2(\R)$ and $\tau\in (0,1)$,
\begin{equation*}
    \normLdeux{a-b}\geq \tau\normLdeux{a}-\dfrac{\tau}{1-\tau}\normLdeux{b},
\end{equation*}
we verify that, up to taking a smaller $\Lambda_2$ (still depending continuously on $c$), we have
\begin{equation*}
    -4\psLdeuxLdeux{M_{c}S\mathcal{H}_{c}\partial_x\etilde}{\etilde} \geq \Lambda_2 \int_\R \eta_c\left((\partial_x\etildeeta)^2 + \etildeeta^2 + \widetilde{e}_v^2\right),
\end{equation*}
whenever $\etilde$ is orthogonal to $\ge_c$. The rest of the argument can be deduced like in the proof of Proposition~8 in~\cite{BetGrSm2}, except for the part where it makes use of the explicit shape of the travelling wave, namely with the part involving the quantities $\partial_c\eta_c$ and $\partial_c v_c$. In this sense, we need to bound the integral
\begin{equation}\label{int_R 1/eta_c (partial_c eta_c^2}
    \int_\R \dfrac{1}{\eta_c}\left((\partial_c\eta_c)^2 + (\partial_c v_c)^2\right)
\end{equation}
by a constant that depending continuously on $c$. Recall that by~\eqref{formule de gu_c classique en fonction eta_c et v_c hydrodynamique}, we get $v_c=\dfrac{c\eta_c}{2(1-\eta_c)}$ so that 
\begin{equation}\label{formule de partial_c v_c}
    \partial_c v_c = \dfrac{\eta_c}{2(1-\eta_c)}+\dfrac{c\partial_c\eta_c}{2(1-\eta_c)^2}.
\end{equation}

Moreover, we invoke Lemma~$2.9$ in~\cite{Bert2} with the same setting apart from the fact that we impose $\sigma=\nu_c(\frac{1}{2}+\delta)$. Taking $\delta$ sufficiently small, we then obtain the decay $\big(\partial_c\eta_c(x)\big)^2\lesssim e^{-\nu_c(1+2\delta)|x|}$. Then by Proposition~\ref{proposition: borne inférieure eta_c(x)}, we obtain 
\begin{equation*}
    \big(\partial_c\eta_c(x)\big)^2\lesssim \dfrac{\eta_c(x)}{m_c} e^{-2\delta\nu_c|x|},
\end{equation*} 
and the same holds for $\partial_c v_c$ by~\eqref{formule de partial_c v_c}. Combining these considerations with Proposition~\ref{proposition: borne uniforme sur les quantités nu_c etc} implies that the integral in~\eqref{int_R 1/eta_c (partial_c eta_c^2} is bounded by a constant depending continuously on $c$.

\underline{Step 2.} We deal with the terms~\eqref{terme n_3(t)},~\eqref{troisième et quatrième terme dans n_1(t)},~\eqref{deux premiers termes dans n_2(t)} and \eqref{troisième et quatrième terme dans n_2(t)}. We show that they can be expressed as $\mathcal{O}\big(\alpha_*^{q}\normX{\etilde(t)}^2\big)$ for some real number $q>0$. Regarding~\eqref{terme n_3(t)},~\eqref{troisième et quatrième terme dans n_1(t)} and~\eqref{troisième et quatrième terme dans n_2(t)}, we use~\eqref{controle des parametres de modulation}, and the estimate~\eqref{controle de varepsilon^* par widetildee} so that we obtain $|\ctilde'(t)|\leq A_*\normX{\varepsilontilde(t)}^2=\grandOde{\alpha_* \normX{\etilde(t)}^2}$. Thus
\begin{align}
    |n_3(t)|&=\grandOde{ \alpha_*\normLii{\partial_c M_{\ctilde}}\normX{\etilde}^2}\label{premier terme des cinq en grandOde alpha_* normx etilde^2}\\
    \big|-2\widetilde{c}'(t)\psLdeuxLdeux{\mu\mathcal{H}_{\ctilde}(\partial_c Q_{\ctilde})}{\etilde}\big|&=\grandOde{\alpha_*\normLdeuxLdeux{\mu\mathcal{H}_{\ctilde}(\partial_c Q_{\ctilde})}\normX{\etilde}^2},\label{deuxième terme des cinq en grandOde alpha_* normx etilde^2}\\
    \big|2\widetilde{c}'(t)\psLdeuxLdeux{\mu\partial_c \mathcal{H}_{\ctilde}(\varepsilontilde)}{\etilde}\big|&= \grandOde{ \alpha_*\normLdeuxLdeux{\mu\partial_c\mathcal{H}_{\ctilde}(\varepsilontilde)}\normX{\etilde}^2},\label{troisième terme des cinq en grandOde alpha_* normx etilde^2}\\
    \big|-2\widetilde{c}'(t)\psLdeuxLdeux{M_{\ctilde}S\mathcal{H}_{\ctilde}(\partial_c Q_{\ctilde})}{\etilde}\big|&= \grandOde{ \alpha_* \normLii{M_{\ctilde}}\normLdeuxLdeux{\mathcal{H}_{\ctilde}(\partial_c Q_{\ctilde})}\normX{\etilde}^2},\label{quatrième terme des cinq en grandOde alpha_* normx etilde^2}\\
    \big|2\widetilde{c}'(t)\psLdeuxLdeux{\partial_c M_{\ctilde}S\mathcal{H}_{\ctilde}(\varepsilontilde)}{\etilde}\big|&= \grandOde{\alpha_*  \normLii{\partial_c M_{\ctilde}}\normLdeuxLdeux{\mathcal{H}_{\ctilde}(\varepsilontilde)}\normX{\etilde}^2}.\label{cinquième terme des cinq en grandOde alpha_* normx etilde^2}
\end{align}

Now, combining Proposition~\ref{proposition: propriété du linéarisé mathcal H_c}, the exponential decay~\eqref{estimée décroissance exponentielle à tout ordre pour eta_c et v_c} and Proposition~\ref{proposition: borne uniforme sur les quantités nu_c etc}, we control~\eqref{premier terme des cinq en grandOde alpha_* normx etilde^2},~\eqref{deuxième terme des cinq en grandOde alpha_* normx etilde^2} and~\eqref{quatrième terme des cinq en grandOde alpha_* normx etilde^2} as $\grandOde{\alpha_* \normX{\etilde}^2}$. We control~\eqref{troisième terme des cinq en grandOde alpha_* normx etilde^2} and~\eqref{cinquième terme des cinq en grandOde alpha_* normx etilde^2} similarly by invoking in addition the estimate~\eqref{decay sur varepsilon dans section liouville} and also the exponential decay~\eqref{estimée décroissance exponentielle à tout ordre pour eta_c et v_c} on the derivative of $\eta_c$ and $v_c$ with respect to $c$. 

The first term of~\eqref{deux premiers termes dans n_2(t)} can be dealt with similarly, using eventually Lemma~\ref{lemme: controle des dérivées de varepsilontilde par normX varepsilon} with $d=0$, which leads to
\begin{align*}
    \Big|2\big(\widetilde{a}'(t)-\ctilde\big)\psLdeuxLdeux{M_{\ctilde}S\mathcal{H}_{\ctilde}\partial_x\varepsilontilde}{\etilde}\big|&\leq 2\sqrt{A_*}M_*\normX{\varepsilontilde}\normLdeuxLdeux{\mathcal{H}_{\ctilde}(\partial_x\varepsilontilde)}\normX{\etilde}\\
    &\leq 2\sqrt{A_*}D_*M_*K_*\normxxk{\varepsilontilde}{2}\normX{\etilde}^2=\grandOde{\alpha_*^{\frac{1}{2}}\normX{\etilde}^2}.
\end{align*}

It remains to deal with the right-hand term of~\eqref{deux premiers termes dans n_2(t)}. We have by the Cauchy-Schwarz inequality
\begin{equation*}
    \Big|2\psLdeuxLdeux{M_{\ctilde}S\mathcal{H}_{\ctilde}J\mathcal{R}_{\ctilde}(\varepsilontilde)}{\etilde}\Big|\leq 2M_*\normLdeuxLdeux{\mathcal{H}_{\ctilde}J\mathcal{R}_{\ctilde}(\varepsilontilde)}\normX{\etilde}.
\end{equation*}

At this moment of the argument, we fix the value of $l_0$. Indeed, we now impose
\begin{equation}\label{value of l_0}
    l_0=12,
\end{equation}

and as a consequence of Corollary~\ref{corollaire: controle de H_c J R_c varepsilon par norme W24,3}, we obtain
\begin{equation*}
    \normxxk{\mathcal{H}_{\ctilde}J\mathcal{R}_{\ctilde}(\varepsilontilde)}{-1}\leq K_*\normX{\varepsilontilde}^{\frac{5}{4}}\max_{\substack{p\in\{1,...,3\}\\ l\in\{0,...,l_0\}}}\normWkp{\varepsilontilde}{l}{p}^{p+2}.
\end{equation*}

By~\eqref{decay sur varepsilon dans section liouville}, each of the previous $W^{l,p}$-norms are uniformly bounded and this implies that
\begin{align*}
    \Big|2\psLdeuxLdeux{M_{\ctilde}S\mathcal{H}_{\ctilde}J\mathcal{R}_{\ctilde}(\varepsilontilde)}{\etilde}\Big|&\leq  \ 2M_*K_*C_*\normX{\varepsilontilde}^\frac{5}{4}\normX{\etilde}=\grandOde{\alpha_*^{\frac{1}{4}}\normX{\etilde}^2}.
\end{align*}
Again, using~\eqref{controle de widetilde varepsilon et ctilde -c_0^* par alpha_*} and Lemma~\ref{lemme: equation sur la variable duale epsilon^*} provides the suitable control and ends the proof of Step 2 with $q=\frac{1}{4}$.\\

\underline{Step 3.} We prove that the terms~\eqref{deux premiers termes dans n_1(t)} and~\eqref{dernier terme dans n_1(t)} can be bounded from below by $\Xi_3\normX{\etilde}^2+\grandOde{\normXboule{\etilde}{R}^2}$ with $\Xi_3$ and $R$ positive constants depending only on $c^*$. We start by developing the first term in~\eqref{deux premiers termes dans n_1(t)}. Recall that $B_c=1-\eta_c$. Then in view of the expression of $\mathcal{H}_{\ctilde}$ in~\eqref{calcul de mathcal(H_c)}, we compute
\begin{align*}
    -4\psLdeuxLdeux{\mu\mathcal{H}_{\ctilde}\big(\partial_x\etilde\big)}{\etilde} =&\int_\R i(x)dx,\notag
\end{align*}
where
\begin{align}
     i&=\Big(\dfrac{3}{2B_{\ctilde}}-\dfrac{\mu\eta_{\ctilde}'}{2B_{\ctilde}^2}\Big)(\partial_x\etildeeta)^2 - \dfrac{2\ctilde(B_{\ctilde}-\mu\eta_{\ctilde}')}{B_{\ctilde}^2}\etildeeta\etildev + 2(B_{\ctilde}-\mu\eta_{\ctilde}')\etildev^2 \label{expression de -4psLdeuxLdeuxmu mathcalH_ctildebig(partial_xetildebig)etilde}\\
    - \bigg(\Big(\dfrac{\eta_{\ctilde}''}{B_{\ctilde}^2}&+\dfrac{3(\eta_{\ctilde}')^2}{2B_{\ctilde}^3}+f'(B_{\ctilde})\Big)+\mu\Big(\dfrac{\eta_{\ctilde}'''}{2B_{\ctilde}^2}+\dfrac{2\eta_{\ctilde}'\eta_{\ctilde}''}{B_{\ctilde}^3}+\dfrac{3(\eta_{\ctilde}')^2}{2B_{\ctilde}^4}-\eta_{\ctilde}'f'(B_{\ctilde})\Big)\bigg)\etildeeta^2\notag .
\end{align}

We compute the following limits
\begin{equation*}
    \dfrac{3}{2B_{\ctilde}}-\dfrac{\mu\eta_{\ctilde}'}{2B_{\ctilde}^2}\underset{|x|\rightarrow +\ii}{\longrightarrow}\dfrac{3}{2},\quad\quad - \dfrac{2\ctilde(B_{\ctilde}-\mu\eta_{\ctilde}')}{B_{\ctilde}^2}\underset{|x|\rightarrow +\ii}{\longrightarrow}-2\ctilde,\quad\quad 2(B_{\ctilde}-\mu\eta_{\ctilde}')\underset{|x|\rightarrow +\ii}{\longrightarrow} 2,
\end{equation*}
and
\begin{equation*}
    - \Big(\dfrac{\eta_{\ctilde}''}{B_{\ctilde}^2}+\dfrac{3(\eta_{\ctilde}')^2}{2B_{\ctilde}^3}+f'(B_{\ctilde})\Big)+\mu\Big(\dfrac{\eta_{\ctilde}'''}{2B_{\ctilde}^2}+\dfrac{2\eta_{\ctilde}'\eta_{\ctilde}''}{B_{\ctilde}^3}+\dfrac{3(\eta_{\ctilde}')^2}{2B_{\ctilde}^4}-\eta_{\ctilde}'f'(B_{\ctilde})\Big)\underset{|x|\rightarrow +\ii}{\longrightarrow} \dfrac{c_s^2}{2}.
\end{equation*}

Set $\tau >0$, we use the exponential decay~\eqref{estimée décroissance exponentielle à tout ordre pour eta_c et v_c} to exhibit $R_1>0$ large enough such that for $|x|\geq R_1$, we have by the Cauchy-Schwarz inequality,
\begin{align*}
    i(x)&\geq \dfrac{3}{2}\big(\partial_x\etildeeta(x)\big)^2 -2\ctilde\etildeeta(x)\etildev(x)+2\etildev(x)^2 +\dfrac{c_s^2}{2}\etildeeta(x)^2-\tau\left(\big(\partial_x\etildeeta(x)\big)^2+\etildeeta(x)^2 + \etildev(x)^2\right)\\
    &\geq \dfrac{3}{2}\big(\partial_x\etildeeta(x)\big)^2 +\dfrac{1}{2}\Big(c_s^2 - \dfrac{\ctilde^2}{1-\delta}\Big)\etildeeta(x)^2 +2\delta\etildev(x)^2-\tau\left(\big(\partial_x\etildeeta(x)\big)^2+\etildeeta(x)^2 + \etildev(x)^2\right).
\end{align*}

We now fix the value of $\delta >0$ such that $c^*+A_*\alpha_* < c_s\sqrt{1-\delta}$. Therefore, taking $\tau$ small enough, we obtain according to~\eqref{controle des parametres de modulation}, a constant $\Xi_3$ depending only on $c^*$ such that for $|x|\geq R_1$,
\begin{equation*}
    i(x)\geq \Xi_3\left(\big(\partial_x\etildeeta(x)\big)^2+\etildeeta(x)^2 + \etildev(x)^2\right),
\end{equation*}

hence
\begin{equation*}
    -4\psLdeuxLdeux{\mu\mathcal{H}_{\ctilde}\big(\partial_x\etilde\big)}{\etilde}=\int_{(\mathcal{B}_{R_1})^c} i(x)dx+\int_{\mathcal{B}_{R_1}} i(x)dx\geq \Xi_3 \normX{\etilde}^2 +\grandOde{\normXboule{\etilde}{R_1}^2}.
\end{equation*}

Next, we deal with the second term in~\eqref{deux premiers termes dans n_1(t)} and~\eqref{dernier terme dans n_1(t)}. On the one hand, take $\delta >0$, and since $r>\frac{5}{2}$, take 
\begin{equation}\label{choix de s dans la preuve de lemme de monotonie viriel}
    s\in \left[\dfrac{3}{2},\frac{r}{2}-\dfrac{1}{4}\right).
\end{equation}
Take also $R_2>0$ large enough such that $|x|\leq \delta (1+|x|^s)$ for any $|x|\geq R_2$. Then, decomposing the line into the part where $|x|\leq R_2$ and its complement, we obtain the estimates
\begin{align*}
    \Big|2\psLdeuxLdeux{\mu\mathcal{H}_{\ctilde}J\mathcal{R}_{\ctilde}\big(\varepsilontilde\big)}{\etilde}\Big| &\leq R_2\normxxk{\mathcal{H}_{\ctilde}J\mathcal{R}_{\ctilde}(\varepsilontilde)}{-1}\normXboule{\etilde}{R_2} +  \delta \normXpoidsunplusx{\mathcal{H}_{\ctilde}J\mathcal{R}_{\ctilde}(\varepsilontilde)}{2s}{-1} \normX{\etilde}.
\end{align*}

Using Young's inequality and then Corollary~\ref{corollaire: controle de la norme de mathcal H_c J mathcal R_c}, we deduce 
\begin{equation*}
    \normXpoidsunplusx{\mathcal{H}_{\ctilde}J\mathcal{R}_{\ctilde}(\varepsilontilde)}{2s}{-1}\leq 2K_*\big(\normxxk{\varepsilontilde}{3}^2+ \normXpoidsx{\varepsilontilde}{2s}{3}^2\big),
\end{equation*}
hence
\begin{align*}
    \Big|2\psLdeuxLdeux{\mu\mathcal{H}_{\ctilde}J\mathcal{R}_{\ctilde}\big(\varepsilontilde\big)}{\etilde}\Big|\leq K_*\Big( R_2\normxxk{\varepsilontilde}{3}^2\normXboule{\etilde}{R_2} +  2\delta \big(\normxxk{\varepsilontilde}{3}^2+\normXpoidsx{\varepsilontilde}{2s}{3}^2\big)\normX{\etilde}\Big) .
\end{align*}

Therefore, by the choice of $s$ given by~\eqref{choix de s dans la preuve de lemme de monotonie viriel}, we can apply Lemma~\ref{lemme: controle des dérivées de varepsilontilde par normX varepsilon} with $d\in\{0,2s\}$ and $l=3$ and there comes finally \begin{align}\label{controle de psLdeuxLdeuxmumathcaH_ctildeJmathcalR_ctildebig(varepsilontildebig)}
    \Big|2\psLdeuxLdeux{\mu\mathcal{H}_{\ctilde}J\mathcal{R}_{\ctilde}\big(\varepsilontilde\big)}{\etilde}\Big|\leq K_*B_*^2 \Big(R_2\normXboule{\etilde}{R_2} +4  \delta  \normX{\etilde}\Big)\normX{\varepsilontilde}.
\end{align}

Regarding~\eqref{dernier terme dans n_1(t)}, we can decompose the integral the same way than for the previous term and deduce that the control in~\eqref{controle de psLdeuxLdeuxmumathcaH_ctildeJmathcalR_ctildebig(varepsilontildebig)} holds also for $\big|2\big(\widetilde{a}'(t)-\ctilde\big)\psLdeuxLdeux{\mu \mathcal{H}_{\ctilde}\big(\partial_x\varepsilontilde\big)}{\etilde}\big|$. To conclude, we use the bound~\eqref{controle de varepsilon^* par widetildee} so that

\begin{equation*}
     K_* B_*^2 \Big(R_2\normXboule{\etilde}{R_2} +  4\delta  \normX{\etilde}\Big)\normX{\varepsilontilde}\leq K_*B_*^2  D_*\Big(\dfrac{R_2^2}{\delta}\normXboule{\etilde}{R_2}^2+5\delta \normX{\etilde}^2\Big).
\end{equation*}

Fixing the value of $\delta$ in the previous estimate such that $20\delta K_* B_*^2 D_*\leq\Xi_3$ and fixing also $R:=\max(R_1,R_2)$ that only depends on $c^*$ as well, we finally conclude that 
\begin{align*}
    -4\psLdeuxLdeux{\mu\mathcal{H}_{\ctilde}\big(\partial_x\etilde\big)}{\etilde} +2\psLdeuxLdeux{\mu\mathcal{H}_{\ctilde}J\mathcal{R}_{\ctilde}\big(\varepsilontilde\big)}{\etilde}+2\big(\widetilde{a}'(t)-\ctilde\big)\big\langle\mu \mathcal{H}_{\ctilde}\big(\partial_x \varepsilontilde\big),\etilde \big\rangle & _{L^2\times L^2}\geq \dfrac{\Xi_3}{2}\normX{\etilde}^2\\
    &+\grandOde{\normXboule{\etilde}{R}^2}
\end{align*}

\underline{Step 4. Conclusion.} We deduce from Proposition~\ref{proposition: borne inférieure eta_c(x)} and all three previous steps the existence of constants $\Xi_2,\Xi_4$ such that 
\begin{align*}
    n'(t)
    &\geq \dfrac{\Xi_3}{2}\normX{\etilde}^2 -\Xi_4\normXboule{\etilde}{R}^2 + \gamma_*\left(\Xi_1\int_\R \eta_c\left((\partial_x\widetilde{e}_\eta)^2 + \etildeeta^2 + \etildev^2\right) - \Xi_2\alpha_*^\frac{1}{2}\normX{\etilde}^2\right)\\
    &\geq \dfrac{\Xi_3}{2}\normX{\etilde}^2+\int_{\mathcal{B}_R}\big(\Xi_1\gamma_* m_{1}e^{-\iota_* R}-\Xi_4\big)\left((\partial_x\widetilde{e}_\eta)^2 + \etildeeta^2 + \etildev^2\right)  - \Xi_2\gamma_*\alpha_*^\frac{1}{4}\normX{\etilde}^2. 
\end{align*}

We set $\gamma_*$ large enough so that $2\Xi_1\gamma_* m_{1}e^{-\iota_* R}\geq\Xi_4$, which yields to 
\begin{equation*}
    n'(t)\geq \min\left(\dfrac{\Xi_3}{2},\dfrac{\Xi_4}{2}\right)\int_\R \left((\partial_x\widetilde{e}_\eta)^2 + \etildeeta^2 + \etildev^2\right) - \Xi_2\gamma_*\alpha_*^\frac{1}{4}\normX{\etilde}^2.
\end{equation*}

Having fixed the value of $\gamma_*$, we can now assume, up to shrinking the value of $\alpha_*$ given by Theorem~\ref{théorème: stabilité orbitale hydro}, that $2\min\left(\frac{\Xi_3}{2},\frac{\Xi_4}{2}\right)\geq\Xi_2\gamma_*\alpha_*^\frac{1}{4}$. We deduce the existence of $\Xi_* >0$ such that
\begin{equation*}
    n'(t)\geq \Xi_*\normX{\etilde}^2.
\end{equation*}
\end{proof}

We now deal with the proof of Claim~\ref{claim: forme quadratique du terme fastidieux}.
\begin{proof}[Proof of Claim~\ref{claim: forme quadratique du terme fastidieux}]
We recall that the expressions of $\mathcal{H}_c$ and $\mathcal{M}_c$ are given in~\eqref{calcul de mathcal(H_c)} and~\eqref{expression de mathcal M_c}. Also, recall the notation $B_c=1-\eta_c$. Since $M_{c}$ and $S$ are symmetric, we have
\begin{align*}
-4\psLdeuxLdeux{M_{c}S\mathcal{H}_{c}\partial_x\etilde}{\etilde}=& -4\psLdeuxLdeux{\mathcal{H}_{c}\etilde}{SM_{c}\etilde}\\
    =&\quad  4\int_\R\bigg(m_{1,c}\etilde_\eta\partial_x\Big(\dfrac{\partial_x^2\etildeeta}{4B_{c}}\Big)-m_{1,c}\mathcal{M}_{c}\etildeeta\partial_x\etildeeta+\dfrac{c m_{1,c}}{2B_{c}}\etildeeta\partial_x\etildev\\
    &-\dfrac{c^2 m_{2,c}}{4B_{c}}\etildeeta\partial_x\etildeeta+ \dfrac{m_{2,c}B_{c}}{2}\etildeeta\partial_x\etildev +\dfrac{c m_{1,c}}{2B_{c}}\etildev\partial_x\etildeeta+m_{1,c}B_{c}\etildev\partial_x\etildev\bigg).
\end{align*}

By doing integrations by parts, we can express the previous quantity as an integral of a quadratic form that is
\begin{align*}
-4\psLdeuxLdeux{M_{c}S\mathcal{H}_{c}\partial_x\etilde}{\etilde}=&\int_\R\bigg(q_{1,c}\etildev^2+q_{2,c} \etildev\etildeeta+q_{3,c}\etildev\partial_x\etildeeta+ \Big(\dfrac{m_{1,c}'}{B_c}\Big)'\etildeeta\partial_x\etildeeta  + \dfrac{m_{1,c}'}{B_c}(\partial_x\etildeeta)^2 \\
    &+\Big(\dfrac{m_{1,c}}{2B_c}\Big)'(\partial_x\etildeeta)^2-4m_{1,c}\mathcal{M}_c\etildeeta\partial_x\etildeeta- 4m_{2,c}\widetilde{B}_c\etildeeta\partial_x\etildeeta\bigg),
\end{align*}

where, by Proposition~\ref{proposition: signe de q_1,c et q_2,c},
\begin{align}
    q_{1,c}&=-2(a_c B_c)'=\dfrac{2\big((\eta_c')^2-\eta_c (1-\eta_c)\eta_c''\big)}{\eta_c^2}\label{expression développée de q_1,c}>0,\\
    q_{2,c}&=4\big(b_c\widetilde{B}_c+a_cB_c)'=\dfrac{2c\big(\eta_c\eta_c''-(\eta_c')^2\big)}{\eta_c^2},\notag \\
    q_{3,c}&=4b_c B_c=-\frac{2c\eta_c'}{1-\eta_c}.\notag\\
\end{align}

Now applying the Gauss reduction for this quadratic form, we obtain 
\begin{align*}
-4\psLdeuxLdeux{M_{c}S\mathcal{H}_{c}\partial_x\etilde}{\etilde}=&\int_\R q_{1,c}\Big(\etildev+\dfrac{q_{2,c}}{2q_{1,c}} \etildeeta+\dfrac{q_{3,c}}{2q_{1,c}}\partial_x\etildeeta\Big)^2\\
     +\int_\R \bigg(-\dfrac{q_{2,c}^2}{4q_{1,c}}\etildeeta^2-&\dfrac{q_{2,c}q_{3,c}}{2q_{1,c}}\etildeeta\partial_x\etildeeta -\dfrac{q_{3,c}^2}{4q_{1,c}}(\partial_x\etildeeta)^2\\
    +\Big(\dfrac{m_{1,c}'}{B_c}\Big)'\etildeeta\partial_x\etildeeta  +& \dfrac{m_{1,c}'}{B_c}(\partial_x\etildeeta)^2+\Big(\dfrac{m_{1,c}}{2B_c}\Big)'(\partial_x\etildeeta)^2-4m_{1,c}\mathcal{M}_c\etildeeta\partial_x\etildeeta- 4m_{2,c}\widetilde{B}_c\etildeeta\partial_x\etildeeta\bigg),\\
\end{align*}

hence
\begin{align*}
    -4\psLdeuxLdeux{M_{c}S\mathcal{H}_{c}\partial_x\etilde}{\etilde}=& \int_\R q_{1,c}\Big(\etildev+\dfrac{q_{2,c}}{2q_{1,c}} \etildeeta+\dfrac{q_{3,c}}{2q_{1,c}}\partial_x\etildeeta\Big)^2 +\int_\R \widetilde{q_{1,c}}\Big(\partial_x\etildeeta +m_{1,c}\etildeeta\Big)^2 \\
     &+\int_\R \bigg(-2m_{1,c}\widetilde{q_{1,c}}\etildeeta\partial_x\etildeeta - m_{1,c}^2\widetilde{q_{1,c}}\etildeeta^2 -\dfrac{q_{2,c}^2}{4q_{1,c}}\etildeeta^2-\dfrac{q_{2,c}q_{3,c}}{2q_{1,c}}\etildeeta\partial_x\etildeeta  \\
    &+\Big(\dfrac{m_{1,c}'}{B_c}\Big)'\etildeeta\partial_x\etildeeta  -4m_{1,c}\mathcal{M}_c\etildeeta\partial_x\etildeeta- 4m_{2,c}\widetilde{B}_c\etildeeta\partial_x\etildeeta\bigg)
\end{align*}

where
\begin{align*}
    \widetilde{q_{1,c}}&=-\dfrac{q_{3,c}^2}{4q_{1,c}}+\dfrac{m_{1,c}'}{B_c}+\Big(\dfrac{m_{1,c}}{2B_c}\Big)'.
\end{align*}

To finish the proof, we show that 
\begin{align*}
    \int_\R \Big(q_{4,c}\etildeeta\partial_x\etildeeta + q_{5,c}\etildeeta^2 \Big) = 0,
\end{align*}

where 
\begin{align*}
    q_{4,c} &=-2m_{1,c}\widetilde{q_{1,c}}-\dfrac{q_{2,c}q_{3,c}}{2q_{1,c}}+\Big(\dfrac{m_{1,c}'}{B_c}\Big)'-4m_{1,c}\mathcal{M}_c- 4m_{2,c}\widetilde{B}_c,\\
    q_{5,c} & =- m_{1,c}^2\widetilde{q_{1,c}}-\dfrac{q_{2,c}^2}{4q_{1,c}}.
\end{align*}

We compute indeed on the one hand \begin{equation*}
    q_{5,c}=\dfrac{3(\eta_c')^2\eta_c''}{2\eta_c^3B_c}-\dfrac{3(\eta_c')^4}{2\eta_c^4B_c}+\dfrac{(\eta_c')^4}{2\eta_c^3B_c^2}+\dfrac{c^2\eta_c''}{2\eta_c B_c}+\dfrac{c^2(\eta_c')^2}{2}\Big(\dfrac{1}{B_c^2}-\dfrac{1}{\eta_c^2}\Big).
\end{equation*}

On the other hand, in order to compute $q_{4,c}$, we first use the equation satisfied by $\eta_c'''$ (obtained by differentiating equation~\eqref{equation -eta'' = g(eta)}), and notice that
\begin{equation*}
    \Big(\dfrac{m_{1,c}'}{B_c}\Big)'-4m_{1,c}\mathcal{M}_c - 4m_{2,c}\widetilde{B}_c= \dfrac{c^2\eta_c'}{\eta_c B_c^2}+\dfrac{3\eta_c'\eta_c''}{\eta_c^2B_c}-\dfrac{2(\eta_c')^2}{\eta_c^3B_c}+\dfrac{(\eta_c')^3}{\eta_c^2B_c^2},
\end{equation*}

then \begin{equation*}
    q_{4,c}=\dfrac{c^2}{q_{1,c}}\Big(\dfrac{2\eta_c'\eta_c''}{\eta_c B_c}-\dfrac{2(\eta_c')^3}{\eta_c^2B_c^2}+\dfrac{\eta_c' q_{1,c}}{\eta_c B_c^2}\Big)+\dfrac{(\eta_c')^3}{\eta_c^3 B_c},
\end{equation*}

and using the expression of $q_{1,c}$ in~\eqref{expression développée de q_1,c}, we obtain 
\begin{equation*}
    q_{4,c}=\dfrac{(\eta_c')^3}{\eta_c^3 B_c}+\dfrac{c^2\eta_c'}{\eta_c B_c}.
\end{equation*}

Finally, in view of the latter expressions of $q_{4,c}$ and $q_{5,c}$, we obtain
\begin{equation*}
    \int_\R \Big(q_{4,c}\etildeeta\partial_x\etildeeta + q_{5,c}\etildeeta^2 \Big) =\int_\R \Big(-\dfrac{q_{4,c}'}{2}+q_{5,c}\Big)\etildeeta^2 = 0,
\end{equation*}
which concludes the proof.
\end{proof}

Finally, we prove Claim~\ref{claim: spectre essentiel de T_c}.
\begin{proof}[Proof of Claim~\ref{claim: spectre essentiel de T_c}]
Recall from~\eqref{condition suffisante pour la stabilité orbitale sur f''(1)+6f'(1)>0} that $k<0$ and from~\eqref{limite à l'infini de q_1,c/eta_c} that
\begin{equation*}
    k_0=2\nu_c^2 - \dfrac{k}{3}>0.
\end{equation*}
From~\eqref{développement asymptotique de widetilde q_1,c}, we derive 
\begin{equation}\label{limite à l'infini de widetilde q_1,c/eta_c}
    \dfrac{\widetilde{q_{1,c}}(x)}{\eta_c(x)}\underset{|x|\rightarrow +\ii}{\longrightarrow}k_1:=-\left(\dfrac{k}{4}+\dfrac{\nu_c^2}{2}+\dfrac{c^2 \nu_c^2}{k_0}\right)\underset{c\rightarrow c_s}{\sim}-\dfrac{k}{4}>0,
\end{equation}
and from the expressions of $q_{2,c}$ and $q_{3,c}$ in Claim~\ref{claim: forme quadratique du terme fastidieux}, we compute
\begin{equation}\label{limite à l'infini de q_3,c}
    \widetilde{q_{3,c}}\underset{|x|\rightarrow +\ii}{\longrightarrow}k_3:=\dfrac{c}{2}\left(\nu_c^2+\dfrac{k}{3}\right),
\end{equation}

which leads to
\begin{equation}\label{limite vers k_3^2/k_0}
    \dfrac{\eta_c\widetilde{q_{3,c}}^2}{q_{1,c}}\underset{|x|\rightarrow +\ii}{\longrightarrow}\dfrac{k_3^2}{k_0}=\dfrac{c^2}{4k_0}\left(\nu_c^4+\dfrac{2\nu_c^2 k}{3}+\dfrac{k^2}{9}\right).
\end{equation}

From the expression of $\widetilde{q_{1,c}}$ in Claim~\ref{claim: forme quadratique du terme fastidieux},~\eqref{limite à l'infini de widetilde q_1,c/eta_c} and from
\begin{equation*}
    m_{1,c}''\underset{|x|\rightarrow +\ii}{\sim}-\dfrac{k\eta_c'}{6},
\end{equation*}

we infer
\begin{equation}\label{limite à l'infini de partial_x (widetilde q_1,c m_1,c/eta_c}
    -\dfrac{3}{2}\partial_x\left(\dfrac{\widetilde{q_{1,c}}m_{1,c}}{\eta_c}\right)=\grandOde{\eta_c}\underset{|x|\rightarrow +\ii}{\longrightarrow}0.
\end{equation}

On the other hand, from the asymptotic behavior of $\eta_c$ and $\eta_c'$, and from~\eqref{limite à l'infini de widetilde q_1,c/eta_c}, we infer
\begin{equation}\label{limite vers 9/4 k_1}
    \frac{9\widetilde{q_{1,c}}m_{1,c}^2}{4\eta_c}\underset{|x|\rightarrow +\ii}{\longrightarrow}\dfrac{9k_1\nu_c^2}{4}.
\end{equation}

As a consequence of~\eqref{limite vers k_3^2/k_0},~\eqref{limite à l'infini de partial_x (widetilde q_1,c m_1,c/eta_c},~\eqref{limite vers 9/4 k_1} and \eqref{limite à l'infini de widetilde q_1,c/eta_c} once again, we derive
\begin{equation}\label{limite vers k_2}
    \dfrac{\eta_c\widetilde{q_{3,c}}^2}{ q_{1,c}}-\dfrac{3}{2}\partial_x \left(\dfrac{\widetilde{q_{1,c}}m_{1,c}}{\eta_c}\right)+\dfrac{9\widetilde{q_{1,c}}m_{1,c}^2}{4\eta_c}\underset{|x|\rightarrow +\ii}{\longrightarrow} k_2 ,
\end{equation}
where the limit $k_2$ can be shown to satisfy
\begin{equation}\label{expression de k_2}
    k_2 =\dfrac{k_3^2}{k_0}+\dfrac{9k_1\nu_c^2}{4}= -\frac{c^2 k}{12} - c^2 \nu_c^2 - \frac{9 k \nu_c^2}{16} - \frac{9 \nu_c^4}{8}\underset{c\rightarrow c_s}{\sim} - \dfrac{c_s^2 k }{12}>0.
\end{equation}

Combining~\eqref{limite à l'infini de q_1,c/eta_c},~\eqref{limite à l'infini de widetilde q_1,c/eta_c},~\eqref{limite à l'infini de q_3,c} and~\eqref{limite vers k_2}, we can formally take the limit $|x|\rightarrow +\ii$ in the operator $T_c$ and obtain the self-adjoint operator \begin{equation*}
    T_\ii(\htilde)=\begin{pmatrix}
        -k_1\partial_x^2 \htildeeta +k_2 \htildeeta + k_3\htildev \\
        k_3 \htildeeta + k_0\htildev.
    \end{pmatrix}.
\end{equation*}

Since $T_\ii$ is a perturbation of $T_c$, then by the Weyl criterion for self-adjoint operators, we have $\sigma_{ess}(T_c)=\sigma_{ess}(T_\ii)$. In the Fourier domain, we observe that $\lambda\in \sigma (T_\ii)$ if and only if there exists $\xi\in\C$ such that $P(\lambda)=\lambda^2-\lambda (k_1|\xi|^2+k_2+k_0)+(k_1|\xi|^2+k_2)k_0 - k_3^2 =0$. The discriminant of the previous polynomial in $\lambda$ is
\begin{align}
    \Delta_{\xi}&=(k_1|\xi|^2+k_2+k_0)^2-4\big((k_1|\xi|^2+k_2)k_0 - k_3^2\big)\label{expression du déterminant}\\
    &=(k_1|\xi|^2+k_2-k_0)^2+4k_3^2 \geq 0 \notag.
\end{align} 

We infer that the polynomial $P$ has two roots
\begin{equation*}
    \lambda_\pm(\xi) = \dfrac{k_1|\xi|^2+k_2+k_0 \pm\sqrt{\Delta_\xi}}{2}.
\end{equation*}

Clearly, $\lambda_-(\xi)<\lambda_+(\xi)$, and we now check by some standard real analysis that the function $\xi\mapsto\lambda_-(\xi)$ increases on $\R_+$ and thus $\lambda_-(0)\leq\lambda_-(\xi)$. Indeed, making the change of variables $y=k_1|\xi|^2 + k_2 + k_0$ leads us to study the variations of the function 
$$h(y):=y - \sqrt{(y-2k_0)^2+4k_3^2}.$$

By shrinking the value of $c_1$ so that both quantities $k_3^2$ and $k_1$ are positive for $c\in (c_1,c_s)$ by~\eqref{signe de k_2k_0 - k_3^2} and~\eqref{limite à l'infini de widetilde q_1,c/eta_c}, we infer in particular that $h$ is well-defined and differentiable on $[k_2+k_0,+\ii)$. In addition, we compute
\begin{equation*}
    h'(y)=1-\dfrac{y-2k_0}{\sqrt{(y-2k_0)^2+4k_3^2}},
\end{equation*}
and we have $0\leq (y-2k_0)^2\leq (y-2k_0)^2+4k_3^2$, then passing to square root function, we deduce
\begin{equation*}
    \dfrac{y-2k_0}{\sqrt{(y-2k_0)^2+4k_3^2}}\leq \dfrac{|y-2k_0|}{\sqrt{(y-2k_0)^2+4k_3^2}}\leq 1.
\end{equation*}

Thus $h'\geq 0$, so that $h(y)\geq h(k_2+k_0)=\lambda_-(0)$. From~\eqref{expression du déterminant}, and the fact that $k_2>0$ and $k_0>0$, we deduce that
\begin{equation*}
    \spess(T_\ii)\subset [\tau_c,+\ii),
\end{equation*}
where we have set \begin{equation*}
    \tau_c=\lambda_-(0)=\dfrac{k_2+k_0-\sqrt{(k_2-k_0)^2+4k_3^2}}{2}=\dfrac{k_2+k_0}{2}\left(1-\sqrt{1+\frac{4(k_3^2-k_0k_2)}{(k_0+k_2)^2}}\right).
\end{equation*}

Finally, by~\eqref{expression de k_2} we deduce that 
\begin{align}\label{signe de k_2k_0 - k_3^2}  
    k_2k_0-k_3^2=\dfrac{9k_0k_1\nu_c^2}{16}
    &\underset{c\rightarrow c_s}{\sim}\dfrac{3k^2\nu_c^2}{16}>0,
\end{align}
which implies that $\tau_c$ is strictly positive with respect to small speeds $c$ and that
\begin{equation*}
    \tau_c\underset{c\rightarrow c_s}{\sim}\dfrac{k_0k_2-k_3^2}{k_0+k_2}\underset{c\rightarrow c_s}{\sim}-\dfrac{k\nu_c^2(c_s^2 + 4)}{4}.
\end{equation*}
\end{proof}

\appendix

\section{Equivalence of the topologies in a neighborhood of the travelling wave}

In this section, we show the correspondence between the topologies for the original and the hydrodynamical framework, at least in a neighborhood of the travelling wave of speed $c^*$. We recall that the metric $d$ is given in~\eqref{métrique d}.

\begin{lem}\label{lemme: équivalence des distances classique et hydro proche d'un soliton}
    There exists a ball $\mathcal{B}_*$ of radius $r_*$ and centered in $\gu_{c^*}$ for the metric $d$ and a positive constant $\widetilde{A}_*$ such that for any $u\in\mathcal{B}_*$, $\normX{Q-Q_{c^*}}\leq \widetilde{A}_* d(u,\gu_{c^*})$.
\end{lem}
\begin{proof}
    First, we have the bound 
\begin{align*}
    \normLii{u}\lesssim \normHun{1-|u|}+1&\lesssim 1 +\normLdeux{1-|u|}+\normLdeux{\partial_x|u|}\leq 1+\normLdeux{1-|u|^2}+\normLdeux{\partial_x u}\\
    &\leq 1+2r_*+\normLdeux{1-|\gu_{c^*}|^2}+\normLdeux{\partial_x \gu_{c^*}}.
\end{align*}
    
Now using this bound, we compute \begin{align}
        \normHun{\eta-\eta_{c^*}}
        & \lesssim \normLdeux{|u|^2-|\gu_{c^*}|^2}+\normLdeux{u.\partial_x u - \gu_{c^*}.\partial_x \gu_{c^*}}\notag\\ \label{estimée du terme u.partial_x u - gu_gc.partial_x gu_gc}
        & \lesssim r_* +\normLii{u}r_* + \normLdeux{\partial_x \gu_{c^*}.(u-\gu_{c^*})}. 
    \end{align}

By the Sobolev embedding theorem, we can define the pointwise value of functions lying in the energy set $\mathcal{X}(\R)$, so that we can write
\begin{align*}
    (u-\gu_{c^*})(x)&=\int_0^x \partial_x(u-\gu_{c^*}) + (u-\gu_{c^*})(0).
\end{align*}

We then deduce that \begin{align*}
    \normLdeux{\partial_x \gu_{c^*}.(u-\gu_{c^*})} & \leq\normLdeux{ \partial_x\gu_{c^*} \sqrt{|x|} }\normLdeux{\partial_x u-\partial_x\gu_{c^*}}+ \normLdeux{\partial_x\gu_{c^*}}\normLiiunun{u-\gu_{c^*}},
\end{align*}

and since we have exponential decay of $\partial_x\gu_{c^*}$, we thus have 
\begin{equation}\label{normHun eta-eta_c leq d(u,u_c)}
    \normHun{\eta-\eta_{c^*}}\lesssim d(u,\gu_{c^*}).
\end{equation}
    
We now turn to the norm $\normLdeux{v-v_{c^*}}$. Since $\gu_{c^*}\in\Nenergyset$, we have a bound from below for this function, and we have just controlled $\normLii{\eta-\eta_{c^*}}$. Up to shrinking the value of $r_*$ we can therefore assume that there exists $m_{c^*} >0$ such that $\inf_\R|u| \geq m_{c^*}$.
    
We conclude by estimating the last norm as \begin{align*}
        \normLdeux{v-v_{c^*}}& \leq \normLdeux{\dfrac{iu.\partial_x u}{|u|^2}-\dfrac{i\gu_{c^*}.\partial_x \gu_{c^*}}{|\gu_{c^*}|^2}}\leq \normLdeux{\dfrac{u.\partial_x u-\gu_{c^*}.\partial_x\gu_{c^*}}{|u|^2}}+\normLdeux{\gu_{c^*}.\partial_x \gu_{c^*}\left(\dfrac{1}{|u|^2}-\dfrac{1}{|\gu_{c^*}|^2}\right)}\\
        & \leq \dfrac{1}{m_{c^*}^2}\normLdeux{u.\partial_x u-\gu_{c^*}.\partial_x\gu_{c^*}}+\dfrac{1}{m_{c^*}^2 \inf_\R|\gu_{c^*}|^2}\normLdeux{\gu_{c^*}.\partial_x\gu_{c^*}}\normLii{\eta-\eta_{c^*}}.
    \end{align*}

The first term in the second-hand side can be controlled similarly as in~\eqref{estimée du terme u.partial_x u - gu_gc.partial_x gu_gc}. As for the second one, we use the one-dimensional Sobolev embedding and~\eqref{normHun eta-eta_c leq d(u,u_c)}, so that we eventually have
\begin{equation*}
     \normLdeux{v-v_{c^*}}\lesssim  d(u,\gu_{c^*}).
\end{equation*}
\end{proof}

\section{Weak-continuity of the classical and hydrodynamical flow}\label{appendix: flow maps and weak continuity}

In this section, we give brief details regarding the Cauchy problem of the equations involved in this paper and how this can be combined to prove the weak continuity of the flow of~\eqref{NLS} (resp.~\eqref{NLShydro}) for the topology given by $d$ (resp. $\normX{.}$). This section relies mainly on two articles concerning the Cauchy problem where, on the one hand~\cite{Gallo3} the global well-posedness of~\eqref{NLS} is proved for dimensions $N\leq 4$ (see also~\cite{Gerard1}), and another article~\cite{Gallo1} where the Cauchy problem is addressed on Zhidkov spaces of order $l$ (see~\eqref{définition Zhitkov space}) for all dimensions $l>\frac{N}{2}$. We also mention that in~\cite{AntHieMar}, the Cauchy problem for~\eqref{NLS} has been handled for dimensions $N\in\{2,3\}$ under less restrictive Kato-type assumptions.

We recall that we have assumed both the following conditions throughout the article:
\begin{itemize}
    \item For all $\rho\in\R$,
\begin{equation}\label{hypothèse de croissance sur F minorant intermediaire}\tag{H1}
    \dfrac{c_s^2}{4} (1-\rho)^2 \leq F(\rho).
\end{equation}
    \item There exist $M\geq 0$ and $q_*\in [2, +\ii)$ such that for all $\rho \geq 2$,
\begin{equation}\label{hypothèse de croissance sur F majorant}\tag{H2}
    F(\rho)\leq M|1-\rho|^{q_*}.
\end{equation}
\end{itemize}

\begin{thm}[Theorem 1.2 in~\cite{Gallo1}]
Let $u_0\in \energyset$. Take $f$ in $ C^2(\R)$ satisfying \eqref{hypothèse de croissance sur F minorant intermediaire}. In addition, assume that there exist $\alpha_1\geq 1$ and $C_0 >0$ such that for all $\rho\geq 1$,

\begin{equation*}
|f''(\rho)|\leq \dfrac{C_0}{\rho^{3-\alpha_1}}.
\end{equation*}

If $\alpha_1 > \frac{3}{2}$, assume moreover that there exists $\alpha_2\in [ \alpha_1-\frac{1}{2},\alpha_1]$ such that for $\rho\geq 2$, $C_0 \rho^{\alpha_2} \leq F(\rho)$.\\
There exists a unique function $w \in  C^0\big(\R,H^1(\R)\big)$ such that $u:=u_0+w$ solves \eqref{NLS}. Moreover, the solution depends continuously on the initial condition, and the energy $E$ and the momentum $p$ are conserved by the flow.
\end{thm}

From the previous global well-posedness result, we derive the weak-continuity for the metric $d$.
\begin{prop}\label{proposition: weak continuity of the flow in original setting}
Let $(\Psi_{n,0})_n\in \mathcal{NX}(\R)$ and $\Psi_{0}\in\Nenergysethydrokk{\!}$ such that 
    \begin{equation}\label{Psi_n,0 ntendfd Psi_0}
        \Psi_{n,0}\ntendfd \Psi_0.
    \end{equation}
Write $\Psi_n$ and $\Psi$ the associated solutions to~\eqref{NLS} given by Theorem~\ref{théorème: local global well-posedness of cauchy problem}. Then for any $t\in\R$,
\begin{equation*}
    \Psi_n(t)\ntendfd \Psi(t).
\end{equation*}
\end{prop}

\begin{proof}
We can follow the main lines of the proof of Proposition A.3 in~\cite{BetGrSm2}. The only remaining details that are to be verified are the compactness argument exhibiting the limit profile satisfying~\eqref{NLS} in the distributional sense and the fact that the limit profile is sufficiently smooth. The construction of the limit profile is based on the boundedness of the sequence of the energies $\big(E(\Psi_{n,0})\big)_n$ in the classical setting. Combining all three convergences in~\eqref{Psi_n,0 ntendfd Psi_0} provides that $(1-|\Psi_{n,0}|^2)_n$ is bounded in $L^2(\R)$ and $\big(\partial_x(1-|\Psi_{n,0}|^2)\big)_n$ as well. The first one is straightforward and to prove the second one, we first write $\partial_x(1-|\Psi_{n,0}|^2)=-2\partial_x\Psi_{n,0}.\Psi_{n,0}$. Now, we set $\chi\in C^\ii_c(\R)$ supported on a compact $K$ and real valued. Using both the weak-convergence of $(\partial_x\Psi_{n,0})_n$ and the convergence in $L^\ii_\loc$ yields to
\begin{equation*}
    \psLdeux{\partial_x\Psi_{n,0}.\Psi_{n,0}}{\chi}=\int_K \partial_x\Psi_{n,0}.(\Psi_{n,0}\chi)\underset{n\rightarrow +\ii}{\longrightarrow}\int_K\partial_x \Psi_{0}.(\Psi_{0}\chi)=\psLdeux{\partial_x\Psi_{0}.\Psi_{0}}{\chi}.
\end{equation*}
    
This proves that $(1-|\Psi_{n,0}|^2)_n$ is bounded in $H^1(\R)$, by some bound denoted by $M_3$. We also know that the sequence of kinetic energies $\big(E_k(\Psi_{n,0})\big)_n$ is bounded, and we still denote $M_3$ such a bound. Regarding the potential energy, we deduce from hypothesis~\eqref{hypothèse de croissance sur F majorant} that there exist two positive constants $M_1$ and $M_2$ such that for any $\rho\in\R$
    \begin{equation*}
        F(\rho)\leq M_1 (1-\rho)^2 + M_2|1-\rho|^{q_*} .
    \end{equation*}
Therefore, by the one-dimensional Sobolev embedding $H^1(\R)\hookrightarrow L^\ii(\R)$, up to taking a larger $M_3$,
\begin{equation*}
    \int_\R F(|\Psi_{n,0}|^2)\leq M_1 \int_\R \big(1-|\Psi_{n,0}|^2\big)^2+M_2\int_\R \big(1-|\Psi_{n,0}|^2\big)^{q_*}\leq M_3 .
\end{equation*}

Then, following the path of the anterior proof, we deduce the existence of a profile $\Phi$ characterized as the weak limit function in $L^\ii\big([0,T],L^\ii(\R)\big)$ satisfying for any $t\in [0,T]$,
\begin{equation*}
    \partial_x\Psi_n(t)\ntendfLdeux\partial_x \Phi(t) \quad\text{ and }\quad  1-|\Psi_n(t)|^2 \ntendfLdeux 1-|\Phi(t)|^2,
\end{equation*}
and satisfying also~\eqref{NLS} with initial data $\Psi_0$ in a distributional sense due again to the Sobolev embedding.

It remains to check that it lies in $C^0\big([0,T],\mathcal{Z}^s(\R)\big)$ where $T$ is a given positive number and $s\in (\frac{1}{2},1)$ designates a smoothness parameter in the Zhidkov space defined by
\begin{equation}\label{définition Zhitkov space}
    \mathcal{Z}^s(\R):=\left\{\psi\in L^\ii(\R)\big|\partial_x\psi\in H^{s-1}(\R)\right\}.
\end{equation}

First, the $L^2$-norms of $\partial_x\Phi$ and $1-|\Phi|^2$ are uniformly bounded with respect to $t\in [0,T]$ and in particular, $(\partial_x\Phi,1-|\Phi|^2)\in L^\ii\big([0,T],L^2(\R)\big)\times L^\ii\big([0,T],H^1(\R)\big)$. Now, since $\Phi$ satisfies~\eqref{NLS} in a distributional sense, we have
\begin{equation}\label{equation d'évolution de partial_x Phi}
    i\partial_t(\partial_x\Phi)=-\partial_x^3\Phi - \partial_x\big(\Phi f(|\Phi|^2) \big).
\end{equation}

Since $\Phi\in L^\ii\big([0,T],L^\ii(\R)\big)$ and $1-|\Phi|^2 \in L^\ii\big([0,T],H^1(\R)\big)$, and since $f$ and $f'$ are continuous functions with $f(1)=0$, we proceed to a first-order Taylor expansion between $|\Phi|^2$ and $1$, in order to provide
\begin{equation}
    \Phi f(|\Phi|^2) = -\Phi (1-|\Phi|^2) \int_0^1 (1-s)f'(s|\Phi|^2)ds \in L^\ii\big([0,T],L^2(\R)\big)
\end{equation}

so that from~\eqref{equation d'évolution de partial_x Phi}, we have $\partial_x\Phi\in W^{1,\ii}\big(|0,T],H^{-2}(\R)\big)$. Using the same argument than in~\cite{BetGrSm2}, we identically infer that $1-|\Phi|^2\in W^{1,\ii}\big(|0,T],H^{-1}(\R)\big)$ and that $\Phi\in C^0\big([0,T],\mathcal{Z}^s(\R)\big)$ as well.

In~\cite{Gallo3}, it is proved that under the condition~\eqref{hypothèse de croissance sur F minorant intermediaire}, the Cauchy problem associated with~\eqref{NLS} is shown to be locally well-posed in the Zhidkov spaces $\mathcal{Z}^s(\R)$ for any $s\in\N^*$ but the proof generalizes to all $s>\frac{1}{2}$ by the one-dimensional Sobolev embedding. The fact that $\Phi$ lies in $\mathcal{Z}^s(\R)$ with $\frac{1}{2}<s<1$ and $\Phi(0)=\Psi_0$ thus implies that $\Phi$ and $\Psi$ coincide, which concludes the proof.
\end{proof}

Now, we state a local well-posedness result in the hydrodynamical framework, proved in~\cite{Gallo3} and extendable to all integer $l>\frac{1}{2}$ according to the remark above Theorem~5.1 in~\cite{Gallo3}.
\begin{thm}[Gallo~\cite{Gallo3}]\label{théorème: local well posedness in hydro}
    We assume that $f\in C^{5}(\R_+)$ is such that for all $\rho\in\R$,
\begin{equation}\label{hypothèse de croissance sur F minorant intermediaire}\tag{H1}
    \dfrac{c_s^2}{4} (1-\rho)^2 \leq F(\rho).
\end{equation}
Let $l\in\{0,...,4\}$ be a natural integer and let $(\eta_0,v_0)\in \Nenergysethydrokk{l}$. There exist $T_{\max} >0$ and a unique solution $(\eta,v)\in C^0\big([0,T_{\max} ),\Nenergysethydrokk{l}\big)$ to~\eqref{NLShydro} with initial datum $(\eta_0,v_0)$. The maximal time $T_{\max}$ is continuous with respect to the initial datum and is characterized by
    \begin{equation*}
        \lim_{t\rightarrow T_{\max}^-}\sup_{x\in\R}\eta(t,x) = 1.
    \end{equation*}

Moreover, the flow map is continuous on $\Nenergysethydrokk{l}$ and the energy and the momentum are conserved along the flow.
\end{thm}

From the local well-posedness in Theorem~\ref{théorème: local well posedness in hydro} and the weak continuity of the flow in the original setting in Proposition~\ref{proposition: weak continuity of the flow in original setting}, we derive the weak-continuity of the flow in the hydrodynamic framework. The proof of this statement can be found in~\cite{BetGrSm2}.

\begin{prop}\label{proposition: weak continuity of the flow in hydrodynamical setting}
    Let $(Q_{n,0})_n\in\Nenergysethydrokk{\!}$ and $Q_{0}\in\Nenergysethydrokk{\!}$ such that 
    \begin{equation*}
        Q_{n,0}\ntendfX Q_0.
    \end{equation*}
Assume in addition that there exists a constant $M\in |0,1)$ such that the associated solutions $Q_n\in C\big([0,T_{\max,n}),\Nenergysethydrokk{\!}\big)$ and $Q\in C\big([0,T_{\max}),\Nenergysethydrokk{\!}\big)$ to~\eqref{NLShydro} given by Theorem~\ref{théorème: local well posedness in hydro} satisfy for some $0<T<T_{\max,n}$ and any $(n,t,x)\in \N\times [-T,T]\times \R$,
\begin{equation*}
    \eta_n(t,x)\leq M .
\end{equation*}

Then $0 < T < T_{\max}$ and for any $t\in [-T,T]$,
\begin{equation*}
    Q_n(t)\ntendfX Q(t).
\end{equation*}
\end{prop}

Finally, we conclude this section by observing that the proof of Proposition~\ref{proposition: convergence faible vers le profil limite le long de l'évolution} that is made in Subsection A.1.1 of~\cite{BetGrSm2} does not depend on the nonlinearity considered but only on the weak-continuity of the hydrodynamic flow as stated in Proposition~\ref{proposition: weak continuity of the flow in hydrodynamical setting}. In particular, we refer to the latter article for the details of the proof of Proposition~\ref{proposition: convergence faible vers le profil limite le long de l'évolution}.

\section{Estimates on the linearized operator $\mathcal{H}_c$}\label{appendix: propriété du linéarizé}
\subsection{Useful estimates}

We recall the expression of the linearized operator, that is 
\begin{align}\label{calcul de mathcal(H_c)}
    \mathcal{H}_c=\begin{pmatrix} \mathcal{L}_c & -\dfrac{c}{2(1-\eta_c)}\\
    -\dfrac{c}{2(1-\eta_c)} & 1-\eta_c
    \end{pmatrix},
\end{align}
with
\begin{equation*}
    \mathcal{L}_c(\varepsilon_\eta)=-\Big(\dfrac{\varepsilon_\eta'}{4(1-\eta_c)}\Big)'+ \mathcal{M}_c(\varepsilon_\eta),\\
\end{equation*}

where
\begin{equation}\label{expression de mathcal M_c}
    \mathcal{M}_c(\varepsilon_\eta)=-\left(\dfrac{\eta_c''}{4(1-\eta_c)^2}+\dfrac{(\eta_c')^2}{4(1-\eta_c)^3}+\dfrac{f'(1-\eta_c)}{2}\right)\varepsilon_\eta .
\end{equation}

The following proposition states several crucial estimates satisfied by the operator $\mathcal{H}_c$. Since we need to have sharp estimates, we introduce for non-negative integers $l$, the set
\begin{equation*}
    \mathcal{Y}_\rho^{l}:=\left\lbrace Q=(\eta,v)\in H^{l+2}(\R)\times H^l(\R) \Big|  \normYpoidsx{Q}{\rho}{l} < +\ii\right\rbrace,
\end{equation*}
where \begin{equation*}
    \normYpoidsx{Q}{\rho}{l}^2 = \sum_{m=0}^{l+2}\int_\R \big(\partial_x^m\eta(x)\big)^2|x|^\rho dx+\sum_{m=0}^{l}\int_\R \big(\partial_x^m v(x)\big)^2|x|^\rho dx.
\end{equation*}

\begin{prop}\label{proposition: propriété du linéarisé mathcal H_c}
There exists an interval $J$ centered in $c^*$ such that for any integer $l\in\{-1,...,2\}$ and $m\in\{0,1\}$, there exists a constant $K_*>0$ only depending on $c^*,l$ and $m$ such that for any $\rho\geq 0,c\in J$ and $\varepsilon\in \dom(\mathcal{H}_c)\cap\mathcal{Y}^{l+1}_\rho$ such that $\partial_c\varepsilon\in \dom(\mathcal{H}_c)\cap\mathcal{Y}^{l+1}_\rho$, 
\begin{equation*}
\normXpoidsx{\partial_c^m\mathcal{H}_{c}(\varepsilon)}{\rho}{l}\leq K_*\left(\normYpoidsx{\varepsilon}{\rho}{l+1}+\normYpoidsx{\partial_c^m\varepsilon}{\rho}{l+1}\right).
    \end{equation*}
\end{prop}

\begin{proof}
The previous estimates directly follows from the expression of $\mathcal{H}_c$ and the exponential decay of the travelling wave and its derivative with respect to $c$ in~\eqref{estimée décroissance exponentielle à tout ordre pour eta_c et v_c} combined with Proposition~\ref{proposition: borne uniforme sur les quantités nu_c etc}.
\end{proof}

In particular, we deduce the useful estimates stated in the following corollary.
\begin{cor}\label{corollaire: controle de la norme de mathcal H_c J mathcal R_c}
Up to taking a larger $K_*>0$ in the previous proposition, we have for any $\rho\geq 0$,
\begin{equation*}
        \normXpoidsx{\mathcal{H}_{\ctilde}J\mathcal{R}_{\ctilde}\big(\varepsilontilde\big)}{{\rho}}{-1}\leq K_{*}\normXpoidsx{\varepsilontilde}{\rho}{3}^{2},
    \end{equation*}
    and
    \begin{equation*}
        \normxxk{\mathcal{H}_{\ctilde}(\etilde)}{-1}\leq K_* \normxxk{\varepsilontilde}{3}
    \end{equation*}
\end{cor}

\begin{proof}
To simplify the notation and because the argument is similar otherwise, we only deal with the case $\rho=0$. In view of the definition of~\eqref{définition mathcal R_c}, we can write
\begin{align*}
    \mathcal{R}_{\ctilde}(\varepsilontilde)=\int_0^1 (1-s)\nabla^3 E(Q_{\ctilde}+s\varepsilontilde)(\varepsilontilde,\varepsilontilde,.)ds,
\end{align*}
where
\begin{align*}
\nabla^3 &E(Q_{\ctilde}+s\varepsilontilde)(\varepsilontilde,\varepsilontilde,.)=\\
&\begin{pmatrix}
\dfrac{3(\partial_x\varepsilontilde_\eta)^2}{4(1-\eta_{\ctilde}-s\varepsilontilde_\eta)^2}+\dfrac{3\partial_x(\eta_{\ctilde}+s\varepsilontilde_\eta)\partial_x\varepsilontilde_\eta \varepsilontilde_\eta}{2(1-\eta_{\ctilde}-s\varepsilontilde_{\eta})^3}+\Big(\dfrac{3\big(\partial_x(\eta_{\ctilde}+s\varepsilontilde_\eta)\big)^2}{4(1-\eta_{\ctilde}-s\varepsilontilde_\eta)^4}+\dfrac{f''(1-\eta_{\ctilde}-s\varepsilontilde_\eta)}{2}\Big)\varepsilontilde_\eta^2 - \varepsilontilde_v^2\\
        -2\varepsilontilde_\eta\varepsilontilde_v
    \end{pmatrix}.
\end{align*}

In view of the expression of the operator $J$ in~\eqref{expression de J}, we obtain
\begin{align}
&J\nabla^3 E(Q_{\ctilde}+s\varepsilontilde)(\varepsilontilde,\varepsilontilde,.)=\label{expression de J nabla^3 E(Q+varepsilon)}\\
2\partial_x &\begin{pmatrix}
    2\varepsilontilde_\eta\varepsilontilde_v\\
    \dfrac{-3(\partial_x\varepsilontilde_\eta)^2}{4(1-\eta_{\ctilde}-s\varepsilontilde_\eta)^2}-\dfrac{3\partial_x(\eta_{\ctilde}+s\varepsilontilde_\eta)\partial_x\varepsilontilde_\eta \varepsilontilde_\eta}{2(1-\eta_{\ctilde}-s\varepsilontilde_{\eta})^3}-\Big(\dfrac{3\big(\partial_x(\eta_{\ctilde}+s\varepsilontilde_\eta)\big)^2}{4(1-\eta_{\ctilde}-s\varepsilontilde_\eta)^4}+\dfrac{f''(1-\eta_{\ctilde}-s\varepsilontilde_\eta)}{2}\Big)\varepsilontilde_\eta^2 +\varepsilontilde_v^2
    \end{pmatrix}\notag.
\end{align}

Using the one-dimensional Sobolev embedding, the Cauchy-Schwarz inequality and the uniform bounds on $Q_{\ctilde}$ and $\normX{\varepsilontilde}$ provided by~\eqref{controle de widetilde varepsilon et ctilde -c_0^* par alpha_*}, we conclude that
 \begin{equation*}
        \normLdeuxLdeux{{\mathcal{H}_{\ctilde}J\mathcal{R}_{\ctilde}\big(\varepsilontilde\big)}}\leq K_{*}\normxxk{\varepsilontilde}{3}^2.
    \end{equation*}

To finish, since $\etilde=S\mathcal{H}_{\ctilde}(\varepsilontilde)$, we can apply Proposition~\ref{proposition: propriété du linéarisé mathcal H_c} twice with $m=0$ and use the continuous embedding $\mathcal{X}_\rho^{l+1}\hookrightarrow \mathcal{Y}_\rho^l$ to eventually deduce the second desired estimate.
\end{proof}

\begin{cor}\label{corollaire: controle de H_c J R_c varepsilon par norme W24,3}
Up to taking a larger $K_*>0$ in the previous proposition, we have
 \begin{equation*}
        \normxxk{\mathcal{H}_{\ctilde}J\mathcal{R}_{\ctilde}\big(\varepsilontilde\big)}{-1}\leq K_{*}\normX{\varepsilontilde}^{\frac{5}{4}}\left(\normWkcarre{\varepsilontilde}{12}{1}^\frac{3}{4}+\normWkcarre{\varepsilontilde}{12}{3}^\frac{5}{4}\right).
    \end{equation*}
\end{cor}

\begin{proof}
By Proposition~\ref{proposition: propriété du linéarisé mathcal H_c} and coming back to the expression~\eqref{expression de J nabla^3 E(Q+varepsilon)}, we compute 
\begin{align*}
    \normxxk{\mathcal{H}_{\ctilde}J\mathcal{R}_{\ctilde}\big(\varepsilontilde\big)}{-1}&\leq K_{*}\normYkk{J\mathcal{R}_{\ctilde}(\varepsilontilde)}{0}.
\end{align*}

This norm can be dealt with essentially the same way than in~\cite{BetGrSm2}. As a matter of example, we handle the term involving the nonlinearity $f$. In particular, by exponential decay~\eqref{estimée décroissance exponentielle à tout ordre pour eta_c et v_c}, by Hölder's inequality, Proposition~\ref{proposition: borne uniforme sur les quantités nu_c etc} and the fact that $f''$ and $f'''$ are continuous, we have for $s\in [0,1]$,
\begin{align*}
    \normLdeux{\partial_x\Big(f''(1-\eta_{\ctilde}-s\varepsilontilde_\eta)\varepsilontilde_\eta^2\Big)}\leq  &  \normLdeux{ \varepsilontilde_\eta^2 f'''(1-\eta_{\ctilde}-s\varepsilontilde_\eta)\partial_x\eta_{\ctilde}}+s\normLdeux{\varepsilontilde_\eta^2\partial_x\varepsilontilde_\eta f'''(1-\eta_{\ctilde}-s\varepsilontilde_\eta)}\\
    &+2\normLdeux{\varepsilontilde_\eta \partial_x\varepsilontilde_\eta f''(1-\eta_{\ctilde}-s\varepsilontilde_\eta)}\\
    \leq & \ I_* \left(K_d \normLdeux{\varepsilontilde_\eta^2} + \normLdeux{\varepsilontilde_\eta^2 \partial_x\varepsilontilde_\eta}+2\normLdeux{\varepsilontilde_\eta \partial_x\varepsilontilde_\eta}\right)\\
    \leq  & \ I_*\left(K_d\left\Vert \varepsilontilde_\eta\right\Vert_{L^4}^2 + \left\Vert \varepsilontilde_\eta\right\Vert_{L^6}^2\left\Vert \partial_x\varepsilontilde_\eta\right\Vert_{L^6} + 2\left\Vert \varepsilontilde_\eta\right\Vert_{L^4}\left\Vert \partial_x\varepsilontilde_\eta\right\Vert_{L^4} \right).
\end{align*}

Eventually, note that the maximum order of differentiation is smaller than $3$ and that the number of factors in the $L^2$-norms never exceed $3$, so that we obtain, up to a larger $K_*$,
\begin{equation}\label{normY JRvarepsilon par norme W3,4}
    \normYkk{J\mathcal{R}_{\ctilde}(\varepsilontilde)}{0}\leq K_*\left(\normWkcarre{\varepsilontilde}{3}{4}^2 + \normWkcarre{\varepsilontilde}{3}{6}^3\right).
\end{equation}

Using the general estimate for $g\in \bigcap _{\substack{m\in\{0,...,12\}\\ p\in\{2,3\}}}W^{m,p}(\R)$,
\begin{equation*}
    \normWkp{g}{m}{p}^{\frac{p}{2}}\leq \normLdeux{g}^{\frac{1}{2}}\normWkp{g}{2m}{2}^{\frac{p-1}{2}},
\end{equation*}

and iterating this estimate read
\begin{equation*}
    \normWkp{g}{m}{p}^{\frac{p}{2}}\leq \normLdeux{g}^\frac{p+1}{4}\normWkp{g}{4m}{2}^{\frac{p-1}{4}}.
\end{equation*}

Using the latter control in~\eqref{normY JRvarepsilon par norme W3,4} leads to
\begin{align*}
    \normYkk{J\mathcal{R}_{\ctilde}(\varepsilontilde)}{0}&\leq K_*\left(\normX{\varepsilontilde}^\frac{5}{4}\normWkcarre{\varepsilontilde}{12}{1}^\frac{3}{4} +\normX{\varepsilontilde}^\frac{7}{4} \normWkcarre{\varepsilontilde}{12}{3}^\frac{5}{4}\right).
\end{align*}

Assuming without loss of generality that $\normX{\varepsilontilde}\leq 1$, we obtain the desired estimate.
\end{proof}

\subsection{Improvement of the orbital stability result}
In this subsection, we refine the orbital stability result in~\cite{Bert2} so that we obtain~\eqref{normXQ(t)-Q_c^*,a(t)^2+normR a'(t)-c(t)^2+normR c'(t)leq A_* alpha_0^2}. We recall that in the previous article~\cite{Bert2}, we have a less restrictive control on the modulation parameters that is
\begin{equation}\label{normXQ(t)-Q_c^*,a(t)^2+normR a'(t)-c(t)+normR c'(t)leq A_* alpha_0}
    \normR{ a'(t)-c(t)}+\normR{ c'(t)}=\grandOde{ \normX{\varepsilon(t)}}.
\end{equation}

From this control, we shall derive the quadratic control on $c'(t)$ stated in~\eqref{normXQ(t)-Q_c^*,a(t)^2+normR a'(t)-c(t)^2+normR c'(t)leq A_* alpha_0^2}. We differentiate the second equality in~\eqref{condition d'orthogonalité sur varepsilon dans décomposition orthogonale} with respect to the time. Using the equation satisfied by $\varepsilon$ that is~\eqref{equation sur variable varepsilon} and the facts that $\nabla p(Q_c)(J\mathcal{H}_c\varepsilon)=\nabla p(Q_c)(\partial_x Q_c)=0$,  provides
\begin{align*}
    c'(t)\dfrac{d}{dc}\Big(p(Q_{c(t)}\Big)=\big(a'(t)-c(t)\big)\nabla &p(Q_{c(t)}).\partial_x\varepsilon(t)+c'(t)\nabla^2 p(Q_{c(t)}).(\partial_c Q_{c(t)},\varepsilon (t))\\
    &+\nabla p(Q_{c(t)})(J\mathcal{R}_{c(t)}\varepsilon(t)).
\end{align*}

Moreover, by~\eqref{normXQ(t)-Q_c^*,a(t)^2+normR a'(t)-c(t)+normR c'(t)leq A_* alpha_0}, we can bound $\big(a'(t)-c(t)\big)\nabla p(Q_{c(t)}).\partial_x\varepsilon(t)$ and $c'(t)\nabla^2 p(Q_c).(\partial_c Q_{c(t)},\varepsilon (t))$ by $\grandOde{\normX{\varepsilon(t)}^2}$. Finally, coming back to the computations made in the proof of Corollary~\ref{corollaire: controle de la norme de mathcal H_c J mathcal R_c}, we can also bound $\nabla p(Q_{c(t)})(J\mathcal{R}_{c(t)}\varepsilon(t))$ by $\grandOde{\normX{\varepsilon(t)}^2}$. Using Proposition~\ref{proposition: borne uniforme sur les quantités nu_c etc} provides the desired quadratic control of $|c'(t)|$.

\section{The transonic regime $c\rightarrow c_s$.}\label{appendix: transonic régime}
In this section, we highlight the properties of the travelling waves in the transonic regime. We recall that one of the fundamental assumption that appears in~\cite{Bert1,Bert2} is that  \begin{equation}
    k:=2f''(1)+6f'(1)\neq 0.
\end{equation}

Recall from~\cite{Bert2} that for $(x,c)\in\R\times (c_0,c_s)$, we have
\begin{equation}\label{eta_c(x) vit dans (0,1)}
    0<\eta_c(x)\leq \xi_c <1,
\end{equation}
where 
\begin{equation}\label{équivalent de xi_c}
    \xi_c=\eta_c(0)\underset{c\rightarrow c_s}{\sim}-\dfrac{3\nu_c^2}{k} .    
\end{equation}

Due to both last statements, we recover the condition
\begin{equation}\tag{H3}
    k:=2f''(1)+6f'(1)< 0.
\end{equation}

The number $k$ plays a major role in the asymptotics of $\eta_c$ when $c\rightarrow c_s$. Recall indeed that the function $\eta_c$ solves the equation 
\begin{equation*}\label{equation -eta'^2=nu(eta)}
    -(\partial_x\eta_c)^2=\mathcal{N}_c(\eta_c),
\end{equation*}
where \begin{equation}\label{expression de N_c}
    \mathcal{N}_c(x)=c^2 x^2 -4(1-x)F(1-x),
\end{equation}
so that
\begin{equation}\label{equation -eta'' = g(eta)}
    -\partial_x^2\eta_c=\dfrac{1}{2}\mathcal{N}_c'(\eta_c).
\end{equation}

Since
\begin{equation}\label{développement asymptotique de n_c(xi) en 0}
    \mathcal{N}_c(\xi)\underset{\xi\rightarrow 0}{=}-\xi^2\big(\nu_c^2+\dfrac{k}{3} \xi + \mathcal{O}(\xi^2)\big),
\end{equation}

and 
\begin{equation}\label{développement asymptotique de n_c'(xi) en 0}
    \mathcal{N}'_c(\xi)\underset{\xi\rightarrow 0}{=}-\xi\big(2\nu_c^2+ k \xi + \mathcal{O}(\xi^2)\big),
\end{equation}
we deduce the following control on $\eta_c$.
\begin{prop}\label{proposition: borne inférieure eta_c(x)}
We have \begin{equation*}\label{convergence uniforme de eta_c vers 0 quand c c_s}
    \eta_c\underset{c\rightarrow c_s}{\longrightarrow}0\quad\text{in }L^\ii(\R).
\end{equation*}

Moreover, there exists $c_1\in (c_0,c_s)$ such that for $c\in (c_1,c_s)$ and any $x\in\R$,
\begin{equation*}\label{équivalent quand x tend vers +ii de eta_c'/eta_c}
    \left|\eta_c'(x)\right|\underset{c\rightarrow c_s}{\sim} \eta_c(x)\nu_c,
\end{equation*}
and there exists $m_c>0$ depending continuously on $c$ and tending to $0$ as $c\rightarrow c_s$ such that for any $x\in\R$,
\begin{equation*}\label{borne inférieure de eta_c}
    |\eta_c(x)|\geq m_c e^{-\nu_c |x|} .
\end{equation*}
\end{prop}

\begin{proof}
The proof can be derived from what is shown in Section~2 in~\cite{Bert2}. However, no control from below was provided for the function $\eta_c$ that is only known to be positive by~\eqref{eta_c(x) vit dans (0,1)}. Therefore, we need to come back to~\cite{Chiron8}, from which we have 
\begin{equation*}
    |\eta_c(x)|\underset{|x|\rightarrow + \ii}{\sim} M_c e^{-\nu_c|x|},
\end{equation*}
where
\begin{equation}\label{expression de M_c}
    M_c=\xi_c e^{\int^{\xi_c}_0\frac{\mathcal{N}_{c_s}(\xi)}{\sqrt{-\xi^2\mathcal{N}_c(\xi)}\big(\sqrt{-\mathcal{N}_c(\xi)}+\sqrt{\nu_c^2\xi^2}\big)}d\xi }. 
\end{equation}

Taking the limit $c\rightarrow c_s$ in the previous integral\footnote{The computations can be properly justified on the basis of those made in~\cite{Bert2}.} yields to
\begin{equation*}
    M_c\underset{c\rightarrow c_s}{\sim}\xi_c e^{\frac{\nu_c^2}{1+\nu_c}}.
\end{equation*}

By~\eqref{équivalent de xi_c}, this provides $c_1\in (c_0,c_s)$ such that $M_c$ is positive for $c\in (c_1,c_s)$. Thus, by~\eqref{eta_c(x) vit dans (0,1)}, we deduce the existence of a constant $m_c$ continuously depending on $c$ such that $\eta_c(x)\geq m_c e^{-\nu_c |x|}$.

\end{proof}

Going back to the quadratic form in the proof of Lemma~\ref{lemme: equation sur la variable duale epsilon^*} we recall the notations $a_c=\frac{\eta_c'}{\eta_c}$ and $\lambda_c=\frac{\eta_c'}{(1-\eta_c)^2}$. We also develop the expressions of

\begin{equation*}
    q_{1,c}=\dfrac{2\big((\eta_c')^2-\eta_c (1-\eta_c)\eta_c''\big)}{\eta_c^2},
\end{equation*}

and of the next quantity that will be shown to be well-defined in the next proposition 
\begin{equation*}
    \widetilde{q_{1,c}}=-\dfrac{3a_c'}{2(1-\eta_c)}-\dfrac{a_c\lambda_c}{2}-\dfrac{c^2(\eta_c')^2}{(1-\eta_c) q_{1,c}}.
\end{equation*}

\begin{prop}\label{proposition: signe de q_1,c et q_2,c}
    There exists $c_1\in (c_0,c_s)$ such that we have for any $(c,x)\in (c_1,c_s)\times \R$, 
    \begin{equation*}
        \dfrac{q_{1,c}(x) }{\eta_c(x)}>0\quad\text{and}\inf_{(c,x)\in (c_1,c_s)\times\R}\quad \dfrac{\widetilde{q_{1,c}}(x)}{\eta_c(x)}>0.
    \end{equation*}
\end{prop}

\begin{proof}
We make a Taylor expansion of $\mathcal{N}_c(x)$ in~\eqref{expression de N_c} and we eventually compute
\begin{equation*}
    q_{1,c}(x)=\eta_c(x)\left(2\nu_c^2 + \int_0^1 (1-t)\big(t-\eta_c(x)\big)\mathcal{N}_c'''\big(t\eta_c(x)\big)dt\right).
\end{equation*}

Since $\mathcal{N}'''_c(0)=-2k >0$, we use the uniform convergence~\eqref{convergence uniforme de eta_c vers 0 quand c c_s} to infer that for $c$ close enough to $c_s$ we have 

\begin{equation}\label{développement asymptotique de q_1,c en grand o de eta_c^2}
    q_{1,c}= \left(2\nu_c^2-\dfrac{k}{3}\right)\eta_c+\grandOde{\eta_c^2}.
\end{equation}

Since $\nu_c$ tends to $0$ as $c$ tends to $c_s$, we infer that there exists $c_1$ close enough to $c_s$ such that for any $x\in\R$,
\begin{equation*}
    q_{1,c}(x)\geq \dfrac{-k\eta_c(x)}{6}.
\end{equation*}

We conclude that $\frac{q_{1,c}(x) }{\eta_c(x)} >0$ for any $(x,c)\in\R\times (c_1,c_s)$, using~\eqref{eta_c(x) vit dans (0,1)}.

Now, we handle $\widetilde{q_{1,c}}$. 
From~\eqref{développement asymptotique de n_c(xi) en 0} and~\eqref{développement asymptotique de n_c'(xi) en 0}, we derive \begin{equation*}
    a_c'=\dfrac{k}{6}\eta_c + \grandOde{\eta_c^2},
\end{equation*}

hence \begin{equation*}
    -\dfrac{3a_c'}{2(1-\eta_c)}=-\dfrac{k}{4}\eta_c + \grandOde{\eta_c^2}.
\end{equation*}

On the other hand, from~\eqref{développement asymptotique de n_c(xi) en 0} again, we compute
\begin{equation*}
    -\dfrac{a_c\lambda_c}{2}=-\dfrac{\nu_c^2}{2}\eta_c + \grandOde{\eta_c^2},
\end{equation*}
and using~\eqref{développement asymptotique de q_1,c en grand o de eta_c^2},
\begin{equation*}
    -\dfrac{c^2(\eta_c')^2}{(1-\eta_c) q_{1,c}}=-\dfrac{c^2\nu_c^2}{2\nu_c^2-\frac{k}{3}}\eta_c+\grandOde{\eta_c^2}.
\end{equation*}

Using all three previous estimates, we derive
\begin{equation}\label{développement asymptotique de widetilde q_1,c}
    \widetilde{q_{1,c}}=-\left(\dfrac{k}{4}+\dfrac{\nu_c^2}{2}+\dfrac{c^2\nu_c^2}{2\nu_c^2-\frac{k}{3}}\right)\eta_c +\grandOde{\eta_c^2}.
\end{equation}

In sum, we conclude that, up to taking $c_1$ even closer to $c_s$, we have for any $(x,c)\in\R\times (c_1,c_s)$, \begin{equation*}
    \dfrac{\widetilde{q_{1,c}}(x)}{\eta_c(x)}\geq -\dfrac{k}{8}>0.
\end{equation*}

\end{proof}

~~\\

\begin{merci}
The author thanks Philippe Gravejat for his well-meaning eye and help in some tedious computations made in this article.
\end{merci}

\bibliographystyle{plain}
\bibliography{bibliography}

\end{document}